\documentclass[twoside, 14pt]{article}
\usepackage{amsmath,amssymb,amsthm,mathrsfs}
\usepackage{ulem}
\usepackage[margin=2.0cm]{geometry} % showframe,
\usepackage{lipsum}
\usepackage{titlesec,hyperref}
\usepackage{fancyhdr}
\usepackage[numbers,sort&compress]{natbib}
\usepackage{color}
\usepackage[titletoc]{appendix}

\fancypagestyle{plain}
{
	\fancyhf{}
	
}

\titleformat{\subsection}{\it}{\thesubsection.\enspace}{1.5pt}{}
\titleformat{\subsubsection}{\it}{\thesubsubsection.\enspace}{1.5pt}{}

\newtheorem{theorem}{Theorem}[section]

\newtheorem{lemma}{Lemma}[section]
\newtheorem{proposition}{Proposition}[section]
\newtheorem{remark}{Remark}[section]
\newtheorem{corollary}{Corollary}[section]

\numberwithin{equation}{section}

\allowdisplaybreaks

\begin{document}
	\title{Global existence and decay estimates of solutions
		for the compressible Prandtl type equations in Sobolev spaces \hspace{-4mm}}
	\author{Yuming Qin$^{\ast}$$^{a,b}$ \quad Junchen Liu$^{b}$ \\[10pt]
		\small {$^a$ Institute for Nonlinear Sciences, Donghua University, 201620, Shanghai, P. R. China}\\
		\small {$^b$ Department of Mathematics, School of Mathematics and Statistics, Donghua University,
	201620, Shanghai, P. R. China}
		\\[5pt]
	}

	\footnotetext{Email: \it $^{\ast}$   yuming@dhu.edu.cn, yuming\_qin@hotmail.com}
	\date{}
	
	\maketitle

	\begin{abstract}
This paper addresses the global existence and large-time decay of solutions to the two-dimensional compressible Prandtl equations in Sobolev spaces. We study small perturbations around a monotone shear flow and establish global weighted Sobolev estimates together with precise time-decay rates for the solutions. The main difficulties are that the normal velocity causes a loss of one tangential derivative, the variable density gives rise to nonstandard boundary and shear-flow coupling terms, and the global-in-time argument requires sufficiently strong decay to close the time-weighted energy estimates. To overcome these difficulties, we combine the cancellation method introduced by Masmoudi and Wong~\cite{MW2015} with the good-quantity approach developed by Paicu and Zhang~\cite{PZ2021}, together with weighted shear-flow estimates and a careful treatment of the natural boundary structure. Under suitable assumptions on the density perturbation, the shear flow, and the function $f(t)$ describing the outer tangential flow, we prove the global-in-time well-posedness of solutions in weighted Sobolev spaces. In this sense, the present result may be viewed as a Sobolev-space counterpart of the recent analytic-framework result obtained by Chen, Huang and Li 
(J. Funct. Anal. 286 (2024), 110322)~\cite{CHL2024}
for the two-dimensional compressible Prandtl-type equations, while requiring a weaker time-decay assumption on the density perturbation.
		
		\vspace*{5pt}
		\noindent{\it {\rm Keywords}}:  Prandtl equations; compressible fluids; Sobolev spaces; 
%		monotone shear flow; 
		energy method.
		
	\vspace*{5pt}
		\noindent{\it {\rm  Mathematics Subject Classification:}}\ {\rm 35Q30; 35Q35; 76N06; 76N20}
	\end{abstract}

	\tableofcontents

\section{Introduction}

The study of boundary layers is one of the most challenging and fundamental topics in mathematical fluid mechanics. Originating from Prandtl's seminal work \cite{P1904} in 1904, boundary layer theory was introduced to resolve the discrepancy between the boundary conditions for ideal inviscid fluids and viscous fluids with small viscosity. Roughly speaking, boundary layers have amplitude $O(1)$ and thickness $O(\sqrt{\nu})$, and thus form a critical transition region in which the flow changes from being governed by the Euler equations in the interior to the Navier--Stokes equations with vanishing viscosity $\nu>0$ near the boundary. At the formal level, this transition is described by the Prandtl equations.

As a basic model for boundary layer dynamics, the Prandtl equations can be formally derived from the incompressible Navier--Stokes equations under the no-slip boundary condition in the inviscid limit. They describe the leading-order behavior of viscous flow near a physical boundary, where the velocity varies rapidly because of the mismatch between viscous and inviscid boundary conditions. Beyond their physical significance, the Prandtl equations also play a central role in the rigorous justification of inviscid limits. Over the years, many powerful analytic tools have been developed to study their well-posedness and structural properties. More recently, some of these ideas have been extended to compressible settings, where density effects make the boundary layer dynamics substantially more involved. In this case, the boundary layer behavior is governed by the following compressible Prandtl-type system:
\begin{equation}\label{prandtl}
	\begin{cases}
		\rho\left(\partial_{t}U+U\partial_{x}U+V\partial_{y}U\right)-\partial_{y}^2U+\partial_{x}p=0,\quad  (t, x, y)\in \mathbb{R}_+ \times \mathbb{R}^2_+, \\
		\partial_{t}\rho+\partial_x \left(\rho U\right) +\partial_y\left(\rho V\right)  =0, \quad  (t, x, y)\in \mathbb{R}_+ \times \mathbb{R}^2_+,\\
		U|_{y=0} = V|_{y=0} =0 , \ \lim\limits_{y\to+\infty} U(t, x, y) = \bar{U}(t,x), \quad  (t, x)\in \mathbb{R}_+ \times \mathbb{R},\\
		U|_{t=0} =U_0 (x,y), \quad (x, y)\in \mathbb{R}^2_+\, ,
	\end{cases}
\end{equation}
where $\mathbb{R}^2_{+}=\mathbb{R}\times \mathbb{R}_+$. Here and throughout, $U$ and $V$ denote the tangential and normal velocities of the boundary layer flow, respectively, while $(\bar U(t,x),\rho(t,x),p(t,x))$ are given functions of $(t,x)$ describing the tangential velocity, density, and pressure of the outer Euler flow. These quantities naturally satisfy the Bernoulli law
\begin{align}
	\rho(\partial_{t}\bar{U}+\bar{U}\partial_{x}\bar{U})+\partial_{x}p=0, \quad (t, x)\in \mathbb{R}_{+} \times \mathbb{R}.
	\label{bernoulli}
\end{align}
The pressure is given by a smooth function $p=p(\rho)$ satisfying $p'(\rho)>0$ for all $\rho>0$.

There is a vast literature on the theoretical, numerical, and experimental study of Prandtl-type boundary layers. Before discussing the compressible case in detail, we first review some representative results for the incompressible Prandtl equations, for which a rather complete theoretical framework has been developed. A central difficulty in the well-posedness theory is that the normal velocity, given by
$
V=-\partial_y^{-1}\partial_xU,
$
causes a loss of one tangential derivative. As a consequence, the standard energy method does not apply directly.

There are essentially two main approaches to overcoming this difficulty. One is to work in H\"older or Sobolev spaces under a monotonicity assumption on the tangential velocity. The other is to study the Prandtl equations in analytic or, more generally, Gevrey-type function spaces, where the loss of tangential derivatives can be compensated by additional regularity in the tangential variable.

We first recall some results in the monotone setting. Early studies of the Prandtl equations mainly focused on the two-dimensional case, for which substantial progress has been made. Oleinik \cite{O1966} first established the existence and uniqueness of classical solutions to the two-dimensional Prandtl equations in H\"older spaces under Oleinik's monotonicity condition. More precisely, solutions exist locally in time on $[0,T]\times [0,L]\times \mathbb{R}_+$ for arbitrary $L$, and globally in time when $L$ is sufficiently small. This result, together with a systematic introduction to boundary layer theory, was presented in the classical monograph \cite{OS1999} by Oleinik and Samokhin. By using the Crocco transformation developed in \cite{O1966,OS1999}, Xin and Zhang \cite{XZ2004} obtained the global existence of BV weak solutions for arbitrary length $L$ under the additional favorable pressure condition $\partial_x P\le 0$.
More recently, Xin, Zhang, and Zhao \cite{XZ2024} significantly refined this result by proving the uniqueness, continuous dependence, and global regularity of such weak solutions, and hence the global existence of smooth solutions.

Using a direct energy method, Alexandre, Wang, Xu, and Yang \cite{AWXY2015} established local-in-time well-posedness in weighted Sobolev spaces via a Nash--Moser iteration, without relying on the Crocco transformation. In contrast, Masmoudi and Wong \cite{MW2015} discovered a key cancellation structure in the convection terms and, under the monotonicity assumption, proved local well-posedness by an energy method that overcomes the loss of tangential derivatives. Building on \cite{MW2015,AWXY2015}, Xu and Zhang \cite{XZ2017} obtained a long-time well-posedness result for the two-dimensional Prandtl equation in Sobolev spaces with polynomial weights on the half-plane, and proved that the lifespan can be made arbitrarily large when the initial data are a small perturbation of a monotone shear profile. More recently, Liu and Zhang \cite{LZ2024} proved the global well-posedness of the Prandtl equations in Sobolev spaces by deriving sufficiently fast large-time decay estimates for suitable weighted energy norms of the tangential velocity. We also refer to \cite{QL2025,CWZ2018,KMVW2014,QW2025,QWL2025} for other results on the incompressible Prandtl equations under monotonicity assumptions.

When Oleinik's monotonicity condition fails, however, blow-up, instability, and ill-posedness phenomena may occur; see, for instance, \cite{G2000,HH2003,GVD2010,LY2017,GVN2010,GN2011,EE1997}. These results indicate that the monotonicity assumption is essential for well-posedness in Sobolev-type spaces, unless one instead works in analytic or Gevrey frameworks. We now briefly review some representative results in these settings.

In the analytic framework, Sammartino and Caflisch first proved the local well-posedness of the two-dimensional Prandtl equations, and their argument also applies to the three-dimensional case, for initial data analytic in both the $x$ and $y$ variables. Subsequently, Lombardo, Cannone, and Sammartino \cite{LCS2003} improved this result by removing the analyticity assumption in the normal variable $y$, making essential use of the regularizing effect of the diffusion operator $\partial_t-\partial_y^2$. For further local well-posedness results under analyticity assumptions, we refer to \cite{KV2013,CLS2001,CLS2013,KMVW2014} and the references therein.

As for the long-time existence of small analytic solutions, Zhang and Zhang \cite{ZZ2016} first proved, by using Littlewood--Paley theory in the tangential variables, that the Prandtl system in $\mathbb{R}_{+}\times \mathbb{R}^{d-1}$, with $d=2,3$, has a unique solution with lifespan $T_{\epsilon} \geq \epsilon^{-\frac{4}{3}}$. Based on the observation that the analytic norm of a suitable good unknown decays faster than that of the original unknown, Ignatova and Vicol \cite{IV2016} proved the existence of an almost global-in-time solution to the two-dimensional Prandtl equation, with lifespan larger than $(\epsilon^{-1}/\log(\epsilon^{-1}))$. By refining this good quantity, Paicu and Zhang \cite{PZ2021} established the global well-posedness of the Prandtl equation in analytic spaces under the assumption that the outer flow is sufficiently small and decays in time.

When the analyticity requirement is relaxed to Gevrey regularity, several important advances have also been made. G\'erard-Varet and Masmoudi \cite{GVMG2015} first proved local well-posedness for the two-dimensional Prandtl equations for initial data in the Gevrey class $\frac{7}{4}$, without analyticity or monotonicity assumptions. This result was later extended in \cite{LY2020} to Gevrey indices $\sigma\in[\frac32,2]$ for data that are small perturbations of a shear flow with a single non-degenerate critical point, where $\sigma=2$ is optimal in view of the ill-posedness results in \cite{GVD2010}. The well-posedness of the linearized Prandtl equation around a non-monotone shear flow was obtained by Chen, Wang, and Zhang \cite{CWZ22018} in the Gevrey class $2-\theta$ for any $\theta>0$. Subsequently, Dietert and G\'erard-Varet \cite{DGV2019} established local well-posedness for initial data with Gevrey regularity of index $2$ in the tangential variable $x$ and Sobolev regularity in the normal variable $y$, thereby improving \cite{GVMG2015} by removing the assumption of a single non-degenerate critical point. Inspired by these local existence results, especially \cite{DGV2019}, Wang, Wang, and Zhang \cite{WWZ2024} proved the global existence of small Gevrey-$2$ solutions, extending the analytic result of \cite{PZ2021} to the optimal Gevrey class.

The results mentioned above mainly concern incompressible boundary layer equations, for which a rather rich mathematical theory has been developed over the past decades. In contrast, the study of compressible boundary layers remains relatively limited. Compared with the incompressible case, the density in compressible boundary layers plays a crucial role, since it not only reflects compressibility effects and possible expansion of the fluid, but also makes the overall structure of the system much more complicated. Moreover, when heat conduction is taken into account, rapid heat transfer near the boundary leads to the so-called thermal boundary layer; see, for example, \cite{K1938,LW2014,LWY2021}. The interaction between viscous and thermal layers then gives rise to a considerably more intricate boundary layer structure. Therefore, understanding the formation and stability of compressible boundary layers remains a challenging problem of both theoretical and practical significance, especially in the study of high-speed flows and shock--boundary layer interactions.

The boundary layer equations for compressible fluids were first derived by Wang and Williams \cite{WW2012}, who investigated the behavior and stability of boundary layers in two-dimensional isentropic compressible flows under Navier boundary conditions, where the slip length depends on the viscosity, and the no-slip case appears as a limiting regime. Later, Wang, Xie, and Yang \cite{WXY2015} established the local existence of solutions to the compressible Prandtl equations with monotone data; see also \cite{FRY2021,GGW2016}. By using a viscous splitting method, the global existence of weak solutions was proved in \cite{DG2017}. More recently, Chen, Huang, and Li \cite{CHL2024} introduced a new good quantity and established the global existence and large-time decay of strong solutions to the two-dimensional compressible Prandtl equations with small initial data that are analytic in the tangential variable. Very recently, the local well-posedness of the three-dimensional compressible boundary layer system has been established in \cite{CRY2025}, under the assumption that the initial data are real analytic in the tangential variables and have Sobolev regularity in the normal variable.

It should be emphasized that the existing results are either local under special assumptions such as monotonicity or analyticity, or global only in analytic-type settings. In particular, the global existence problem for the two-dimensional compressible Prandtl equation in Sobolev spaces remains open. This is precisely the problem considered in the present paper. For simplicity, throughout this paper we take $\bar{U}(t,x)$ in system \eqref{prandtl} to be $\varepsilon f(t)$ with $f(0)=0$, which together with \eqref{bernoulli} implies
$$
\partial_x p=-f'(t), \qquad \partial_t\rho=-\varepsilon f(t)\partial_x\rho.
$$
Motivated by the study in \cite{CHL2024} for the compressible Prandtl equations, we introduce a cut-off function $\chi\in C^\infty[0,\infty)$ satisfying
\begin{equation}\label{chi}
	\chi(y) = \begin{cases} 
		1 & \text{if } y \geq 2, \\
		0 & \text{if } y \leq 1, 
	\end{cases}
	\quad \text{and} \quad
	\int_{0}^{\infty} (1 - \chi(y)) \, dy = 0. 
\end{equation}
The second condition in \eqref{chi} may appear somewhat unusual, but it is crucial for the proof of a key Lemma \ref{uG} below. Let us briefly explain why such a function exists. Choose any $y_1,y_2\in[1,2]$ with $y_1<y_2$. One can construct $\chi$ so that $\chi(y)\ge 1$ on $[y_1,y_2]$ and
$$
\int_1^2 (1-\chi(y))\,dy=-1 ,
$$
which, together with the conditions $\chi(y)=0$ for $y\le 1$ and $\chi(y)=1$ for $y\ge 2$,  yields
$$
\int_0^\infty (1-\chi(y))\,dy=0.
$$
Hence, a smooth cut-off function satisfying \eqref{chi} does indeed exist.

From now on, we denote
$$
W \buildrel\hbox{\footnotesize def}\over = U-\varepsilon f(t)\chi(y).
$$
Then $W$ solves
\begin{equation}\label{wprandtl}
\left\{
\begin{array}{ll}
\partial_t W + (W + \varepsilon f(t)\chi(y)) \partial_x W + V \partial_y (W + \varepsilon f(t)\chi(y)) \\
\quad -(a+1)\partial_y^2 W - a\varepsilon f(t)\chi''(y) = \varepsilon m, \quad (t,x,y) \in \mathbb{R}_+ \times \mathbb{R}_+^2, \\
\partial_x (\rho W) + \partial_y (\rho V) = \varepsilon f(t)\partial_x \rho (1 - \chi(y)), \quad (t,x,y) \in \mathbb{R}_+ \times \mathbb{R}_+^2, \\
W|_{y=0} = V|_{y=0} = 0 \quad \text{and} \quad \lim\limits_{y \to +\infty} W(t,x,y) = 0, \quad (t,x) \in \mathbb{R}_+ \times \mathbb{R}, \\
W|_{t=0} = U_0, \quad (x,y) \in \mathbb{R}_+^2,
\end{array}
\right.
\end{equation}
where
$a \buildrel\hbox{\footnotesize def}\over =\frac{1}{\rho}-1$ and 
  $m(t,y)\buildrel\hbox{\footnotesize def}\over =(1-\chi(y))f'(t)+ f(t)\chi''(y).$

Let $u^s(t, y)$ satisfy the equation
\begin{equation}\label{shear flow}
	\left\{
	\begin{array}{ll}
		\partial_t u^s - \partial_y^2 u^s = \varepsilon m, \quad (t,y) \in \mathbb{R}_+ \times \mathbb{R}_+, \\
		u^s|_{y=0} = 0 \quad \text{and} \quad \lim\limits_{y \to +\infty} u^s(t,y) = 0, \quad t \in \mathbb{R}_+, \\
		u^s|_{t=0} = u_0^{s}, \quad u_0^s(0)=0, \quad y \in \mathbb{R}_+.
	\end{array}
	\right.
\end{equation}

Consider the solution of equation \eqref{wprandtl} as a perturbation of the shear flow $(u^s, 0)$, and write $(u, v) \buildrel\hbox{\footnotesize def}\over=(W-u^s, V)$.
Thus, from \eqref{wprandtl} and \eqref{shear flow}, we derive
\begin{equation}\label{uprandtl}
	\left\{
	\begin{array}{ll}
		\partial_t u + (u + u^s + \varepsilon f(t)\chi(y)) \partial_x u + v \partial_y \left( u + u^s + \varepsilon f(t)\chi(y) \right) \\
		\quad-(a+1)\partial_y^2 u - a \partial_y^2  u^s - a\varepsilon f(t)\chi''(y) = 0, \quad (t,x,y) \in \mathbb{R}_+ \times \mathbb{R}_+^2, \\
		\partial_x (\rho u) + \partial_y (\rho v) = \varepsilon f(t)\partial_x \rho (1 - \chi(y)) - \partial_x \rho u^s, \quad (t,x,y) \in \mathbb{R}_+ \times \mathbb{R}_+^2, \\
		u|_{y=0} = v|_{y=0} = 0 \quad \text{and} \quad \lim\limits_{y \to +\infty} u(t,x,y) = 0, \quad (t,x) \in \mathbb{R}_+ \times \mathbb{R}, \\
		u|_{t=0} \stackrel{\text{def}}{=} u_0 = U_0-u^s_0, \quad (x,y) \in \mathbb{R}_+^2.
	\end{array}
	\right.
\end{equation}
Meanwhile, since $\partial_x (\rho u) + \partial_y (\rho v) = \varepsilon f(t)\partial_x \rho (1 - \chi(y)) - \partial_x \rho u^s$, there exists a potential function $\varphi$ on $\mathbb{R}_+^2$ such that
\begin{align}\label{potential function}
u \stackrel{\text{def}}{=} \partial_y \varphi \quad \text{and} \quad v \stackrel{\text{def}}{=} -\partial_x \varphi - \frac{\partial_x a}{a+1} \int_y^\infty \left[ u + u^s - \varepsilon f(t)(1 - \chi(y')) \right] dy'.
\end{align}
In fact, observe that $\varphi(t, x, y)=-\int_{y}^\infty u(t, x, y')dy'$, by integrating \eqref{uprandtl} over the interval $[y, +\infty)$, one obtains

\begin{equation}\label{varphi uprandtl}
	\left\{
\begin{array}{ll}
	\partial_t \varphi + (u + u^s + \varepsilon f(t)\chi(y)) \partial_x \varphi + 2 \int_y^\infty \left( \partial_{y'} \left( u + u^s + \varepsilon f(t)\chi(y') \right) \partial_x \varphi \right) dy' \\
	\quad -\partial_y^2 \varphi = F_1 + F_2, \quad (t,x,y) \in \mathbb{R}_+ \times \mathbb{R}_+^2, \\
	\partial_y \varphi|_{y=0} = 0 \quad \text{and} \quad \lim\limits_{y \to +\infty} \varphi(t,x,y) = 0, \quad (t,x) \in \mathbb{R}_+ \times \mathbb{R}, \\
	\varphi|_{t=0} \stackrel{\text{def}}{=} \varphi_0, \quad (x,y) \in \mathbb{R}_+^2,
\end{array}
	\right.
\end{equation}
where
\begin{align}
\begin{split}
	F_1 &\stackrel{\text{def}}{=} a\partial_y u^s + a\varepsilon f(t)\chi'(y), \\
F_2 &\stackrel{\text{def}}{=} -\frac{\partial_x a}{a+1} \int_y^\infty \left[ u + u^s - \varepsilon f(t)(1 - \chi(y')) \right] \left[ u + u^s + \varepsilon f(t)\chi(y') \right] dy' \\
&\quad + \frac{\partial_x a}{a+1} \int_y^\infty \left[ u + u^s - \varepsilon f(t)(1 - \chi(y')) \right] dy' \left[ u + u^s + \varepsilon f(t)\chi(y') \right].
\end{split}
\label{f1f2}
\end{align}

On the one hand, it follows from \eqref{potential function} that 
\begin{align*}
 v \stackrel{\text{def}}{=} -\int_{y}^\infty \partial_x u(t, x, y')dy' - \frac{\partial_x a}{a+1} \int_y^\infty \left[ u (t, x, y') + u^s (t, y')- \varepsilon f(t)(1 - \chi(y')) \right] dy',
\end{align*}
which yields the classical loss of one tangential derivative in the standard energy method. To overcome this difficulty, we introduce the good unknown
\begin{align}\label{g}
	g_m=\partial_x^m w- \frac{\partial_y(w + w^s + \varepsilon f(t)\chi'(y))}{w + w^s + \varepsilon f(t)\chi'(y)}\partial_{x}^{m}u
	=
	(w + w^s + \varepsilon f(t)\chi'(y))\partial_y\left( \frac{\partial_x^m u}{w + w^s + \varepsilon f(t)\chi'(y)}\right) ,
\end{align}
inspired by  \cite{MW2015}. Here, the  vorticity $w=\partial_y u$  satisfies the following vorticity system:
\begin{equation}\label{vprandtl}
	\left\{
	\begin{array}{ll}
		\partial_t w + (u + u^s + \varepsilon f(t)\chi(y)) \partial_x w + v \partial_y \left( w + w^s + \varepsilon f(t)\chi'(y) \right) \\
		\quad -(a+1)\partial_y^2 w - a \partial_y^2  w^s - a\varepsilon f(t)\chi'''(y) = F_3, \quad (t,x,y) \in \mathbb{R}_+ \times \mathbb{R}_+^2,  \\
\bigl((a+1)\partial_y w +a\partial_y w^s+a\varepsilon f(t)\chi''(y)\bigr)|_{y=0}  = 0 , \quad (t,x) \in \mathbb{R}_+ \times \mathbb{R}, \\
		w|_{t=0} \stackrel{\text{def}}{=} w_0, \quad (x,y) \in \mathbb{R}_+^2,
	\end{array}
	\right.
\end{equation}
where 
\begin{align*}
F_3\stackrel{\text{def}}{=}\frac{\partial_x a}{a+1} \left[ u + u^s - \varepsilon f(t)(1 - \chi(y)) \right] (w + w^s + \varepsilon f(t)\chi'(y)).
\end{align*}

On the other hand, in order to derive a faster decay rate for the solution $u$ of system \eqref{uprandtl}, which is needed to control the lower-order derivatives in Sobolev spaces, we recall from \cite{CHL2024} the good quantity
\begin{align}\label{G}
	G \overset{\text{def}}{=} u + \frac{y}{2\langle t\rangle} \varphi ,
	\qquad
	\partial_y G = \partial_y u + \frac{y}{2\langle t\rangle} u + \frac{\varphi}{2\langle t\rangle},
\end{align}
which was introduced in \cite{IV2016,PZ2021}. Here and below,
$$
\langle t\rangle \overset{\text{def}}{=} 1+t.
$$
Then $G$ solves
\begin{equation}\label{G uprandtl}
	\left\{
	\begin{array}{ll}
		\partial_t G - (a+1)\partial_y^2 G + \langle t \rangle^{-1} G + (u + u^s + \varepsilon f(t)\chi(y)) \partial_x G + v \partial_y G + v \partial_y \left( u^s + \varepsilon f(t)\chi(y) \right) \\ 
		\quad 
		-\frac{1}{2} \langle t \rangle^{-1} v \partial_y (y\varphi) 
		- a \partial_y^2 u^s 
		+ \frac{y}{\langle t \rangle} \int_y^\infty 
		\left( \partial_{y'} (u + u^s + \varepsilon f(t)\chi(y')) \partial_x \varphi \right) dy' 
		\\ 
		\quad -a \varepsilon f(t)\chi''(y) - \frac{y}{2\langle t \rangle}(F_1 + F_2) = 0, \quad (t, x, y) \in \mathbb{R}_+ \times \mathbb{R}_+^2, \\ 
		G\big|_{y=0} = 0 \quad \text{and} \quad \lim\limits_{y \to +\infty} G(t, x, y) = 0, \quad (t, x) \in \mathbb{R}_+ \times \mathbb{R}, \\ 
		G\big|_{t=0} \stackrel{\text{def}}{=} G_0 = u_0 + \frac{y}{2}\varphi_0, \quad (x, y) \in \mathbb{R}_+^2,
	\end{array}
	\right.
\end{equation}
where $F_1$ and $F_2$ are defined in \eqref{f1f2}.

\textbf{Notations.}
Before stating the main result, we introduce some notation that will be used throughout the paper.
We write $\|\cdot\|_{L^2}$ and $\langle \cdot , \cdot \rangle_{L^2}$ for the norm and inner product in
$
L^2=L^2(\mathbb{R}_+^2)=L^2(\mathbb{R}\times \mathbb{R}_+).
$
More generally, for $1\le p,q\le \infty$, we denote by
$
L_x^pL_y^q=L^p(\mathbb{R};L^q(\mathbb{R}_+))
$
the usual mixed Lebesgue space.
In particular, we write
$
L_x^p=L^p(\mathbb{R}),
$
and
$
L_y^q=L^q(\mathbb{R}_+).
$
For $s\in\mathbb{N}$, we denote by $H_x^s=H^s(\mathbb{R})$ the usual Sobolev space with respect to the tangential variable $x$.
We also use the anisotropic Sobolev space $H^{m,k}(\mathbb{R}_+^2)$.
For $m,k\in\mathbb{N}$, it is defined by
$$
\| f \|_{H^{m,k}}^2 \overset{\text{def}}{=} \sum_{\substack{0 \leq j \leq m \\ 0 \leq \ell \leq k}} \left\| \partial_x^j \partial_y^\ell f \right\|_{L^2}^2.
$$
Finally, for two operators $A$ and $B$, we write $[A,B]$ for their commutator.

Throughout this paper, we  use the time-dependent Gaussian weight
\begin{align}
	\Psi=\Psi(t,y)\overset{\text{def}}{=}\frac{y^2}{8\langle t\rangle }.
	\label{Psi}
\end{align}

We are now in a position to state our main result.
\begin{theorem}\label{t1}
	Fix $0<\kappa_0<\frac14$, and let $\gamma_0\in(0,1)$, $\kappa\in(\kappa_0,\frac12)$, and $c(\kappa_0)$ be a positive constant such that
$
		0<c(\kappa_0)<3\kappa_0.
$
	Suppose that the following conditions are fulfilled.
	
	\smallskip
	
	\noindent
	$(i)$ The function $f\in H^1(\mathbb{R}_+)$ satisfies
	\begin{align}
		\begin{split}
		C_f &\overset{\mathrm{def}}{=}
\sup_{t\in[0,+\infty)} \left( |f(t)|+|f'(t)|\right) 
+ \left(\int_0^\infty \langle t\rangle^{\frac92-2\kappa}f^2(t)\,dt\right)^{\frac12}
<+\infty.
		\end{split}
		\label{Cf}
	\end{align}
	
	\smallskip
	
	\noindent
	$(ii)$ The shear flow $u^s$ determined by \eqref{shear flow} satisfies
	\begin{align}
		\left( \partial_y u^s(t,y)+\varepsilon f(t)\chi'(y)\right)
		\geq c_0 \langle t \rangle^{-2+\kappa} e^{-\frac{\gamma_0 y^2}{8 \langle t \rangle}},
		\quad
\left\| e^{\frac{\gamma_0 y^2}{8}} \partial_y^2 u_0^s \right\|_{L_y^2} \leq M
		,
		\label{CM}
	\end{align}
	and
	\begin{align}
		\left|  \partial_y^2 u^s (t,y) \right| <  C_M \left\langle t\right\rangle^{-2+\kappa}  e^{-\frac{(2-\gamma) y^2}{8\left\langle t\right\rangle }}, \qquad \left( \int_{0}^{\infty} \left\|\left\langle t\right\rangle^{\frac{9}{4}-\kappa}  e^{\frac{\gamma y^2}{8}} \partial_y^3 u^s (t)\right\| _{L^2_y}^2 dt\right)^{\frac{1}{2
		}}
		\leq C_M.
		\label{CM1}
	\end{align}
	
	\smallskip
	
	\noindent
	$(iii)$ The density $\rho=\rho(t,x)$ satisfies
	\begin{equation}\label{rho-upper-new}
		1\le \rho(t,x)\le \bar\rho<
		\frac{1-c(\kappa_0)-\kappa_0}{1-4\kappa_0},
	\end{equation}
and the function $a$ satisfies
	\begin{align}
		C_a \overset{\mathrm{def}}{=}
		\sup_{t\in[0,+\infty)}
		\left\|\langle t\rangle^{\frac12+\kappa_0} a(t)\right\|_{H^{m+1}_x}
		<+\infty,
		\label{Ca}
	\end{align}
	together with
	\begin{align}
C_{a,\infty} \overset{\mathrm{def}}{=}	\max_{1\le i\le m}
		\Bigl\{
		C_m^i \|\partial_x^i a\|_{L_x^\infty}
		\Bigr\}
		\leq
		\frac{1-c(\kappa_0)-\kappa_0-(1-4\kappa_0)\bar\rho}{6(m-1)\bar\rho},
		\label{Cainfty}
	\end{align}
where 	$a \overset{\mathrm{def}}{=} \frac1\rho-1$.
	\smallskip
	
	\noindent
	$(iv)$ The initial data $u_0$ satisfy the compatibility conditions
$
		u_0(x,0)=0,$ and $ \int_0^\infty u_0\,dy=0.
$
	
	Then there exist a sufficiently small constant $\varepsilon_0>0$ and a positive integer $m_0$, depending only on $\kappa$, such that, for any 	$
	m\ge m_0>\frac{17-24\kappa}{2(1-2\kappa)},
	$, if
	\begin{align}
		\varepsilon \le \varepsilon_0,\qquad
		C_a\le \varepsilon_0,
		\qquad \text{and} \qquad
		\left\| e^{\frac{y^2}{8}}\partial_y^2 u_0 \right\|_{H^{m,0}} \leq \varepsilon_0,
		\label{u0}
	\end{align}
	then system \eqref{uprandtl} admits a unique global-in-time solution $u$ satisfying
		\begin{align*}
		&\sup_{t \in [0, \infty)} \left( \left\| \langle t \rangle^{\frac{5}{4} - \kappa} \mathrm{e}^{\Psi} G(t) \right\|_{H^{5,0}}^2 
		+ \left\| \langle t \rangle^{\frac{5}{4} - \kappa} \mathrm{e}^{\gamma \Psi}  u(t) \right\|_{H^{5,0}}^2
		+ \left\| \langle t \rangle^{\frac{7}{4} - \kappa} \mathrm{e}^{\gamma \Psi} \partial_y u(t) \right\|_{H^{4,0}}^2 
		+ \left\| \langle t \rangle^{\frac{9}{4} - \kappa} \mathrm{e}^{\gamma \Psi} \partial_y^2 u(t) \right\|_{H^{3,0}}^2 \right) \\
		&+ \int_0^\infty \left( \left\| \langle t \rangle^{\frac{5}{4} - \kappa} \mathrm{e}^{\Psi} \partial_y G(t) \right\|_{H^{5,0}}^2 
		+ \left\| \langle t \rangle^{\frac{5}{4} - \kappa} \mathrm{e}^{\gamma \Psi} \partial_y u(t) \right\|_{H^{5,0}}^2
		+ \left\| \langle t \rangle^{\frac{7}{4} - \kappa} \mathrm{e}^{\gamma \Psi} \partial_y^2 u(t) \right\|_{H^{4,0}}^2 
		+ \left\| \langle t \rangle^{\frac{9}{4} - \kappa} \mathrm{e}^{\gamma \Psi} \partial_y^3 u(t) \right\|_{H^{3,0}}^2 \right) dt \\
		&\leq C \varepsilon_0^2,
	\end{align*}
	for any $\gamma \in [\gamma_0,1 )$.
\end{theorem}
\begin{remark}
	For each fixed $\kappa_0\in(0,\frac14)$, the density assumption \eqref{rho-upper-new} is uniform for all $\kappa\in(\kappa_0,\frac12)$. Since $\kappa_0$ is arbitrary in $(0,\frac14)$, these intervals cover the whole range $(0,\frac12)$:
	$
	\bigcup_{0<\kappa_0<1/4}(\kappa_0,\tfrac12)=(0,\tfrac12).
	$
	Hence, for any prescribed $\kappa\in(0,\frac12)$, one may choose $0<\kappa_0<\min\{\kappa,\frac14\}$, and the above argument remains valid.
\end{remark}
\begin{remark}
	Notice that, by the definition of $m$ and equation \eqref{shear flow},
	$
	\partial_{y} (u^s(t,y)+\varepsilon f(t)\chi(y))
$
	satisfies the homogeneous heat equation
$
	(\partial_t-\partial_y^2)\partial_{y} (u^s(t,y)+\varepsilon f(t)\chi(y))=0.
$
	Therefore, it is natural to impose the monotonicity assumption $\left( \partial_y u^s(t,y)+\varepsilon f(t)\chi'(y)\right)
	\geq c_0 \langle t \rangle^{-2+\kappa} e^{-\frac{\gamma_0 y^2}{8 \langle t \rangle}}.
$
\end{remark}

\begin{remark}
 If we choose
	$
	c(\kappa_0)=2\kappa_0,$ and
	$
	\bar\rho=1+\frac{\kappa_0}{2(1-4\kappa_0)},
	$
	then the condition \eqref{Cainfty}
	reduces to
	$
	C_{a,\infty}
	\le
	\frac{\kappa_0(1-4\kappa_0)}{6(m-1)(2-7\kappa_0)}.
	$
	Since  $m$ satisfies
	$
	m\ge m_0>\frac{17-24\kappa}{2(1-2\kappa)},
	$
hence  $m-1\geq 8$, it is sufficient to impose the stronger but simpler condition
	$
	C_{a,\infty}
	\le
	\frac{\kappa_0(1-4\kappa_0)}{48(2-7\kappa_0)}.
	$
Although this condition is restrictive, it is compatible with the present small-data framework, since
$
a=\frac1\rho-1
$
denotes the density perturbation.
%	Although this condition is restrictive, it is natural in the present perturbative setting, since
%	$$
%	a=\frac1\rho-1
%	$$
%	measures the deviation of the density from the constant state.
\end{remark}
\begin{remark}
	Under the assumption \eqref{Ca}, the density perturbation $a$   decays like $\langle t\rangle^{-\frac12-\kappa_0}$ in the Sobolev framework.  Since $\kappa_0$ can be chosen strictly smaller than $1/4$, this decay is still weaker than the one imposed in the analytic framework of Chen--Huang--Li in \cite{CHL2024}, namely,
	$$
	\sup_{t\in[0,+\infty)}
	\left\|\langle t\rangle^{\frac34} e^{\frac{\delta}{2}|D_x|}a(t)\right\|_{B^{1}_{1/2}\cap B^{3}_{1/2}}
	<+\infty.
$$
	The small additional factor $\langle t\rangle^{-\kappa_0}$ beyond the borderline decay $\langle t\rangle^{-1/2}$ is used to make the density--shear-flow coupling terms time integrable in the Sobolev energy estimates.  In this sense, the present result may be viewed as a Sobolev-space counterpart of their result for the two-dimensional compressible Prandtl-type equations, obtained under a weaker temporal decay assumption on the density perturbation.
\end{remark}

The remainder of this paper is organized as follows. In Section \ref{s2}, we explain the main analytical difficulties and collect the auxiliary results used throughout the paper. In Section \ref{s3}, we derive the weighted $L^2$ estimate for $e^{\Psi}g_m$, where the good unknown $g_m$ is introduced to remove the loss of one tangential derivative. Sections \ref{s4}, \ref{s5}, and \ref{s6} are devoted respectively to the weighted decay estimates for $G$, $w$, and $\partial_y w$. In Section \ref{s7}, we combine these a priori estimates to close the bootstrap argument and prove Theorem \ref{t1}.  Appendix \ref{ba} derives the evolution equation and the boundary condition satisfied by $g_m$.

\section{Preliminaries}\label{s2}
\subsection{Difficulties and outline of our approach}
In this subsection, we explain the main difficulties in the proof of Theorem \ref{t1}
and outline the strategy used to overcome them.

\textbf{The first difficulty is to overcome the tangential derivative loss while obtaining enough large-time decay to close a global Sobolev argument.}
Indeed, from the definition of the potential function $\varphi$, the normal velocity is given by
$$
v=-\partial_x\varphi-\frac{\partial_x a}{a+1}\int_y^\infty
\left[u+u^s-\varepsilon f(t)(1-\chi(y'))\right] \,dy',
$$
so that $v$ contains one tangential derivative of $u$ through $\partial_x\varphi$.
Consequently, if one performs a direct $H^{m,0}$ energy estimate on the equation for $u$
or on the vorticity equation for $w=\partial_yu$, then the convection terms involving
$$
v\partial_y\bigl(u+u^s+\varepsilon f(t)\chi(y)\bigr)
\quad\text{and}\quad
v\partial_y\bigl(w+w^s+\varepsilon f(t)\chi'(y)\bigr)
$$
lead to a loss of one $x$-derivative, so the standard Sobolev energy method does not close at top order.
Moreover, removing the derivative loss alone is still not enough for a global result: one also needs sufficiently fast time decay in order to make the coefficients in the weighted energy inequalities integrable and to control the lower-order terms for all time.

To deal with these two aspects simultaneously, we introduce two complementary good quantities.
First, we use the good unknown
$$
g_m
=
\partial_x^m w
-
\frac{\partial_y(w+w^s+\varepsilon f(t)\chi'(y))}{w+w^s+\varepsilon f(t)\chi'(y)}
\partial_x^m u,
$$
which is inspired by the cancellation structure discovered by Masmoudi and Wong in \cite{MW2015}.
The role of $g_m$ is to cancel the most singular contribution generated by the normal velocity,
thereby yielding a closed weighted estimate at the highest tangential order without derivative loss.
Second, we introduce the good quantity
$$
G=u+\frac{y}{2\langle t\rangle}\varphi,
$$
which was proposed in \cite{IV2016,PZ2021}.
The equation \eqref{G uprandtl} satisfied by $G$ contains the favorable damping term
$
\langle t\rangle^{-1}G,
$
which is absent from the equation for $u$ itself.
Combined with the Gaussian weight $\Psi$, this structure provides faster large-time decay for $G$ than for $u$ and gives the extra lower-order control needed in the global argument.
In short, $g_m$ is used to recover the top-order Sobolev estimate, while $G$ is used to recover the decay needed to propagate the global-in-time bounds.

\textbf{The second difficulty lies in the nonstandard boundary structure of both the vorticity equation \eqref{wprandtl} and the good-unknown equation \eqref{gm prandtl}.}
Unlike the standard incompressible Prandtl setting, the vorticity equation \eqref{vprandtl} does not satisfy a simple homogeneous boundary condition such as $\partial_y w|_{y=0}=0$.
Instead, because of the variable density coefficient, the shear-flow contribution, and the cut-off correction, the natural boundary condition is
$$
\bigl((a+1)\partial_y w+a\partial_y w^s+a\varepsilon f(t)\chi''(y)\bigr)\big|_{y=0}=0,
$$
which reduces to
$$
\bigl((a+1)\partial_y w+a\partial_y w^s\bigr)\big|_{y=0}=0,
$$
since $\chi''(0)=0$.
After introducing the good unknown $g_m$, the boundary relation becomes a Robin-type condition involving not only $\partial_y g_m$, but also the lower-order boundary terms produced by the commutator $[\partial_x^m,a+1]\partial_y w$ and the shear-flow contribution $\partial_x^m a\,\partial_y w^s$.
Therefore, the boundary difficulty is not only that the boundary condition is non-homogeneous at the vorticity level, but also that the good-unknown equation inherits additional coupled boundary terms that cannot be ignored in the weighted energy method.

The key point is that this difficulty can still be overcome by preserving the exact flux structure of the equation.
Rather than imposing an artificial homogeneous boundary condition, we retain the natural boundary relations inherited from \eqref{vprandtl} and \eqref{gm prandtl}, and then organize all boundary contributions arising from the diffusion term, the derivative term hidden in the definition of $g_m$, and the commutator terms in a unified manner.
In particular, these pieces can be combined into a single boundary contribution, which can then be controlled by trace inequalities together with the lower-order weighted estimates.
This allows us to close the Sobolev energy estimate without losing the dissipation, and it is one of the distinctive features of the present argument.

A further point is that the boundary difficulty at the level of
$\partial_y w$ is handled by keeping the natural combined structure of the
boundary condition. More precisely, instead of estimating $\partial_y w$
directly, we first regard
$$
\partial_y w+\frac{a}{a+1}\partial_y w^s
$$
as one single quantity. This is consistent with the boundary relation
$$
\left((a+1)\partial_y w+a\partial_y w^s\right)\big|_{y=0}=0,
$$
or equivalently,
$$
\left(
\partial_y w+\frac{a}{a+1}\partial_y w^s
\right)\bigg|_{y=0}=0.
$$
In this way, the boundary contribution coming from the diffusion term can
be treated within the exact flux structure of the equation. Once the
weighted estimate for the above combined quantity is obtained, the
estimate for $\partial_y w$ itself is recovered by using the decay of $a$
and the weighted shear-flow estimates. Consequently, under the assumptions
of Theorem \ref{t1}, we obtain
$$
\|\langle t\rangle^{\frac94-\kappa}e^{\gamma\Psi} \partial_y w\|_{L_T^\infty(H^{3,0})}^2
+
\int_0^T \langle t\rangle^{\frac92-2\kappa}
\left\|e^{\gamma\Psi}\partial_y^2 w\right\|_{H^{3,0}}^2\,dt
\le c_1 .
$$
This final recovery is one of the key steps in closing the continuity
argument.

\textbf{The third difficulty is specific to compressibility and comes from the density-dependent couplings.}
Since the density is not constant and depends only on $(t,x)$, the system is no longer divergence-free,
and several additional terms involving the density coefficient
$
a=\frac1\rho-1
$
appear in the perturbation equations.
For instance, the equation \eqref{uprandtl} for $u$ contains
$
a\partial_y^2u^s$
and
$a\varepsilon f(t)\chi''(y),
$
the equation \eqref{wprandtl} for $w$ contains the terms $
a\partial_y^2 w^s$ and $F_3$,
and the equation \eqref{G uprandtl} for $G$ contains
$
-a\partial_y^2u^s,$ $-a\varepsilon f(t)\chi''(y)$,
and $ -\frac{y}{2\langle t\rangle}(F_1+F_2)$.

These terms have no analogue in the incompressible case and create extra source terms in the weighted Sobolev estimates.
Moreover, the factor
$
\frac{\partial_x a}{a+1}
$
in the expression of $v$ generates additional nonlinear couplings among the density, the perturbation, and the shear flow, which makes the compressible problem substantially more delicate.

Another important point is the choice of the cut-off function $\chi$.
As in the analytic compressible case \cite{CHL2024}, the condition
$
\int_0^\infty (1-\chi(y))\,dy=0
$
is crucial for recovering the boundary relation $\varphi(t,x,0)=0$
from the continuity equation and the positivity of $\rho$.
This fact will be repeatedly used in the derivation of the weighted estimates for $\varphi$ and $G$.

To handle the density-dependent terms together with the shear-flow contribution,
we use the structural assumptions imposed on the shear flow directly.
More precisely, the lower-order pointwise bound on \(\partial_y^2u^s\) in \eqref{CM1}, together with the elementary Gaussian tail estimate, yields the weighted decay of
\[
\|e^\Psi w^s(t)\|_{L_y^2}^2
\lesssim \langle t\rangle^{-\frac52+2\kappa},
\qquad w^s=\partial_yu^s.
\]
This estimate is used to control the lower-order shear-flow terms coupled with the density coefficient \(a\), for example the terms involving
\[
\langle t\rangle^{\frac52-2\kappa}
\Bigl(\|a\|_{H_x^6}^2+
\Bigl\|\frac{\partial_xa}{a+1}\Bigr\|_{H_x^6}^2\Bigr)
\|e^\Psi w^s\|_{L_y^2}^2 .
\]
Here the strengthened decay assumption \eqref{Ca}, namely
\(\|a(t)\|_{H_x^{m+1}}\lesssim \langle t\rangle^{-\frac12-\kappa_0}\),
provides the additional time integrability needed to avoid the borderline logarithmic divergence.
The terms involving the cut-off correction are controlled by the time-decay assumptions on \(f\) in \eqref{Cf}.

At the highest-order level, however, the situation is different.
In the weighted energy estimate for
\[
q:=\partial_yw+\frac{a}{a+1}\partial_yw^s,
\]
the differentiated vorticity equation produces the stronger source term
\[
a\partial_y^2w^s=a\partial_y^3u^s.
\]
Since the coefficient \(a=a(t,x)\) carries no Gaussian weight in the normal variable,
there is no mechanism to lower the weight on \(\partial_y^2w^s\) from \(e^{\gamma\Psi}\) to \(e^{\gamma_0\Psi}\).
Consequently, this term has to be controlled exactly at the same Gaussian level as the unknown in the corresponding highest-order energy estimate, namely at the level \(e^{\gamma\Psi}\).
This explains why the highest-order shear-flow assumption in \eqref{CM1} is imposed on
\(e^{\gamma\Psi}\partial_y^3u^s\), rather than only at the weaker lower-order weight \(e^{\gamma_0\Psi}\).

The density coefficient $a=\frac1\rho-1$ produces several variable-coefficient terms in the equations for $w$ and $\partial_y w$, in particular
$$
(a+1)\partial_y^2 w,\qquad a\partial_y^2u^s,\qquad a\partial_y^2w^s,
$$
as well as the forcing terms generated by $a\varepsilon f(t)\chi''(y)$ and $a\varepsilon f(t)\chi'''(y)$.
To handle these terms, we shall repeatedly use the following weighted heat estimate and its higher-order consequence.

\begin{lemma}\label{heat pro}
	Let $f \in H^{0,1}$ and $g \in L^2$ satisfy
	\begin{align*}
		\left( \partial_t -(a+1)\partial_y^2\right)f =g, \qquad f\partial_y f \big|_{y=0}=0,
	\end{align*}
	where
	$$
	a \buildrel\hbox{\footnotesize def}\over =\frac{1}{\rho}-1,
	\qquad
	\rho=\rho(t,x)\ge 1,
	$$
	and $\rho$ is independent of $y$. Then for any $\gamma \in [0,1]$, we have
	\begin{align}
		\frac12\frac{d}{dt}
		\left\| e^{\gamma\Psi} f\right\|_{L^2}^2
		+	\int_{\mathbb{R}_+^2}
		\left(\frac1{\rho}-\frac{\gamma}{2\rho^2}\right)
		e^{2\gamma\Psi} |\partial_y f|^2\,dxdy
		\leq
		\left\langle e^{\gamma\Psi}f , e^{\gamma\Psi}g \right\rangle_{L^2},
		\label{heat1a}
	\end{align}
	and
	\begin{align}
		\frac12\frac{d}{dt}
		\left\| e^{\gamma\Psi} f\right\|_{L^2}^2
		+\left(1-\frac{\gamma}{2}\right)
		\left\| \frac1{\sqrt\rho}e^{\gamma\Psi}\partial_y f\right\|_{L^2}^2
		\leq
		\left\langle e^{\gamma\Psi}f , e^{\gamma\Psi}g \right\rangle_{L^2}.
		\label{heat1}
	\end{align}
\end{lemma}

\begin{proof}
	Since
	$$
	\partial_t\Psi+2(\partial_y\Psi)^2=0,
	$$
	we first compute
	\begin{align}
		\int_{\mathbb{R}_+^2} e^{2\gamma\Psi} f \partial_t f \,dxdy
		&=
		\frac12\frac{d}{dt}\left\| e^{\gamma\Psi} f\right\|_{L^2}^2
		-\gamma\int_{\mathbb{R}_+^2}\partial_t\Psi\, e^{2\gamma\Psi} f^2\,dxdy.
		\label{heat1.00}
	\end{align}
	On the other hand, since $\rho=\rho(t,x)$ is independent of $y$, integrating by parts gives
	\begin{align}
		\begin{split}
			-\int_{\mathbb{R}_+^2}\frac1{\rho}e^{2\gamma\Psi} f\,\partial_y^2 f\,dxdy
			&=
			\int_{\mathbb{R}_+^2}\frac1{\rho}e^{2\gamma\Psi} |\partial_y f|^2\,dxdy
			+2\gamma\int_{\mathbb{R}_+^2}\frac1{\rho}\partial_y\Psi\, e^{2\gamma\Psi} f\,\partial_y f\,dxdy.
		\end{split}
		\label{heat1.0}
	\end{align}
	Applying Young's inequality to the last term on the right-hand side of \eqref{heat1.0}, we obtain
	\begin{align}
		\begin{split}
			2\gamma\frac1{\rho}\partial_y \Psi\, f\,\partial_y f
			&\ge
			-\frac{\gamma}{2\rho^2}|\partial_y f|^2
			-2\gamma |\partial_y\Psi|^2 f^2.
		\end{split}
		\label{heat1.05}
	\end{align}
	Substituting \eqref{heat1.05} into \eqref{heat1.0}, we arrive at
	\begin{align}
		\begin{split}
			-\int_{\mathbb{R}_+^2}\frac1{\rho}e^{2\gamma\Psi} f\,\partial_y^2 f\,dxdy
			&\ge
			\int_{\mathbb{R}_+^2}
			\left(\frac1{\rho}-\frac{\gamma}{2\rho^2}\right)
			e^{2\gamma\Psi} |\partial_y f|^2\,dxdy
			\\
			&\quad
			-2\gamma\int_{\mathbb{R}_+^2}(\partial_y\Psi)^2 e^{2\gamma\Psi} f^2\,dxdy.
		\end{split}
		\label{heat1.06}
	\end{align}
	Combining \eqref{heat1.00} and \eqref{heat1.06}, and using
	$$
	\partial_t\Psi+2(\partial_y\Psi)^2=0,
	$$
	we deduce that
	\begin{align}
		\begin{split}
			\left\langle e^{\gamma\Psi}f , e^{\gamma\Psi}g \right\rangle_{L^2}
			&\ge
			\frac12\frac{d}{dt}\left\| e^{\gamma\Psi} f\right\|_{L^2}^2
			+\int_{\mathbb{R}_+^2}
			\left(\frac1{\rho}-\frac{\gamma}{2\rho^2}\right)
			e^{2\gamma\Psi} |\partial_y f|^2\,dxdy.
		\end{split}
		\label{heat1.1}
	\end{align}
	Since $\rho\ge1$ and $\gamma\in[0,1]$, we have
	\begin{align*}
		\frac1{\rho}-\frac{\gamma}{2\rho^2}
		=
		\frac1{\rho}\left(1-\frac{\gamma}{2\rho}\right)
		\ge
		\frac1{\rho}\left(1-\frac{\gamma}{2}\right).
	\end{align*}
	Therefore,
	\begin{align*}
		\left\langle e^{\gamma\Psi}f , e^{\gamma\Psi}g \right\rangle_{L^2}
		\ge
		\frac12\frac{d}{dt}\left\| e^{\gamma\Psi} f\right\|_{L^2}^2
		+\left(1-\frac{\gamma}{2}\right)
		\left\| \frac1{\sqrt\rho}e^{\gamma\Psi}\partial_y f\right\|_{L^2}^2,
	\end{align*}
	which proves \eqref{heat1}.
\end{proof}

\begin{corollary}\label{heat cor}
	Under the assumptions of Lemma \ref{heat pro}, let $m\in \mathbb{N}$. Assume in addition that
	\[
	\partial_x^j f\,\partial_y\partial_x^j f\big|_{y=0}=0,
	\qquad 0\le j\le m .
	\]
	Then, for any $\gamma\in [0,1]$, one has
	\begin{align}
		\frac12\frac{d}{dt}\|e^{\gamma\Psi}f\|_{H_x^{m,0}}^2
		+\Bigl(1-\frac{\gamma}{2}\Bigr)
		\sum_{j=0}^{m}
		\left\|
		\frac1{\sqrt{\rho}}e^{\gamma\Psi}\partial_y\partial_x^j f
		\right\|_{L^2}^2
		\le
		\left\langle e^{\gamma\Psi}f,e^{\gamma\Psi}g\right\rangle_{H_x^{m,0}}
		+\sum_{j=1}^m \mathcal R_j,
		\label{heat cor1}
	\end{align}
	where
	\begin{align}
		\mathcal R_j
		:=
		\left\langle
		e^{\gamma\Psi}\partial_x^j f,\,
		e^{\gamma\Psi}\bigl[\partial_x^j,(a+1)\bigr]\partial_y^2 f
		\right\rangle_{L^2}.
		\label{heat cor2}
	\end{align}
\end{corollary}
Combined with  the strengthened time decay assumptions \eqref{Ca} and the smallness \eqref{Cainfty} on $a$,
Lemma \ref{heat pro} and Corollary \ref{heat cor} allow us to control
the variable-coefficient diffusion and commutator terms produced by compressibility.

\subsection{Some useful lemmas}
In this subsection, we collect several auxiliary lemmas that will be repeatedly used in Sections \ref{s3}--\ref{s6} to derive the uniform weighted energy estimates in the Sobolev framework.

We begin with a weighted Poincar\'e-type inequality, which allows us to control weighted $L^2$ norms by weighted normal derivatives and will be used throughout the weighted energy estimates.

\begin{lemma}\label{po}
	Let $f$ be a sufficiently smooth function on $\mathbb{R}_+^2$ that decays sufficiently fast as $y\to+\infty$. Then, for any $\gamma\in [0,1]$ and $\zeta \in (0,1)$, one has
	\begin{align}\label{ptype1}
		\frac{\gamma}{2\left\langle t \right\rangle} \int_{\mathbb{R}_+^2} e^{2\gamma\Psi} |f(x,y)|^2  dxdy \leq \int_{\mathbb{R}_+^2} e^{2\gamma\Psi} |\partial_y f(x,y)|^2  dxdy, 
	\end{align}
	and
	\begin{align}\label{ptype2}
		\frac{\zeta\gamma}{2\left\langle t\right\rangle } \int_{\mathbb{R}_+^2} e^{2\gamma\Psi} |f(x,y)|^2  dxdy + \frac{\zeta(1-\zeta)\gamma}{4\left\langle t\right\rangle ^2} \int_{\mathbb{R}_+^2} y^2 e^{2\gamma\Psi}  |f(x,y)|^2  dxdy \leq \int_{\mathbb{R}_+^2} e^{2\gamma\Psi} |\partial_y f(x,y)|^2  dxdy.
	\end{align}
\end{lemma}
The proof follows from a standard  argument; for similar details, we refer to Refs.~\cite{CHL2024,PZ2021}.
The next lemma relates the good quantity $G$ to the original unknowns $\varphi$, $u$, and $w$.
It allows us to transfer the weighted decay estimate of $G$ to the corresponding weighted estimates of $\varphi$, $u$, and $\partial_yu$.

\begin{lemma}\label{uG}
	Let $\varphi$, $u$  and $w$ be sufficiently smooth solutions of  \eqref{varphi uprandtl}, \eqref{uprandtl} and $\eqref{vprandtl}$ on $\mathbb{R}_+^2$, respectively.  Then for any $\gamma \in (0, 1)$ and $w  \stackrel{\text{def}}{=} \partial_y u$, it holds that
	\begin{align*}
		&\left\| e^{\gamma\Psi} \varphi\right\|_{H^{m,0}} \leq C\left\langle t\right\rangle ^{\frac12}
		\left\| e^{\Psi} G\right\|_{H^{m,0}},\\
		&\left\| e^{\gamma\Psi}  u\right\|_{H^{m,0}} \leq C
		\left\| e^{\Psi} G\right\|_{H^{m,0}},\\
		&\left\| e^{\gamma\Psi} w\right\|_{H^{m,0}} \leq C
		\left\| e^{\Psi}\partial_y G\right\|_{H^{m,0}}.
	\end{align*}
\end{lemma}
\begin{proof}
	Before proving this lemma, we first verify the boundary condition $\varphi|_{y=0}=0$.
	In fact, due to $\eqref{uprandtl}_2$, we get
	\begin{align*}
		\partial_x \left( \rho \int_{\mathbb{R}_+} u \,dy'\right) +\rho v \big|_{0}^{\infty}
		=\varepsilon f(t)\partial_{x}\rho \int_{\mathbb{R}_+} (1-\chi(y)) \,dy'-\partial_x\rho\int_{\mathbb{R}_+} u^s \,dy',
	\end{align*}
	which, together with the boundary condition $v|_{y=0}=0$ and \eqref{chi}, implies that 
	\begin{align}
		\partial_{x}\left( \rho(t,x) \int_{\mathbb{R}_+} u(t,x, y) \,dy \right) =0.
	\end{align}
	Therefore, we have
	\begin{align}
		\rho(t,x)  \int_{\mathbb{R}_+} u(t,x, y) \,dy  =C(t).
	\end{align}
	Since $u$ decays to zero as $|x|\to\infty$ and $\rho(t,x)>0$, we deduce that $C(t)=0$, which yields
	\begin{align*}
		\varphi |_{y=0}=-\int_{0}^\infty u(t,x, y) \,dy=0.
	\end{align*}
	Thanks to the definition of $G$ in \eqref{G}, by solving the following ODE system for $\varphi$:
	\begin{equation}\label{ODE}
		\left\{
		\begin{array}{ll}
			\left( \partial_y+\frac{y}{2\left\langle t \right\rangle } \right) \varphi (t, x, y) =G(t, x, y), \\
			\varphi(t, x, 0)=0,
		\end{array}
		\right.
	\end{equation}
	we deduce that 
	\begin{align*}
		\varphi (t, x, y)= {e^{-\frac{ y^2}{4\left\langle t\right\rangle }}} \int_0^y {e^{\frac{ (y')^2}{4\left\langle t\right\rangle }}} G(t, x, y') \,dy'.
	\end{align*}
	Then for any $0 \leq \gamma \leq 1$ and $m \in \mathbb{N}$, we have
	\begin{align}
		e^{\frac{ \gamma y^2}{8\left\langle t\right\rangle }} \partial^{m}_{x}\varphi(t, x, y)=e^{\frac{ (\gamma-1) y^2}{8\left\langle t\right\rangle }} \int_{0}^{y} e^{\frac{  ((y')^2-y^2)}{8\left\langle t\right\rangle }}
		e^{\frac{  (y')^2}{8\left\langle t\right\rangle }}\partial_{x}^{m} G(t, x, y') \,dy'.
		\label{uG1}
	\end{align}
	On the other hand, notice that
	\begin{align*}
		e^{-\zeta^2}\int_{0}^\zeta e^{z^2} \,dz < \infty
	\end{align*}
	for any $\zeta > 0$. Then we conclude that
	\begin{align}
		\int_{0}^{y} e^{\frac{  ((y')^2-y^2)}{4\left\langle t\right\rangle }} \,dy'
		=2\left\langle t\right\rangle^{\frac12} e^{-\frac{y^2}{4\left\langle t\right\rangle }}
		\int_{0}^{\frac{y}{2\sqrt{\left\langle t\right\rangle }}} e^{z^2} \,dz \leq C\left\langle t\right\rangle^{\frac{1}{2}}.
		\label{uG2}
	\end{align}
	As a result, by Cauchy-Schwarz inequality and in view of \eqref{uG1}-\eqref{uG2}, we have
	\begin{align*}
		\left|\frac{y}{2\left\langle t\right\rangle }e^{\gamma\Psi} \partial^{m}_{x}\varphi(t, x, y)\right|^2
		&\leq \frac{y^2}{4\left\langle t\right\rangle^2 } e^{\frac{ (\gamma-1) y^2}{4\left\langle t\right\rangle }}
		\left(\int_{0}^{y} e^{\frac{  ((y')^2-y^2)}{4\left\langle t\right\rangle }} \,dy' \right)
		\left(\int_{0}^{y} e^{\frac{  (y')^2}{4\left\langle t\right\rangle }} |\partial_{x}^{m} G(t, x, y')|^2 \,dy' \right)\\
		&\leq C\left\langle t\right\rangle^{\frac12}
		\frac{y^2}{4\left\langle t\right\rangle^2 } e^{\frac{ (\gamma-1) y^2}{4\left\langle t\right\rangle }}
		\int_{0}^{y} e^{\frac{  (y')^2}{4\left\langle t\right\rangle }} |\partial_{x}^{m} G(t, x, y')|^2 \,dy'.
	\end{align*}
	Integrating the above inequality over $\mathbb{R}_+^2$ and using Fubini's theorem, we obtain
	\begin{align*}
		\left\|\frac{y}{2\left\langle t\right\rangle }e^{\gamma\Psi} \partial^{m}_{x}\varphi\right\|_{L^2}^2
		&\leq C\left\langle t\right\rangle^{\frac12}
		\int_{\mathbb{R}}\int_0^\infty
		\frac{y^2}{4\left\langle t\right\rangle^2 } e^{\frac{ (\gamma-1) y^2}{4\left\langle t\right\rangle }}
		\int_{0}^{y} e^{\frac{  (y')^2}{4\left\langle t\right\rangle }} |\partial_{x}^{m} G(t, x, y')|^2 \,dy' \,dy\,dx\\
		&= C\left\langle t\right\rangle^{\frac12}
		\int_{\mathbb{R}}\int_0^\infty
		e^{\frac{  (y')^2}{4\left\langle t\right\rangle }} |\partial_{x}^{m} G(t, x, y')|^2
		\left(\int_{y'}^\infty \frac{y^2}{4\left\langle t\right\rangle^2 } e^{\frac{ (\gamma-1) y^2}{4\left\langle t\right\rangle }} \,dy \right)dy' \,dx\\
		&\leq C_{\gamma}\left\| e^{\Psi}\partial_x^m G\right\|_{L^2}^2.
	\end{align*}
	Hence,
	\begin{align*}
		\left\|\frac{y}{2\left\langle t\right\rangle }e^{\gamma\Psi}\varphi\right\|_{H^{m,0}} \leq C_{\gamma}\left\| e^{\Psi}G\right\|_{H^{m,0}} .
	\end{align*}
	Therefore, by the definition of $G$ in \eqref{G}, we get
	\begin{align*}
		\left\| e^{\gamma\Psi}u\right\|_{H^{m,0}}
		&\leq \left\| e^{\gamma\Psi}G\right\|_{H^{m,0}}
		+\left\|\frac{y}{2\left\langle t\right\rangle }e^{\gamma\Psi}\varphi\right\|_{H^{m,0}}\\
		&\leq C_{\gamma}\left\| e^{\Psi}G\right\|_{H^{m,0}} .
	\end{align*}
	Using \eqref{ptype1}, we also obtain
	\begin{align*}
		\left\| e^{\gamma\Psi}\varphi\right\|_{H^{m,0}}
		\leq  \left( \frac{2\left\langle t \right\rangle}{\gamma} \right)^\frac{1}{2}\left\| e^{\gamma\Psi}\partial_y\varphi\right\|_{H^{m,0}}
		= \left( \frac{2\left\langle t \right\rangle}{\gamma} \right)^\frac{1}{2}\left\| e^{\gamma\Psi}u\right\|_{H^{m,0}}
		\leq C_{\gamma} \left\langle t\right\rangle^{\frac12}\left\| e^{\Psi}G\right\|_{H^{m,0}}.
	\end{align*}
	Therefore, in view of the definition of $G$ in \eqref{G}, we have
	\begin{align}
		\left\| e^{\gamma\Psi} \partial_y u\right\|_{H^{m,0}} \leq 
		\left\| e^{\gamma\Psi}\partial_y G\right\|_{H^{m,0}}+\frac{1}{2\left\langle  t \right\rangle }\left\| e^{\gamma\Psi} \varphi \right\|_{H^{m,0}}+\frac{1}{2\left\langle t\right\rangle }
		\left\| e^{\gamma\Psi}yu\right\|_{H^{m,0}}.
		\label{uG3}
	\end{align}
	For the second term on the right-hand side of the above inequality, it follows from \eqref{ptype1} that
	\begin{align}
		\frac{1}{2\left\langle  t \right\rangle }\left\| e^{\gamma\Psi} \varphi \right\|_{H^{m,0}} \leq C_{\gamma}\left\langle t\right\rangle^{-\frac{1}{2}}
		\left\| e^{\Psi}G\right\|_{H^{m,0}} \leq C_{\gamma}\left\| e^{\Psi}\partial_y G\right\|_{H^{m,0}}.
		\label{uG4}
	\end{align}
	Moreover, for any $0 \leq \gamma <1$, we use the fact that $\left\| y^2 e^{-\zeta y^2}\right\|_{L^\infty} \leq C_{\zeta}$ for any positive constant $\zeta$, to obtain
	\begin{align*}
		\frac{1}{2\left\langle t\right\rangle }
		\left\| e^{\gamma\Psi}yu\right\|_{H^{m,0}} 
		&\leq 
		\frac{1}{2\left\langle t\right\rangle }
		\left\| y^2 e^{-(1-\gamma)\Psi}\right\|_{L^{\infty}_y(\mathbb{R}_+)}
		\left\| e^{\frac{1+\gamma}{2}\Psi}u\right\|_{H^{m,0}}
		\\
		&\leq
		C_{\gamma}
		\frac{1}{\left\langle t\right\rangle^{\frac12} }
		\left\| e^{\Psi}G\right\|_{H^{m,0}}
		\leq
		C_{\gamma}
		\left\|  e^{\Psi}\partial_{y}G\right\|_{H^{m,0}},
	\end{align*}
	which, together with \eqref{uG3}-\eqref{uG4}, gives immediately
	\begin{align*}
		\left\| e^{\gamma\Psi} w\right\|_{H^{m,0}} \leq C_{\gamma}
		\left\| e^{\Psi}\partial_y G\right\|_{H^{m,0}}.
	\end{align*}
	We complete the proof of Lemma \ref{uG}.
\end{proof}

In addition, we give Gaussian weighted $L^\infty$-type estimates, which will play an essential role in controlling the nonlinear terms involving $u$ and its normal derivatives arising in the Sobolev energy estimates.
\begin{lemma}\label{infty}
	Let $\Psi$ be defined by \eqref{Psi}. Assume that $f=f(x,y)$ is sufficiently smooth, that
	$f(x,\cdot)\in H^1(\mathbb R_+)$ for each fixed $x$, and that $f(x,y)$ decays to zero
	sufficiently fast as $y\to+\infty$. Then:
	
	\smallskip
	\noindent
	{\rm(i)} For any $\delta>0$, we have
	\begin{align}
		\|e^{\delta\Psi}f\|_{L^\infty}
		\le C_\delta \langle t\rangle^{\frac12}
		\|e^{\delta\Psi}\partial_y f\|_{L^\infty}.
		\label{infty1}
	\end{align}
	
	\noindent
{\rm(ii)} For any \(0<\delta\le \bar\delta\) or any \(0\le \delta<\bar\delta\), we have
\begin{align}
	\|e^{\delta\Psi}f\|_{L_y^\infty}
	\le C_{\delta,\bar\delta}\langle t\rangle^{\frac14}
	\|e^{\bar\delta\Psi}\partial_y f\|_{L_y^2}.
	\label{infty2}
\end{align}
	
	\noindent
{\rm(iii)} For any \(0<\delta\le \bar\delta\) or any \(0\le \delta<\bar\delta\), we have
\begin{align}
	\|e^{\delta\Psi}f\|_{L^\infty}
	\le C_{\delta,\bar\delta}\langle t\rangle^{\frac14}
	\|e^{\bar\delta\Psi}\partial_y f\|_{H^{1,0}}.
	\label{infty3}
\end{align}
\end{lemma}
\begin{remark}
	We emphasize that the endpoint case \(\delta=\bar\delta\) is admissible provided that \(\delta=\bar\delta>0\). Indeed, in this case the Gaussian decay still yields
	\[
	\sup_{y\ge0} e^{\delta\Psi(t,y)}
	\left(\int_y^\infty e^{-2\delta\Psi(t,y')}\,dy'\right)^{1/2}
	\le C_\delta \langle t\rangle^{1/4}.
	\]
	However, the case \((\delta,\bar\delta)=(0,0)\) must be excluded, since then
	\[
	\int_y^\infty 1\,dy'=+\infty,
	\]
	and the corresponding estimate fails in general.
\end{remark}
\begin{remark}
	The preceding lemma can be applied to the shear flow $u^s$ and its normal derivatives. 
	Indeed, by the boundary condition in \eqref{shear flow}, one has
	$
	\lim_{y\to+\infty}u^s(t,y)=0.
	$
	Moreover, since
	$
	w^s=\partial_yu^s,
	$
	and
	$ 
	\partial_yw^s=\partial_y^2u^s,
	$
	the assumption on $\partial_y^2u^s$ gives
	$$
	|\partial_yw^s(t,y)|
	\le
	C_s\langle t\rangle^{-2+\kappa}
	\exp\left\{
	-\frac{(2-\gamma)y^2}{8\langle t\rangle}
	\right\}.
	$$
	Hence $\partial_yw^s(t,\cdot)\in L^1(\mathbb R_+)$ and 
	$\partial_yw^s(t,y)\to0$ as $y\to+\infty$. It follows that 
	$w^s(t,y)$ has a finite limit as $y\to+\infty$. This limit must be zero; otherwise, if
	$
	\lim_{y\to+\infty}w^s(t,y)=\ell(t)\neq0,
	$
	then $u^s(t,y)=u^s(t,Y)+\int_Y^y w^s(t,z)\,dz$ would be unbounded as
	$y\to+\infty$, contradicting $\lim_{y\to+\infty}u^s(t,y)=0$. Therefore,
	$$
	\lim_{y\to+\infty}w^s(t,y)=0.
	$$
	Consequently, $u^s$, $w^s$, and $\partial_yw^s$ all satisfy the far-field decay
	condition required in Lemma \ref{infty}. Hence Lemma \ref{infty} and Lemma \ref{po} may be used for
	these shear-flow quantities, provided that the corresponding weighted norm of their
	normal derivative is finite.
\end{remark}

To handle the product and commutator terms arising from the nonlinearities, we also recall the following standard Moser-type inequalities.
\begin{lemma}\label{moser}
	Suppose that $f, g \in L^\infty \cap H^s(\mathbb{R}^d)$ with $s \geq d/2$. Then for $1 \leq m \leq s$ (where $m$ denotes the order of the partial derivative), we have
	\begin{align}
		\|\partial_x^m(fg)\|_{L^2} \leq C\|f\|_{L^\infty}\|\partial_x^m g\|_{L^2} + C\|g\|_{L^\infty}\|\partial_x^m f\|_{L^2},\label{m1}
	\end{align}
	and
	\begin{align}
		\|\partial_x^m(fg) - f \partial_x^m g\|_{L^2} \leq C\|\partial_xf\|_{L^\infty}\|\partial_x^{m-1} g\|_{L^2}+C\|g\|_{L^\infty}\|\partial_x^m f\|_{L^2} .\label{m2}
	\end{align}
\end{lemma}

Finally, we state the following Gagliardo--Nirenberg inequality, which will be used to interpolate between different orders of derivatives in the subsequent analysis.
\begin{lemma}(Gagliardo-Nirenberg Inequality)\label{GN}
	For $0 \leq m, \alpha \leq \ell$ and $2 \leq p \leq \infty$, we have
	\begin{align}
		\|\nabla^\alpha f\|_{L^p} \lesssim \|\nabla^m f\|_{L^2}^{1-\theta} \|\nabla^\ell f\|_{L^2}^\theta, \label{GNI}
	\end{align}
	where $0 \leq \theta \leq 1$ and $\alpha$ satisfies
	\[
	\alpha + 3\left( \frac{1}{2} - \frac{1}{p} \right) = m(1-\theta) + \ell\theta.
	\]
	In particular, when $p = \infty$, we require $0 < \theta < 1$, $m \leq \alpha + 1$ and $\ell \geq \alpha + 2$ here.
\end{lemma}

\section{Energy estimate of $g_m$}\label{s3}

By the existing local well-posedness theory \cite{FRY2021} for the compressible Prandtl equations, it is natural to derive a priori estimates on the lifespan of a local smooth solution.
Accordingly, let $T^*>0$ be the maximal existence time of a sufficiently smooth solution to system \eqref{uprandtl}.
The aim of this section is to derive a closed weighted $L^2$ estimate for $e^\Psi g_m$ on $[0,T^\#]$ by the energy method.
This will allow us to overcome the loss of one tangential derivative caused by the normal velocity.

Recall that the good unknown $g_m$ is defined by \eqref{g}.
For each $0\le j\le m$, we define $g_j$ in the same way as in \eqref{g}, with $m$ replaced by $j$.
Moreover, $g_m$ solves the following system (see Appendix \ref{ba} for the derivation):
\begin{equation}\label{gm prandtl}
	\left\{
	\begin{array}{ll}
		\partial_t g_m -(a+1)\partial_y^2 g_m+(u + u^s + \varepsilon f(t)\chi(y))\partial_x g_m
		\\
		\quad
		-2(a+1)g_m  \partial_y\left(  \dfrac{\partial_y  (w+  w^s+\varepsilon f(t)\chi'(y))}{w + w^s + \varepsilon f(t)\chi'(y)}\right)
		\\
		=\partial_x^m u \left( \partial_x w+\partial_y\left( \dfrac{\partial_y(w + w^s + \varepsilon f(t)\chi'(y))}{w + w^s + \varepsilon f(t)\chi'(y)} v\right) -\partial_y\left( \dfrac{ F_3 }{w + w^s + \varepsilon f(t)\chi'(y)} \right)\right)
		\\
		\quad+(w + w^s + \varepsilon f(t)\chi'(y))\partial_y\left( \dfrac{\left[ u, \partial_{x}^m\right] \partial_x u
			+\left[ w, \partial_x^m\right] v
			+[a, \partial_x^m]\partial_y^2 (u + u^s)}{w + w^s + \varepsilon f(t)\chi'(y)}\right)
		\\
		\quad+
		\partial_{x}^{m}a\varepsilon f(t)
		\left( \chi'''(y)-\dfrac{\partial_y(w + w^s + \varepsilon f(t)\chi'(y))}{w + w^s + \varepsilon f(t)\chi'(y)}\chi''(y)\right)
		+\partial_{x}^{m}F_3,
		\\
		\Biggl[
		(a+1)\left(
		\partial_y g_m
		+\dfrac{\partial_y(w+w^s+\varepsilon f(t)\chi'(y))}{w+w^s+\varepsilon f(t)\chi'(y)}\,g_m
		\right)
		\\
		\qquad\qquad
		+[\partial_x^m,a+1]\partial_y w
		+\partial_x^m a\,\partial_y w^s
		\Biggr]\bigg|_{y=0}=0,
		\\
		g_{m}\big|_{t=0}=g_{m,0}.
	\end{array}
	\right.
\end{equation}

In order that the good unknown $g_m$ be well defined, we need to preserve the positivity of
$w+w^s+\varepsilon f(t)\chi'(y)$.
To this end, we introduce
\begin{equation}
	T^{\#} \stackrel{\text{def}}{=} \sup\left\{ T \in [0, T^*) :  w(x,t,y)+w^s(t,y)+ \varepsilon f(t)\chi'(y) \geq \frac{c_0}{2} \langle t \rangle^{-2+\kappa} e^{-\gamma_0 \Psi}, \forall\ t \in [0,T] \right\}.
	\label{T}
\end{equation}

We first establish the relation between $g_m$ and $\partial_x^m w$.

\begin{lemma} \label{gmw}
	Let $g_m$ be given by \eqref{g}. For
	any $t \leq T^{\#}$ and $0<\gamma<1$,  we have
	\begin{align}
		\|e^\Psi g_m\|_{L^2 } \leq C_{c_0, \gamma} \left(1 + \langle t \rangle^{\frac{5}{2}-\kappa} \|e^{\gamma_0 \Psi} \partial_y(w + w^s + \varepsilon f(t)\chi'(y))\|_{L^\infty }\right) \|e^\Psi \partial_x^m w\|_{L^2 }.\label{gmin1}
	\end{align}
	and
	\begin{align}
		\|e^{\gamma\Psi} \partial_x^m w\|_{L^2 }\leq C_{c_0, \gamma_0, \gamma} \left(1 + \langle t \rangle^{\frac{5}{2}-\kappa} \|e^{\gamma_0 \Psi} \partial_y(w + w^s + \varepsilon f(t)\chi'(y))\|_{L^\infty }\right) \|e^\Psi g_m\|_{L^2 },
		\label{gmin2}
	\end{align}
\end{lemma}
\begin{proof}
	Recalling the definition of $g_m$, we have
	$$
	g_m=\partial_x^m w- \frac{\partial_y(w + w^s + \varepsilon f(t)\chi'(y))}{w + w^s + \varepsilon f(t)\chi'(y)}\partial_{x}^{m}u.
	$$
	Hence, for any $t \leq T^{\#}$,
	\begin{align*}
		\begin{split}
			\|e^\Psi g_m\|_{L^2 }
			&\leq \|e^\Psi \partial_x^m w\|_{L^2 }
			\\
			&\quad + \|e^{-\gamma_0 \Psi}(w + w^s + \varepsilon f(t)\chi'(y))^{-1}\|_{L^\infty }
			\|e^{\gamma_0 \Psi} \partial_y (w+w^s+\varepsilon f(t)\chi'(y))\|_{L^\infty }
			\|e^\Psi \partial_x^m u\|_{L^2 }.
		\end{split}
	\end{align*}
	By the definition of $T^\#$, we have
	$$
	\|e^{-\gamma_0 \Psi}(w + w^s + \varepsilon f(t)\chi'(y))^{-1}\|_{L^\infty }
	\leq
	C_{c_0}\langle t\rangle^{2-\kappa}.
	$$
	Therefore,
	\begin{align*}
		\begin{split}
			\|e^\Psi g_m\|_{L^2 }
			&\leq
			\|e^\Psi \partial_x^m w\|_{L^2 }
			+
			C_{c_0}\langle t\rangle^{2-\kappa}
			\|e^{\gamma_0 \Psi} \partial_y (w+w^s+\varepsilon f(t)\chi'(y))\|_{L^\infty }
			\|e^\Psi \partial_x^m u\|_{L^2 }.
		\end{split}
	\end{align*}
	Applying \eqref{ptype1} to $\partial_x^m u$ and using $\partial_y\partial_x^m u=\partial_x^m w$, we obtain
	$$
	\|e^\Psi \partial_x^m u\|_{L^2 }\leq C_\gamma \langle t\rangle^{1/2}\|e^\Psi \partial_x^m w\|_{L^2 }.
	$$
	Substituting this bound into the previous inequality yields \eqref{gmin1}.
	
	On the other hand, using the alternative expression of $g_m$ in \eqref{g},
	$$
	g_m=
	(w + w^s + \varepsilon f(t)\chi'(y))\partial_y\left( \frac{\partial_x^m u}{w + w^s + \varepsilon f(t)\chi'(y)}\right),
	$$
	we deduce, since $\partial_x^m u|_{y=+\infty}=0$, that
	\begin{align*}
		\partial_x^m u
		=
		-(w + w^s + \varepsilon f(t)\chi'(y))
		\int_y^{\infty}\frac{g_m(y')}{w(y') + w^s(y') + \varepsilon f(t)\chi'(y')} \, dy'.
	\end{align*}
	Therefore,
	\begin{align*}
		\partial_x^m w
		=
		g_m-\partial_y(w + w^s + \varepsilon f(t)\chi'(y))
		\int_y^{\infty}\frac{g_m(y')}{w(y') + w^s(y') + \varepsilon f(t)\chi'(y')} \, dy'.
	\end{align*}
	Hence,
	\begin{align*}
		\begin{split}
			\|e^{\gamma\Psi} \partial_x^m w\|_{L^2 }
			&\leq \|e^{\gamma\Psi} g_m\|_{L^2 }
			\\
			&\quad +C_{c_0}\langle t\rangle^{2-\kappa}
			\|e^{\gamma_0 \Psi} \partial_y (w+w^s+\varepsilon f(t)\chi'(y))\|_{L^\infty }
			\left\| e^{(\gamma-\gamma_0)\Psi}\int_y^{\infty}e^{\gamma_0\Psi(y')}{g_m(y')} \, dy'\right\|_{L^2 }.
		\end{split}
	\end{align*}
	To estimate the last term, we use the facts that for any $\gamma' < 0$, $\gamma \in (0, 1)$, and $\alpha \in \mathbb{N}$,
	$$
	\sup_{y \in \mathbb{R}_+} \left( e^{-\gamma'\Psi} \int_y^\infty e^{{\gamma'}\Psi} \, dy' \right) \leq C_{\gamma'} \langle t \rangle^{\frac{1}{2}},
	$$
	and
	$$
	\left\| y^{\alpha}e^{-\frac{1-\gamma}{2}\Psi}\right\|_{L^2_y(\mathbb{R}_+)} \leq C_\gamma \left\langle t\right\rangle^{\frac{1}{4}+\frac{\alpha}{2}}.
	$$
	Then, by the Cauchy--Schwarz inequality,
	\begin{align*}
		& \left\| e^{(\gamma-\gamma_0)\Psi}\int_y^{\infty}e^{\gamma_0\Psi(y')}{g_m(y')} \, dy'\right\|_{L^2 }
		\\
		&\lesssim \left\| e^{(\gamma-1)\Psi} (e^{2(1-\gamma_0)\Psi})^{1/2}
		\left( \int_y^{\infty}e^{2(\gamma_0-1)\Psi(y')} \, dy'\right)^{1/2}
		\left( \int_y^{\infty}\left| e^{\Psi(y')}{g_m(y')}\right|^2 \, dy'\right)^{1/2}
		\right\|_{L^2 }
		\\
		&\leq C_{\gamma_0}\left\langle t\right\rangle^{\frac{1}{4}}
		\left\| e^{(\gamma-1)\Psi}
		\left( \int_0^{\infty}\left| e^{\Psi(y')}{g_m(y')}\right|^2 \, dy'\right)^{1/2}
		\right\|_{L^2 }
		\\
		&\leq C_{\gamma_0}\left\langle t\right\rangle^{\frac{1}{4}} \left\| e^{{\gamma-1}\Psi}\right\|_{L^2_y(\mathbb{R}_+)} \left\| e^{\Psi}{g_m}\right\|_{L^2 }
		\\
		&\leq C_{\gamma_0, \gamma}\left\langle t\right\rangle^{\frac{1}{2}}  \left\| e^{\Psi}{g_m}\right\|_{L^2 }.
	\end{align*}
	Substituting this bound into the previous estimate, we obtain \eqref{gmin2}.
	This completes the proof.
\end{proof}

With the above equivalence estimate at hand, we are now in a position to derive
a closed weighted energy estimate for $g_m$.

\begin{proposition}\label{pgm}
	Fix $0<\kappa_0<\frac14$, and let $c(\kappa_0)$ be a positive constant such that
	$
		0<c(\kappa_0)<3\kappa_0.
	$
	For any $\kappa\in(\kappa_0,\frac12)$, suppose that the density function $\rho=\rho(t,x)$ satisfies
	\begin{equation*}
		1\le \rho(t,x)\le \bar\rho<
		\frac{1-c(\kappa_0)-\kappa_0}{1-4\kappa_0}.
	\end{equation*}
	Moreover, choose $C_{a,\infty}$ such that
	\begin{align*}
		\max_{1 \leq i \leq m}{  C_{m}^{i}\left\| {\partial_x^{i} a} \right\|_{L^\infty_x}}=	C_{a,\infty}\le
		\frac{1-c(\kappa_0)-\kappa_0-(1-4\kappa_0)\bar\rho}{6(m-1)\bar\rho}.
	\end{align*}
	Let $g_m$ be determined by \eqref{g}. Then there exists a positive constant $C$ such that for all
	$
	0<2\delta\le \gamma_0\le \gamma<1,\,
	t\leq T^{\#},
	$
	one has
	\begin{align}
		\begin{split}
			&\frac{d}{dt} \sum_{j=0}^m\left\| e^{\Psi} g_j\right\|_{L^2}^2 +\kappa_0\sum_{j=0}^m
			\left\| \frac{1}{\sqrt\rho}e^{\Psi} \partial_y g_j\right\|_{L^2}^2 
			\\
			&
			\leq 		C \theta_1 \sum_{j=0}^m\left\| e^{\Psi} g_j\right\|_{L^2}^2
			+C\left( \left( 1+\left\langle t\right\rangle^{\frac32-\kappa} \right) \left\| a\right\|_{H^{m+1}_x }^2\left\| e^{\Psi}\partial_y w^s\right\|_{L^2_y }^2
			+
			\varepsilon^2 \left\langle t\right\rangle |f(t)|^2 \right),
		\end{split}
		\label{g0}
	\end{align}
	where
	\begin{align}
		\begin{split}
			\theta_1(t) 
			&=\left\langle t\right\rangle ^{\frac{9}{2}-2\kappa}\left( 1+\left\| a\right\|_{H^{m+1}_x}^2 \right) 
			\left\|e^{\gamma_0 \Psi} \partial_y^2(w + w^s + \varepsilon    f(t)\chi'(y))\right\|_{L^\infty_xL^2_y }^2
			\\
			&
			\quad
			+
			\left\langle t\right\rangle ^{4-2\kappa}\left( 1+\left\| a\right\|_{H^{m+1}_x}^2 \right) 
			\left\|e^{\gamma_0 \Psi} \partial_y(w + w^s + \varepsilon    f(t)\chi'(y))\right\|_{L^\infty }^2
			\\
			&\quad+
			\left\langle t\right\rangle^{\kappa-\frac{3}{2}}\left( \left\langle t\right\rangle^{{\frac94}-\kappa}\left\| e^{\gamma\Psi} \partial_y w\right\|_{H^{1,0} }\right) +\left\langle t\right\rangle^{-1}\left\| a\right\|_{H^{m+1}_x}^2 
			\\
			&\quad
			+ \left(1 + \langle t \rangle^{\frac{5}{2}-\kappa} \|e^{\gamma_0 \Psi} \partial_y(w + w^s + \varepsilon f(t)\chi'(y))\|_{L^\infty }\right)^2\left( 1+ \left\| \frac{\partial_x a}{a+1} \right\|_{H^m_x }^2\right) 
			\\
			&\qquad\cdot 
			\left( \left\langle t\right\rangle^{\frac{3}{4}}\left\| e^{\gamma\Psi} \partial_y w \right\|_{H^{2,0} }+\left\langle t\right\rangle^{\frac52}\left\| e^{\gamma\Psi} \partial_y w \right\|_{H^{2, 0} }^2+
			\left\langle t\right\rangle \left\| e^{\gamma_0\Psi} \partial_y w^s \right\|_{L^\infty_y } \right. 
			\\
			&\qquad \quad
			\left. +\left\langle t\right\rangle^{3}\left\| e^{\gamma_0\Psi}\partial_y w^s \right\|_{L^\infty_y }^2+\langle t\rangle^{-4+2\kappa}
			+
			\langle t\rangle^{-\frac32+\kappa}\right)
			\\
			&
			\quad+ \left(1 + \langle t \rangle^{\frac{5}{2}-\kappa} \|e^{\gamma_0 \Psi} \partial_y(w + w^s + \varepsilon f(t)\chi'(y))\|_{L^\infty }\right)^2  
			\\
			&\qquad
			\cdot 
			\left(   \varepsilon \left\langle t\right\rangle^{-\frac12}\left\| \frac{\partial_x a}{a+1} \right\|_{L^\infty_x } |f(t)| +\varepsilon^2\left\| \frac{\partial_x a}{a+1} \right\|_{L^\infty_x }^2 |f(t)|^2 \right). 
		\end{split}
		\label{theta1}
	\end{align}
\end{proposition}

\begin{proof}
Multiplying the equation $\eqref{gm prandtl}_1$ by $e^{2\Psi}g_m$, and using the integration by parts in the $y$-variable, and Lemma \ref{heat pro}, we have
\begin{align}
	\begin{split}
		&\frac12\frac{d}{dt}
		\left\| e^{\Psi} {g_m}\right\|_{L^2}^2 + \frac{1}{2}\left\| \frac{1}{\sqrt{\rho}}e^{\Psi} \partial_y {g_m}\right\|_{L^2}^2 \\
		& \leq I_{\mathrm{boundary}} -\int_{\mathbb{R}_+^2} \partial_y \left( e^{2\Psi}g_m(w + w^s + \varepsilon f(t)\chi'(y))\right) \left( \frac{[\partial_x^m, (a+1)]\partial_y^2 u }{w + w^s + \varepsilon f(t)\chi'(y)}\right)
		\\&\quad  -\int_{\mathbb{R}_+^2} e^{2\Psi}g_m (u + u^s + \varepsilon f(t)\chi(y))\partial_x g_m \\
		& \quad- 2\int_{\mathbb{R}_+^2} \partial_y \left(e^{2\Psi}g_m (a+1)g_m  \right) \left(  \frac{\partial_y  (w+  w^s+\varepsilon f(t)\chi'(y))}{w + w^s + \varepsilon f(t)\chi'(y)}\right) \\
		&\quad+ \int_{\mathbb{R}_+^2} e^{2\Psi}g_m\partial_x^m u \left\{\partial_x w+\partial_y\left( \frac{\partial_y(w + w^s + \varepsilon f(t)\chi'(y))}{w + w^s + \varepsilon f(t)\chi'(y)} v\right) \right\}
		\\
		&\quad- \int_{\mathbb{R}_+^2} e^{2\Psi}g_m\partial_x^m u \partial_y\left( \frac{ F_3 }{w + w^s + \varepsilon f(t)\chi'(y)} \right)
		\\
		&\quad- \int_{\mathbb{R}_+^2} \partial_y \left( e^{2\Psi}g_m(w + w^s + \varepsilon f(t)\chi'(y))\right) \left( \frac{\left[ u, \partial_{x}^m\right] \partial_x u
			+\left[ w, \partial_x^m\right] v
			+\partial_x^m a\partial_y^2 u^s}{w + w^s + \varepsilon f(t)\chi'(y)}\right)
		\\
		&\quad+ \int_{\mathbb{R}_+^2} e^{2\Psi}g_m 
		\partial_{x}^{m}a\varepsilon f(t)
		\left( \chi'''(y)-\frac{\partial_y(w + w^s + \varepsilon f(t)\chi'(y))}{w + w^s + \varepsilon f(t)\chi'(y)}\chi''(y)\right) 
		+\int_{\mathbb{R}_+^2} e^{2\Psi}g_m \partial_{x}^{m}F_3\\
		&=: I_{\mathrm{boundary}}+\sum_{i=1}^{8}.
	\end{split}
	\label{gm1.1}
\end{align}
Here
\begin{align}
	I_{\mathrm{boundary}}
	:=
	I_{\mathrm{diff,boundary}}
	+I_{1,\mathrm{boundary}}
	+I_{3,\mathrm{boundary}}
	+I_{6,\mathrm{boundary}},
	\label{gm-bdry-total}
\end{align}
with
\begin{align*}
	I_{\mathrm{diff,boundary}}
	&:=-\int_{\mathbb R}(a+1)g_m\partial_y g_m\big|_{y=0}\,dx,
	\\
	I_{1,\mathrm{boundary}}
	&:=-\int_{\mathbb R} g_m\,[\partial_x^m,(a+1)]\partial_y w\big|_{y=0}\,dx,
	\\
	I_{3,\mathrm{boundary}}
	&:=-2\int_{\mathbb R}(a+1)\frac{\partial_y(w + w^s + \varepsilon f(t)\chi'(y))}{w + w^s + \varepsilon f(t)\chi'(y)}|g_m|^2\big|_{y=0}\,dx,
	\\
	I_{6,\mathrm{boundary}}
	&:=-\int_{\mathbb R} g_m
	\Bigl(
	[u,\partial_x^m]\partial_xu+[w,\partial_x^m]v+\partial_x^m a\,\partial_y^2u^s
	\Bigr)\big|_{y=0}\,dx.
\end{align*}
Since $u|_{y=0}=0$ and $v|_{y=0}=0$, one has
$$
[u,\partial_x^m]\partial_xu\big|_{y=0}=0,\qquad
[w,\partial_x^m]v\big|_{y=0}=0,
$$
and thus
$$
I_{6,\mathrm{boundary}}
=
-\int_{\mathbb R} g_m\,\partial_x^m a\,\partial_y w^s\big|_{y=0}\,dx.
$$
Moreover, by the boundary condition in \eqref{gm prandtl}, we have
\begin{align*}
	I_{\mathrm{diff,boundary}}
	&=
	\int_{\mathbb R}(a+1)\frac{\partial_y(w + w^s + \varepsilon f(t)\chi'(y))}{w + w^s + \varepsilon f(t)\chi'(y)}|g_m|^2\big|_{y=0}\,dx
	+\int_{\mathbb R}g_m\,[\partial_x^m,a+1]\partial_y w\big|_{y=0}\,dx
	\\
	&\quad
	+\int_{\mathbb R}g_m\,\partial_x^m a\,\partial_y w^s\big|_{y=0}\,dx.
\end{align*}
Hence
\begin{align}
	\begin{split}
		I_{\mathrm{boundary}}
		&=
		-\int_{\mathbb R}(a+1)\frac{\partial_y(w + w^s + \varepsilon f(t)\chi'(y))}{w + w^s + \varepsilon f(t)\chi'(y)}|g_m|^2\big|_{y=0}\,dx
		.
	\end{split}
	\label{gm-bdry-combined}
\end{align}
By the trace inequality
$$
\|g_m(\cdot,0)\|_{L_x^2}^2
\le C \|e^\Psi g_m\|_{L^2}\|e^\Psi\partial_y g_m\|_{L^2},
$$
it follows that
\begin{align}
	\begin{split}
		|I_{\mathrm{boundary}}|
		&\le \frac{1}{6C_1}\left\| \frac{1}{\sqrt{\rho}}e^{\Psi}\partial_y g_m\right\|_{L^2 }^2
		+C\left\langle t\right\rangle ^{4-2\kappa}
		\left\|e^{\gamma_0 \Psi} \partial_y(w + w^s + \varepsilon f(t)\chi'(y))\right\|_{L^\infty }^2 
		\left\| e^{\Psi} g_m\right\|_{L^2 } ^2
		.
	\end{split}
	\label{gm-bdry-est}
\end{align}

Now, we estimate the right-hand side of \eqref{gm1.1} term by term as follows.

\textbf{Dealing with $I_1$  term :}
For $I_1$, we have
\begin{align}
	\begin{split}
		I_1
		&\leq
		\left\| \frac{1}{\sqrt{\rho}} \frac{1}{a+1}e^{\Psi} [\partial_x^m, (a+1)]\partial_y^2 u\right\|_{L^2 }
		\left( \left\| \frac{1}{\sqrt{\rho}} e^{\Psi} \partial_y g_m\right\|_{L^2} +\frac{1}{2\left\langle t\right\rangle }\left\|\frac{1}{\sqrt{\rho}} y e^{\Psi}  g_m\right\|_{L^2} \right)
		\\
		&\quad+\left\| e^{\Psi}\frac{\partial_y(w + w^s + \varepsilon f(t)\chi'(y))}{w + w^s + \varepsilon f(t)\chi'(y)}[\partial_x^m, (a+1)]\partial_y^2 u\right\|_{L^2 }
		\left\| e^{\Psi}  g_m\right\|_{L^2 }
		.
	\end{split}
	\label{gm11.22}
\end{align}
On the one hand, using  \eqref{ptype2} in Lemma \ref{po} and the relation
\begin{align*}
	\partial_x^n \partial_yw &=	\partial_yg_n+ \partial_y\left( \frac{\partial_y(w + w^s + \varepsilon f(t)\chi'(y))}{w + w^s + \varepsilon f(t)\chi'(y))}\partial_{x}^{n}u\right) 
	\\
	&
	=\partial_yg_n+  \frac{\partial_y(w + w^s + \varepsilon f(t)\chi'(y))}{w + w^s + \varepsilon f(t)\chi'(y)}g_n+
	\partial_y^2(w + w^s + \varepsilon f(t)\chi'(y))\partial_y^{-1}\left(\frac{ g_n}{w + w^s + \varepsilon f(t)\chi'(y))}\right) ,
\end{align*}
we obtain
\begin{align*}
	&\left\| \frac{1}{\sqrt{\rho}} \frac{1}{a+1}e^{\Psi} [\partial_x^m, (a+1)]\partial_y^2 u\right\|_{L^2 }
	\left( \left\| \frac{1}{\sqrt{\rho}} e^{\Psi} \partial_y g_m\right\|_{L^2} +\frac{1}{2\left\langle t\right\rangle }\left\|\frac{1}{\sqrt{\rho}}e^{\Psi} y e^{\Psi}  g_m\right\|_{L^2} \right)
	\\
	&\leq 3\left\| \frac{1}{\sqrt{\rho}}e^{\Psi} e^{\Psi} \partial_y g_m\right\|_{L^2}  C_{m}^{n}\bar\rho\left\| {\partial_x^{m-n} a} \right\|_{L^\infty_x} \sum_{n=0}^{m-1}\left\| \frac{1}{\sqrt{\rho}}e^{\Psi}  \partial_y\partial_x^{n} w\right\|_{L^2}
	\\
	&
	\leq 3C_{a,\infty} \bar\rho\left\| \frac{1}{\sqrt{\rho}} e^{\Psi} \partial_y g_m\right\|_{L^2} \sum_{n=0}^{m-1}   \left\| \frac{1}{\sqrt{\rho}}e^{\Psi}\partial_yg_{n}\right\|_{L^2} +  \frac{1}{6C_1}\left\|\frac{1}{\sqrt{\rho}} e^{\Psi}\partial_yg_{m}\right\|_{L^2}^2
	\\
	&\quad+C_{C_1, \bar{\rho}} \left\langle t\right\rangle ^{4-2\kappa}\left\| a\right\|_{H^{m+1}_x}^2 
	\|e^{\gamma_0 \Psi} \partial_y(w + w^s + \varepsilon f(t)\chi'(y))\|_{L^\infty }^2 \sum_{n=0}^{m-1}
	\left\| e^{\Psi} g_n\right\|_{L^2 } ^2
	\\
	&\quad+C_{C_1, \bar{\rho}} \left\langle t\right\rangle ^{\frac{9}{2}-2\kappa}\left\| a\right\|_{H^{m+1}_x}^2 
	\|e^{\gamma_0 \Psi} \partial_y^2(w + w^s + \varepsilon    f(t)\chi'(y))\|_{L^\infty_xL^2_y }^2 \sum_{n=0}^{m-1}
	\left\| e^{\Psi} g_n\right\|_{L^2 } ^2,
\end{align*}
where we have used the following fact:
\begin{align*}
	&\left\| e^{\Psi}\partial_y^2(w + w^s + \varepsilon f(t)\chi'(y))\partial_y^{-1}\left(\frac{ g_n}{w + w^s + \varepsilon f(t)\chi'(y))}\right)\right\|_{L^2} 
	\\
	&\leq
	\left\| e^{\gamma_0\Psi}\partial_y^2(w + w^s + \varepsilon f(t)\chi'(y))\right\|_{L_x^\infty L_y^2} \left\| e^{(1-\gamma_0)\Psi}\partial_y^{-1}\left(\frac{ g_n}{w + w^s + \varepsilon f(t)\chi'(y))}\right)\right\|_{L^2_xL^\infty_y} 
	\\
	&\leq \left\langle t\right\rangle ^{\frac{9}{4}-\kappa} \left\| e^{\gamma_0\Psi}\partial_y^2(w + w^s + \varepsilon    f(t)\chi'(y))\right\|_{L_x^\infty L^2_y} \left\| e^{\Psi}g_n\right\|_{L^2}.
\end{align*}
On the other hand by a similar method, we can find that
\begin{align*}
	&\left\| e^{\Psi}\frac{\partial_y(w + w^s + \varepsilon f(t)\chi'(y))}{w + w^s + \varepsilon f(t)\chi'(y)}[\partial_x^m, (a+1)]\partial_y^2 u\right\|_{L^2 }
	\left\| e^{\Psi}  g_m\right\|_{L^2 }
	\\
	&
	\leq C_{C_a,C_1,\bar{\rho}} \left\langle t\right\rangle ^{4-2\kappa}
	\|e^{\gamma_0 \Psi} \partial_y(w + w^s + \varepsilon f(t)\chi'(y))\|_{L^\infty }^2 
	\left\| e^{\Psi} g_m\right\|_{L^2 } ^2 +  \frac{1}{6C_1}\left\|\frac{1}{\sqrt{\rho}} e^{\Psi}\partial_yg_{m}\right\|_{L^2}^2
	\\
	&\quad+C_{C_1, \bar{\rho}} \left\langle t\right\rangle ^{4-2\kappa}\left\| a\right\|_{H^{m+1}_x}^2 
	\|e^{\gamma_0 \Psi} \partial_y(w + w^s + \varepsilon f(t)\chi'(y))\|_{L^\infty }^2 \sum_{n=0}^{m-1}
	\left\| e^{\Psi} g_n\right\|_{L^2 } ^2
	\\
	&\quad+C_{C_1, \bar{\rho}} \left\langle t\right\rangle ^{\frac{9}{2}-2\kappa}\left\| a\right\|_{H^{m+1}_x}^2 
	\|e^{\gamma_0 \Psi} \partial_y^2(w + w^s + \varepsilon    f(t)\chi'(y))\|_{L^\infty_x L^2_y }^2 \sum_{n=0}^{m-1}
	\left\| e^{\Psi} g_n\right\|_{L^2 } ^2.
\end{align*}
Adding the above two estimates to \eqref{gm11.22}, we can directly obtain
\begin{align}
	\begin{split}
		I_1
		&
		\leq 3C_{a,\infty}\bar\rho \left\| \frac{1}{\sqrt{\rho}}e^{\Psi} \partial_y g_m\right\|_{L^2} \sum_{n=0}^{m-1}   \left\|\frac{1}{\sqrt{\rho}} e^{\Psi}\partial_yg_{n}\right\|_{L^2} +  \frac{1}{3C_1}\left\| \frac{1}{\sqrt{\rho}}e^{\Psi}\partial_yg_{m}\right\|_{L^2}^2
		\\
		&\quad+C \left\langle t\right\rangle ^{4-2\kappa}\left( 1+\left\| a\right\|_{H^{m+1}_x}^2 \right) 
		\left\|e^{\gamma_0 \Psi} \partial_y(w + w^s + \varepsilon    f(t)\chi'(y))\right\|_{L^\infty }^2 \sum_{n=0}^{m}
		\left\| e^{\Psi} g_n\right\|_{L^2 } ^2.
		\\
		&\quad+C \left\langle t\right\rangle ^{\frac{9}{2}-2\kappa}\left( 1+\left\| a\right\|_{H^{m+1}_x}^2 \right) 
		\left\|e^{\gamma_0 \Psi} \partial_y^2(w + w^s + \varepsilon    f(t)\chi'(y))\right\|_{L^\infty_x L^2_y }^2 \sum_{n=0}^{m}
		\left\| e^{\Psi} g_n\right\|_{L^2 } ^2.
	\end{split}
	\label{gm1.2}
\end{align}
	
	\textbf{Dealing with $I_2$ term:}
	Using the integration by parts in the $x$-variable for $I_2$, we can deduce from \eqref{infty1} and \eqref{infty3} in Lemma \ref{infty} that
	\begin{align}
		\begin{split}
			I_2
			&=
			-\frac{1}{2}\int_{\mathbb{R}_+^2} e^{2\Psi}\partial_x(g_m)^2 (u + u^s + \varepsilon f(t)\chi(y))
			\\
			&=
			\frac{1}{2}\int_{\mathbb{R}_+^2} e^{2\Psi}(g_m)^2 \partial_x u 
			\\
			&\leq \left\| \partial_x u\right\|_{L^\infty }\left\| e^{\Psi}g_m\right\|_{L^2 }^2
			\\
			&\leq C_\delta\left\langle t\right\rangle^{\frac12} \left\|e^{\delta\Psi} \partial_x w\right\|_{L^\infty }
			\left\| e^{\Psi}g_m\right\|_{L^2 }^2 
			\\
			&\leq C_\gamma\left\langle t\right\rangle^{\frac34}\left\| e^{\gamma\Psi} \partial_y w\right\|_{H^{1,0} }\left\| e^{\Psi}g_m\right\|_{L^2 }^2 
			\\
			& \leq C_\gamma \left\langle t\right\rangle^{\kappa-\frac{3}{2}}\left( \left\langle t\right\rangle^{{\frac94}-\kappa}\left\| e^{\gamma\Psi} \partial_y w\right\|_{H^{1,0} }\right) \left\| e^{\Psi}g_m\right\|_{L^2 }^2 .
		\end{split}
		\label{gm1.3}
	\end{align}
	
	\textbf{Dealing with $I_3$ term:}
	Using the integration by parts in the $y$-variable for  $I_3$, and using \eqref{ptype2} ($\gamma=1, \zeta=\frac12$) in Lemma \ref{po}, we have
	\begin{align}
		\begin{split}
			I_3
			&\leq 
			4
			\left\| a+1\right\|_{L^\infty_x }^{\frac{3}{2}}
			\left\| e^{\Psi} g_m\right\|_{L^2 }
			\left( \left\|\frac{1}{\sqrt{\rho}} e^{\Psi}\partial_y g_m\right\|_{L^2 }+\frac{1}{2\left\langle t\right\rangle } \left\|\frac{1}{\sqrt{\rho}} e^{\Psi}y g_m\right\|_{L^2 }\right)
			\\
			&\quad \cdot
			\left\| \frac{\partial_y(w + w^s + \varepsilon f(t)\chi'(y))}{w + w^s + \varepsilon f(t)\chi'(y)}\right\|_{L^\infty }
			\\
			&\leq 12 
			C_{\bar{\rho}}
			\left\| e^{\Psi} g_m\right\|_{L^2 }
			\left\| \frac{1}{\sqrt{\rho}}e^{\Psi} \partial_y g_m\right\|_{L^2 }
			\left\| \frac{\partial_y(w + w^s + \varepsilon f(t)\chi'(y))}{w + w^s + \varepsilon f(t)\chi'(y)}\right\|_{L^\infty }
			\\
			&\leq \frac{1}{6C_1}\left\| \frac{1}{\sqrt{\rho}}e^{\Psi} \partial_y g_m\right\|_{L^2 }^2+
			C \left\langle t\right\rangle ^{4-2\kappa}
			\left\|e^{\gamma_0 \Psi} \partial_y(w + w^s + \varepsilon f(t)\chi'(y))\right\|_{L^\infty }^2 
			\left\| e^{\Psi} g_m\right\|_{L^2 } ^2. 
		\end{split}
		\label{gm1.4}
	\end{align}
	
\textbf{Dealing with $I_4$  term :}
Since 
$$v {=} -\partial_x \varphi - \frac{\partial_x a}{a+1} \int_y^\infty \left[ u + u^s - \varepsilon f(t)(1 - \chi(y')) \right] dy',$$ 
we move to estimate $$I_4=\int_{\mathbb{R}_+^2} e^{2\Psi}g_m\partial_x^m u \left\{\partial_x w+\partial_y\left( \frac{\partial_y(w + w^s + \varepsilon f(t)\chi'(y))}{w + w^s + \varepsilon f(t)\chi'(y)} v\right) \right\}$$ involving three terms
\begin{align*}
	I_4^1=\int_{\mathbb{R}_+^2} e^{2\Psi}g_m\partial_x^m u \partial_x w,
\end{align*}
\begin{align*}
	I_4^2=-\int_{\mathbb{R}_+^2} e^{2\Psi}g_m\partial_x^m u \partial_y\left( \frac{\partial_y(w + w^s + \varepsilon f(t)\chi'(y))}{w + w^s + \varepsilon f(t)\chi'(y)} \partial_x \varphi\right), 
\end{align*}
and
\begin{align*}
	I_4^3=-\int_{\mathbb{R}_+^2} e^{2\Psi}g_m\partial_x^m u \partial_y\left( \frac{\partial_y(w + w^s + \varepsilon f(t)\chi'(y))}{w + w^s + \varepsilon f(t)\chi'(y)} \frac{\partial_x a}{a+1} \int_y^\infty \left[ u + u^s - \varepsilon f(t)(1 - \chi(y')) \right] dy'\right).
\end{align*}
For $I_4^1$, we use
\eqref{ptype1} and \eqref{gmin2} to find that
\begin{align}
	\begin{split}
		I_4^1 &\leq  \left\| e^{(1-\delta)\Psi} \partial_x^m u\right\|_{L^2 }
		\left\| e^{\delta\Psi}  \partial_x w\right\|_{L^\infty }\left\| e^{\Psi}  g_m\right\|_{L^2 }
		\\
		&\leq C_\delta \left\langle t\right\rangle^{\frac{1}{2}} \left\| e^{\delta\Psi}  \partial_x w\right\|_{L^\infty }
		\left( 1+\left\langle t\right\rangle^{\frac{5}{2}-\kappa}\left\| e^{\gamma_0\Psi}\partial_y (w + w^s + \varepsilon f(t)\chi'(y))\right\|_{L^\infty }\right) 
		\left\| e^{\Psi}  g_m\right\|_{L^2 }^2.
	\end{split}
	\label{gm1.5}
\end{align}
To estimate $I_4^2$, we recall from the definition of $g_m$ that
\begin{align*}
	&\left\| e^{(1-\delta)\Psi}\frac{\partial_y(w + w^s + \varepsilon f(t)\chi'(y))}{w + w^s + \varepsilon f(t)\chi'(y)}\partial_{x}^{m}u \right\|_{L^2 }\\
	& \leq \left\| e^{(1-\delta)\Psi} g_m \right\|_{L^2 }+\left\| e^{(1-\delta)\Psi} \partial_x^m w\right\|_{L^2 } \\
	&\leq C_{c_0, \gamma_0, \delta} \left(1 + \langle t \rangle^{\frac{5}{2}-\kappa} \|e^{\gamma_0 \Psi} \partial_y(w + w^s + \varepsilon f(t)\chi'(y))\|_{L^\infty }\right) \|e^\Psi g_m\|_{L^2 },
\end{align*}
so, applying \eqref{gmin2} in Lemma \ref{gmw}, we have
\begin{align*}
	I_4^2&=\int_{\mathbb{R}_+^2}\frac{y}{2\left\langle t\right\rangle } e^{2\Psi}g_m\partial_x^m u \left( \frac{\partial_y(w + w^s + \varepsilon f(t)\chi'(y))}{w + w^s + \varepsilon f(t)\chi'(y)} \partial_x \varphi\right)
	+\int_{\mathbb{R}_+^2} e^{2\Psi}\partial_yg_m\partial_x^m u \left( \frac{\partial_y(w + w^s + \varepsilon f(t)\chi'(y))}{w + w^s + \varepsilon f(t)\chi'(y)} \partial_x \varphi\right)\\
	&\quad+
	\int_{\mathbb{R}_+^2} e^{2\Psi}g_m\partial_x^m w \left( \frac{\partial_y(w + w^s + \varepsilon f(t)\chi'(y))}{w + w^s + \varepsilon f(t)\chi'(y)} \partial_x \varphi\right)
	\\
	&\leq \frac{1}{12C_1}\left\|\frac{1}{\sqrt{\rho}}e^\Psi \partial_yg_m\right\|_{L^2 }^2
	+C_{c_0, \gamma_0, \delta} \left(1 + \langle t \rangle^{\frac{5}{2}-\kappa} \|e^{\gamma_0 \Psi} \partial_y(w + w^s + \varepsilon f(t)\chi'(y))\|_{L^\infty }\right)^2 \left\| e^{\delta\Psi}\partial_x \varphi \right\|_{L^\infty }^2 \|e^\Psi g_m\|_{L^2 }^2
	\\
	&\quad+
	C_{c_0, \gamma_0, \delta}
	\left\langle t\right\rangle^{{\frac52}-\kappa}\left\| e^{\gamma_0\Psi}\partial_y (w + w^s + \varepsilon f(t)\chi'(y))\right\|_{L^\infty }
	\left(1 + \langle t \rangle^{\frac{5}{2}-\kappa} \|e^{\gamma_0 \Psi} \partial_y(w + w^s + \varepsilon f(t)\chi'(y))\|_{L^\infty }\right)
	\\
	&\qquad \cdot
	\left\| e^{(1-\delta)\Psi} g_m\right\|_{L^2 }
	\left( \left\langle t\right\rangle^{-\frac{1}{2}}\left\| e^{\delta\Psi}\partial_x \varphi \right\|_{L^\infty }\right) 
	\left\| e^{\Psi} g_m\right\|_{L^2 },
\end{align*}
which yields
\begin{align}
	\begin{split}
		I_4^2 &\leq \frac{1}{12C_1}\left\|\frac{1}{\sqrt{\rho}}e^\Psi \partial_yg_m\right\|_{L^2 }^2\\
		&\quad+ \left(1 + \langle t \rangle^{\frac{5}{2}-\kappa} \|e^{\gamma_0 \Psi} \partial_y(w + w^s + \varepsilon f(t)\chi'(y))\|_{L^\infty }\right)^2
		\\
		&\qquad\cdot 
		\left( \left\langle t\right\rangle^{-\frac{1}{2}}\left\| e^{\delta\Psi}\partial_x \varphi \right\|_{L^\infty }+\left\| e^{\delta\Psi}\partial_x \varphi \right\|_{L^\infty }^2\right)
		\left\| e^{\Psi} g_m\right\|_{L^2 }^2 \\
		&\leq \frac{1}{12C_1}\left\|\frac{1}{\sqrt{\rho}}e^\Psi \partial_yg_m\right\|_{L^2 }^2\\
		&\quad+ C_{c_0, \gamma_0, \delta}\left(1 + \langle t \rangle^{\frac{5}{2}-\kappa} \|e^{\gamma_0 \Psi} \partial_y(w + w^s + \varepsilon f(t)\chi'(y))\|_{L^\infty }\right)^2
		\\
		&\qquad\cdot 
		\left( \left\langle t\right\rangle^{\frac{1}{2}}\left\| e^{\delta\Psi}\partial_x w \right\|_{L^\infty }+\left\langle t\right\rangle^2\left\| e^{\delta\Psi}\partial_x w \right\|_{L^\infty }^2\right)
		\left\| e^{\Psi} g_m\right\|_{L^2 }^2.
	\end{split}
	\label{gm1.6}
\end{align}
For $I_4^3$, we proceed similarly as for $I_4^2$ due to their analogous structure. Thus, it suffices to estimate the following norm
\begin{align*}
	\left\| \frac{\partial_x a}{a+1} \int_y^\infty \left[ u + u^s - \varepsilon f(t)(1 - \chi(y')) \right] dy'\right\|_{L^\infty },
\end{align*}
and then substitute the resulting norm bound in place of $\left\| e^{\delta\Psi}\partial_x \varphi \right\|_{L^\infty }$
appearing in the third line of \eqref{gm1.6}. It follows from \eqref{infty2} that
\begin{align*}
	&\left\| e^{\delta\Psi}\frac{\partial_x a}{a+1} \int_y^\infty \left[ u + u^s - \varepsilon f(t)(1 - \chi(y')) \right] dy'\right\|_{L^\infty }\\
	&\leq C_\delta \left\| \frac{\partial_x a}{a+1} \right\|_{L^\infty_x }
	\left( \left\langle t\right\rangle\left\|  e^{\delta\Psi} w \right\|_{L^{\infty} }
	+ \left\langle t\right\rangle^{\frac{3}{2}}\left\|  e^{\delta\Psi} \partial_y w^s \right\|_{L^{\infty}_y }
	+ C_{\chi}\varepsilon |f(t)|
	\right) ,
\end{align*}
which implies that
\begin{align}
	\begin{split}
		I_4^3 &\leq  \frac{1}{12C_1}\|\frac{1}{\sqrt{\rho}} e^\Psi \partial_yg_m\|_{L^2 }^2
		\\
		&
		\quad+ C_{c_0, \gamma_0, \delta}\left(1 + \langle t \rangle^{\frac{5}{2}-\kappa} \|e^{\gamma_0 \Psi} \partial_y(w + w^s + \varepsilon f(t)\chi'(y))\|_{L^\infty }\right)^2 \left\| \frac{\partial_x a}{a+1} \right\|_{L^\infty_x } 
		\\
		&\qquad\cdot 
		\left(  \left\langle t\right\rangle^{\frac{1}{2}}\left\| e^{\delta\Psi} w \right\|_{L^\infty }+\left\langle t\right\rangle \left\| e^{\delta\Psi}  \partial_y w^s \right\|_{L^\infty_y } +\varepsilon \left\langle t\right\rangle^{-\frac12}|f(t)| \right)
		\left\| e^{\Psi} g_m\right\|_{L^2 }^2
		\\
		&
		\quad+ C_{c_0, \gamma_0, \delta}\left(1 + \langle t \rangle^{\frac{5}{2}-\kappa} \|e^{\gamma_0 \Psi} \partial_y(w + w^s + \varepsilon f(t)\chi'(y))\|_{L^\infty }\right)^2  \left\| \frac{\partial_x a}{a+1} \right\|_{L^\infty_x }^2 
		\\
		&\qquad\cdot 
		\left( \left\langle t\right\rangle^2 \left\| e^{\delta\Psi} w \right\|_{L^\infty }^2+\left\langle t\right\rangle^{3} \left\| e^{\delta\Psi}  \partial_y w^s \right\|_{L^\infty_y }^2+\varepsilon^2 |f(t)|^2 \right) 
		\left\| e^{\Psi} g_m\right\|_{L^2 }^2.
	\end{split}
	\label{gm1.7}
\end{align}
Combining  estimates \eqref{gm1.5}-\eqref{gm1.7} for $I_4^1-I_4^3$, and using the fact that 
\begin{align*}
	\left\| e^{\delta\Psi}  w\right\|_{L^\infty }+	\left\| e^{\delta\Psi}\partial_x w  \right\|_{L^\infty }  \leq 	C_\gamma \left\langle t\right\rangle^\frac{1}{4} \left\| e^{\gamma\Psi} \partial_y w \right\|_{H^{2,0} },
\end{align*}
we estimate $I_4$ as follows
\begin{align}
	\begin{split}
		I_4 &\leq  \frac{1}{6C_1}\left\|\frac{1}{\sqrt{\rho}}e^\Psi \partial_yg_m\right\|_{L^2 }^2 \\
		&\quad+ C_{c_0, \gamma_0, \gamma, \delta}\left(1 + \langle t \rangle^{\frac{5}{2}-\kappa} \|e^{\gamma_0 \Psi} \partial_y(w + w^s + \varepsilon f(t)\chi'(y))\|_{L^\infty }\right)^2\left( 1+ \left\| \frac{\partial_x a}{a+1} \right\|_{L^\infty_x }^2\right) 
		\\
		&\qquad\cdot 
		\left( \left\langle t\right\rangle^{\frac{3}{4}}\left\| e^{\gamma\Psi} \partial_y w \right\|_{H^{2,0} }+\left\langle t\right\rangle^{\frac52}\left\| e^{\gamma\Psi} \partial_y w \right\|_{H^{2, 0} }^2+
		\left\langle t\right\rangle\left\| e^{\gamma_0\Psi} \partial_y w^s \right\|_{L^\infty_y }\right. 
		\\
		&\qquad \quad\left. +\left\langle t\right\rangle^{3}\left\| e^{\gamma_0 \Psi} \partial_y w^s \right\|_{L^\infty_y }^2
		\right)
		\left\| e^{\Psi} g_m\right\|_{L^2 }^2
		\\
		&
		\quad+ C_{c_0, \gamma_0, \delta}\left(1 + \langle t \rangle^{\frac{5}{2}-\kappa} \|e^{\gamma_0 \Psi} \partial_y(w + w^s + \varepsilon f(t)\chi'(y))\|_{L^\infty }\right)^2 
		\\
		&\qquad\cdot 
		\left(   \varepsilon \left\langle t\right\rangle^{-\frac12}\left\| \frac{\partial_x a}{a+1} \right\|_{L^\infty_x } |f(t)| +\varepsilon^2\left\| \frac{\partial_x a}{a+1} \right\|_{L^\infty_x }^2 |f(t)|^2 \right) 
		\left\| e^{\Psi} g_m\right\|_{L^2 }^2.
	\end{split}
	\label{gm1.8}
\end{align}

	\textbf{Estimate of the term $I_5$.}
	By the definition of $F_3$, we have
	\begin{align*}
		\partial_y\left( \frac{F_3}{w+w^s+\varepsilon f(t)\chi'(y)}\right) 
		=
		\frac{\partial_x a}{a+1}
		\left( w+w^s+\varepsilon f(t)\chi'(y) \right).
	\end{align*}
	Then, we deduce that
	\begin{align}
		\begin{split}
			I_5
			&=
			-\int_{\mathbb{R}_+^2}
			e^{2\Psi} g_m \partial_x^m u
			\frac{\partial_x a}{a+1}
			\left( w+w^s+\varepsilon f(t)\chi'(y) \right)
			\\
			&\leq
			C\left\langle t\right\rangle^{-\frac32+\kappa} \left\|\frac{\partial_x a}{a+1}\right\|_{L_x^\infty}
			\left(
			1+\langle t\rangle^{\frac52-\kappa}
			\left\|
			e^{\gamma_0\Psi}
			\partial_y(w + w^s + \varepsilon f(t)\chi'(y))
			\right\|_{L^\infty}
			\right)^2
			\left\| e^\Psi g_m \right\|_{L^2}^2 .
		\end{split}
		\label{gm1.9}
	\end{align}
	
	\textbf{Dealing with $I_6$ term:}
We estimate that $$I_6=-\int_{\mathbb{R}_+^2} \partial_y\left(  e^{2\Psi}g_m(w + w^s + \varepsilon f(t)\chi'(y))\right)  \frac{\left[ u, \partial_{x}^m\right] \partial_x u
	+\left[ w, \partial_x^m\right] v
	+\partial_x^m a\partial_y^2 u^s}{w + w^s + \varepsilon f(t)\chi'(y)}$$ can be handled analogously to $I_4$, provided that the complex terms 
$$\left[ u, \partial_{x}^m\right] \partial_x u
+\left[ w, \partial_x^m\right] v
$$
in $I_6$ are appropriately addressed. According to compressible condition and \eqref{m2} in Lemma \ref{moser}, we know
\begin{align*}
	&\left\| e^{\Psi}(\left[ u, \partial_{x}^m\right] \partial_x u
	+\left[ w, \partial_x^m\right] v)\right\|_{L^2 }\\
	&\lesssim 2\left\| e^{\delta\Psi}\partial_x u \right\|_{L^\infty }
	\left\| e^{(1-\delta)\Psi}\partial_x^m u \right\|_{L^2 }
	+\left\| e^{\delta\Psi}\partial_x w \right\|_{L^\infty }
	\left\| e^{(1-\delta)\Psi}\partial_x^{m-1} v \right\|_{L^2 }
	\\
	&\quad+\left\| e^{\delta\Psi}v \right\|_{L^\infty }
	\left\| e^{(1-\delta)\Psi}\partial_x^{m} w \right\|_{L^2 } .
\end{align*}
Moreover, it follows from the definition of $v$ and Lemma \ref{infty} that
\begin{align*}
	&\left\| e^{(1-\delta)\Psi}\partial_x^{m-1} v \right\|_{L^2 }
	\\
	&\leq
	\left\| e^{(1-\delta)\Psi}\partial_x^{m-1}\left( \partial_x\varphi - \frac{\partial_x a}{a+1} \int_y^\infty \left[ u (t, x, y') + u^s (t, y')- \varepsilon f(t)(1 - \chi(y')) \right] dy'\right)  \right\|_{L^2 }
	\\
	&\lesssim \left( 1+\left\| \frac{\partial_xa}{a+1} \right\|_{H^{m}_x }\right) \left\langle t\right\rangle\left( \left\| e^{(1-\delta)\Psi}\partial_x^{m} w \right\|_{L^2 }+\left\| e^{(1-\delta)\Psi}\partial_x^{m-1} w \right\|_{L^2 }\right),
\end{align*}
then we can deduce that
\begin{align*}
	\begin{split}
		&\left\| e^{\Psi}(\left[ u, \partial_{x}^m\right] \partial_x u
		+\left[ w, \partial_x^m\right] v)\right\|_{L^2 }
		\\
		&\lesssim  \left( 1+\left\| \frac{\partial_xa}{a+1} \right\|_{H^{m}_x }\right) \left\langle t\right\rangle 
		\left\| e^{\delta\Psi}\partial_x w \right\|_{L^\infty }
		\left(1 + \langle t \rangle^{\frac{5}{2}-\kappa} \|e^{\gamma_0 \Psi} \partial_y(w + w^s + \varepsilon f(t)\chi'(y))\|_{L^\infty }\right) \left\| e^{\Psi} g_m\right\|_{L^2 }
		\\
		&\lesssim  \left( 1+\left\| \frac{\partial_xa}{a+1} \right\|_{H^{m}_x }\right) \left\langle t\right\rangle^{\frac{5}{4}} 
		\left\| e^{\delta\Psi}\partial_y w \right\|_{H^{2,0} }
		\left(1 + \langle t \rangle^{\frac{5}{2}-\kappa} \|e^{\gamma_0 \Psi} \partial_y(w + w^s + \varepsilon f(t)\chi'(y))\|_{L^\infty }\right) \left\| e^{\Psi} g_m\right\|_{L^2 },
	\end{split}
\end{align*}
which, together with the previous estimation methods, implies 
	\begin{align}
		\begin{split}
			I_6
			&\lesssim
			\left\| e^{\Psi} \left( \left[ u, \partial_{x}^m\right] \partial_x u
			+\left[ w, \partial_x^m\right] v
			+\partial_x^m a\partial_y w^s\right) \right\|_{L^2 }
			\left\| e^{\Psi} \partial_y g_m\right\|_{L^2 }
			\\
			&\quad+\left\| e^{\Psi}\frac{\partial_y(w + w^s + \varepsilon f(t)\chi'(y))}{w + w^s + \varepsilon f(t)\chi'(y)} \left( \left[ u, \partial_{x}^m\right] \partial_x u
			+\left[ w, \partial_x^m\right] v
			+\partial_x^m a\partial_y w^s\right) \right\|_{L^2 }
			\left\| e^{\Psi}  g_m\right\|_{L^2 }
			\\
			&\leq\frac{1}{6C_1}\left\|\frac{1}{\sqrt{\rho}}e^\Psi \partial_yg_m\right\|_{L^2 }^2
			+C\left( 1+\left\langle t\right\rangle^{\frac32-\kappa} \right) \left\| e^{\Psi}\partial_x^m a\partial_y w^s\right\|_{L^2 }^2
			\\
			&\quad+ C_{c_0, \gamma_0, \delta}\left( 1+\left\| \frac{\partial_xa}{a+1} \right\|_{H^{m}_x }\right) ^2\left(1 + \langle t \rangle^{\frac{5}{2}-\kappa} \|e^{\gamma_0 \Psi} \partial_y(w + w^s + \varepsilon f(t)\chi'(y))\|_{L^\infty }\right)^2
			\\
			&\qquad\cdot 
			\left(\left\langle t\right\rangle^{\frac{5}{2}} 
			\left\| e^{\gamma\Psi}\partial_y w \right\|_{H^{2,0} }^2+\left\langle t\right\rangle^{\frac{3}{4}} 
			\left\| e^{\gamma\Psi}\partial_y w \right\|_{H^{2,0} }+\left\langle t\right\rangle^{-\frac32+\kappa}\right)
			\left\| e^{\Psi} g_m\right\|_{L^2 }^2.
		\end{split}
		\label{gm1.10}
	\end{align}
	
	\textbf{Dealing with $I_7$ term:}
	We directly deal with $I_7$ as follows:
	\begin{align}
		\begin{split}
			I_7 
			&\lesssim   \left\|    \partial_x^m a\right\|_{L^2_x }  \varepsilon|f(t)|\left\|e^{ \Psi} \chi'''(y)\right\|_{L^2_y } \left\| e^{\Psi} g_m\right\|_{L^2 }
			\\
			&\quad+  \left\|    \partial_x^m a\right\|_{L^2_x } \varepsilon|f(t)| \langle t \rangle^{2-\kappa}\left\| e^{\gamma_0 \Psi}\partial_y(w + w^s + \varepsilon f(t)\chi'(y))\right\|_{L^\infty } 
			\left\| e^{ \Psi}\chi''(y)\right\|_{L^2_y } \left\| e^{\Psi} g_m\right\|_{L^2 }
			\\
			&\lesssim
			\left\|    \partial_x^m a\right\|_{L^2_x } \varepsilon|f(t)| \left( 1+\langle t \rangle^{2-\kappa}\left\| e^{\gamma_0 \Psi}\partial_y(w + w^s + \varepsilon f(t)\chi'(y))\right\|_{L^\infty } \right) 
			\left\| e^{\Psi} g_m\right\|_{L^2 }
			\\
			&\lesssim
			\varepsilon^2 \left\langle t\right\rangle |f(t)|^2+  \left( \left\langle t\right\rangle^{-1} \left\|    \partial_x^m a\right\|_{L^2_x }^2+\left\|    \partial_x^m a\right\|_{L^2_x }^2\langle t \rangle^{4-2\kappa}\left\| e^{\gamma_0 \Psi}\partial_y(w + w^s + \varepsilon f(t)\chi'(y))\right\|_{L^\infty } \right)^2 
			\left\| e^{\Psi} g_m\right\|_{L^2 }^2
			.
		\end{split}
		\label{gm1.11}
	\end{align}
	
	\textbf{Dealing with $I_8$ term:}
	It is very easy to find that
	\begin{align*}
		\left\| e^{ \Psi} \partial_x^m F_3 \right\|_{L^2 }	
		&\leq \left\| e^{ \Psi}  F_3 \right\|_{H^{m,0}_{x}}	
		\\
		&\leq
		\left\|  \frac{\partial_x a}{a+1}\right\|_{H^{m}_x}   \left\| e^{(1-\delta) \Psi}  u \right\|_{H^{m,0}_x}  \left\|e^{\delta\Psi} (w + w^s + \varepsilon f(t)\chi'(y)) \right\|_{L^\infty }
		\\
		&\quad
		+\left\|  \frac{\partial_x a}{a+1}\right\|_{H^{m}_x} \left\| e^{(1-\delta) \Psi}  w \right\|_{H^{m,0}_x}  \left\| e^{\delta\Psi}\left( u + u^s - \varepsilon f(t)(1 - \chi(y)) \right)  \right\|_{L^\infty }.
	\end{align*}
	Then directly estimating $I_8$ yields
	\begin{align}
		\begin{split}
			I_8&\lesssim \left\| \frac{\partial_xa}{a+1} \right\|_{H^{m}_x }
			\langle t \rangle^{-\frac{3}{2}+\kappa} 
			\left(1 + \langle t \rangle^{\frac{5}{2}-\kappa} \|e^{\gamma_0 \Psi} \partial_y(w + w^s + \varepsilon f(t)\chi'(y))\|_{L^\infty }\right)^2
			\left\| e^{\Psi} g_m\right\|_{L^2 }^2.
		\end{split}
		\label{gm1.12}
	\end{align}
	
Inserting \eqref{gm1.2}, \eqref{gm1.3}, \eqref{gm1.4}, \eqref{gm1.8}, \eqref{gm1.9}, \eqref{gm1.10}, \eqref{gm1.11}, and \eqref{gm1.12} into \eqref{gm1.1}, and noting that the case $j=0$ is easier since no interaction term appears, we only treat the case $1\le j\le m$. For each $1\le j\le m$, we have
\begin{align*}
	&\frac{1}{2}\frac{d}{dt} \left\| e^{\Psi} g_j\right\|_{L^2}^2
	+\frac{1}{2}\left\|\frac{1}{\sqrt{\rho}} e^{\Psi} \partial_y g_j\right\|_{L^2}^2
	\\
	&\leq
	C \theta_1 \sum_{n=0}^{m}\left\| e^{\Psi} g_n\right\|_{L^2}^2
	+C\left( \left( 1+\left\langle t\right\rangle^{\frac{3}{2}-\kappa}\right)
	\left\| a\right\|_{H^{m+1}_x }^2\left\| e^{\Psi}\partial_y w^s\right\|_{L^2_y }^2
	+\varepsilon^2 \left\langle t\right\rangle |f(t)|^2 \right)
	\\
	&\quad
	+3C_{a,\infty}\bar{\rho}
	\left\|\frac{1}{\sqrt{\rho}} e^{\Psi} \partial_y g_j\right\|_{L^2}
	\sum_{n=1}^{j-1}
	\left\| \frac{1}{\sqrt{\rho}}e^{\Psi}\partial_y g_{n}\right\|_{L^2}
	+\frac{1}{C_1}\left\|\frac{1}{\sqrt{\rho}} e^{\Psi}\partial_y g_{j}\right\|_{L^2}^2.
\end{align*}

Set
$$
A_j(t):=\left\|\frac{1}{\sqrt{\rho}}e^\Psi \partial_y g_j(t)\right\|_{L^2},
\qquad 1\le j\le m.
$$
Summing the above inequalities with respect to $j$ from $1$ to $m$,  we obtain
\begin{align*}
	\frac12\frac{d}{dt}\sum_{j=1}^m \|e^\Psi g_j\|_{L^2}^2
	+\frac12\sum_{j=1}^m A_j^2(t)
	&\le
	\sum_{j=1}^m 3C_{a,\infty}\bar\rho\, A_j(t)\sum_{n=1}^{j-1}A_n(t)
	+\frac1{C_1}\sum_{j=1}^m A_j^2(t)
	+C\theta_1\sum_{j=0}^m \|e^\Psi g_j\|_{L^2}^2
	\\
	&\quad
	+C\left(\left(1+\langle t\rangle^{\frac32-\kappa}\right)
	\left\| a\right\|_{H^{m+1}_x }^2\left\| e^{\Psi}\partial_y w^s\right\|_{L^2_y }^2
	+\left\langle t\right\rangle \right).
\end{align*}

Since
$$
\sum_{j=1}^m A_j\sum_{n=1}^{j-1}A_n
=
\sum_{1\le n<j\le m}A_jA_n,
$$
we infer from Young's inequality that
\begin{align*}
	\sum_{j=1}^m 3C_{a,\infty}\bar\rho\, A_j\sum_{n=1}^{j-1}A_n
	&=
	3C_{a,\infty}\bar\rho\sum_{1\le n<j\le m}A_jA_n
	\\
	&\le
	\frac32 C_{a,\infty}\bar\rho\sum_{1\le n<j\le m}(A_j^2+A_n^2).
\end{align*}
For each fixed $k\in\{1,\dots,m\}$, the quantity $A_k^2$ appears exactly $m-1$ times in the above double sum. Therefore,
$$
\sum_{1\le n<j\le m}(A_j^2+A_n^2)
=
(m-1)\sum_{k=1}^m A_k^2.
$$
Hence,
\begin{align*}
	\sum_{j=1}^m 3C_{a,\infty}\bar\rho\, A_j\sum_{n=1}^{j-1}A_n
	\le
	\frac{3(m-1)}{2}C_{a,\infty}\bar\rho\sum_{k=1}^m A_k^2.
\end{align*}
Substituting this estimate into the previous inequality yields
\begin{align*}
	&\frac{d}{dt}\sum_{j=1}^m \|e^\Psi g_j\|_{L^2}^2
	+\left(
	1-\frac{2}{C_1}-3(m-1)C_{a,\infty}\bar\rho
	\right)\sum_{k=1}^m A_k^2(t) \\
	&\leq
	C\theta_1\sum_{j=0}^m \|e^\Psi g_j\|_{L^2}^2
	+C\left(\left(1+\langle t\rangle^{\frac32-\kappa}\right)
	\left\| a\right\|_{H^{m+1}_x }^2\left\| e^{\Psi}\partial_y w^s\right\|_{L^2_y }^2
	+\varepsilon^2 \left\langle t\right\rangle  |f(t)|^2\right).
\end{align*}

Now choose $C_1>0$ such that
$$
\frac{2}{C_1}\le c(\kappa_0).
$$
By the assumption on $C_{a,\infty}$, we have
\begin{align*}
	3(m-1)C_{a,\infty}\bar\rho
	&\le
	\frac12\Bigl(1-c(\kappa_0)-\kappa_0-(1-4\kappa_0)\bar\rho\Bigr).
\end{align*}
Therefore,
\begin{align*}
	1-\frac{2}{C_1}-3(m-1)C_{a,\infty}\bar\rho
	&\ge
	1-c(\kappa_0)-\frac12\Bigl(1-c(\kappa_0)-\kappa_0-(1-4\kappa_0)\bar\rho\Bigr)
	\\
	&=
	\frac12\Bigl(1-c(\kappa_0)+\kappa_0+(1-4\kappa_0)\bar\rho\Bigr).
\end{align*}
Hence,
\begin{align*}
	1-\frac{2}{C_1}-3(m-1)C_{a,\infty}\bar\rho-\kappa_0
	&\geq
	\frac12\Bigl(1-c(\kappa_0)-\kappa_0+(1-4\kappa_0)\bar\rho\Bigr)
	\\
	&>
	\frac12\Bigl(1-c(\kappa_0)-\kappa_0\Bigr)
>
	\frac12\bigl(1-4\kappa_0\bigr)
	>0,
\end{align*}
where we used $0<c(\kappa_0)<3\kappa_0$ and $0<\kappa_0<\frac14$. Consequently,
\begin{align*}
	\frac{d}{dt}\sum_{j=1}^m \|e^\Psi g_j\|_{L^2}^2
	+\kappa_0\sum_{j=1}^m A_j^2(t)
	&\le
	C\theta_1\sum_{j=0}^m \|e^\Psi g_j\|_{L^2}^2
	\\
	&\quad
	+C\left(\left(1+\langle t\rangle^{\frac32-\kappa}\right)
	\left\| a\right\|_{H^{m+1}_x }^2\left\| e^{\Psi}\partial_y w^s\right\|_{L^2_y }^2
	+\varepsilon^2 \left\langle t\right\rangle  |f(t)|^2\right).
\end{align*}
Combining this with the analogous estimate for $j=0$, we arrive at \eqref{g0}.
\end{proof}

\section{Energy estimate of $G$}\label{s4}
In this section, we derive such a priori decay estimate
of $G$, namely, the following decay-in-time estimate stated in Proposition \eqref{pG} below.
\begin{proposition}\label{pG}
	Fix $0<\kappa_0<\frac14$, and let $c(\kappa_0)$ be a positive constant such that
$
		0<c(\kappa_0)<3\kappa_0.
$
	For any $\kappa\in(\kappa_0,\frac12)$, suppose that the density function $\rho=\rho(t,x)$ satisfies
	\begin{align*}
		1\le \rho(t,x)\le \bar\rho<
		\frac{1-c(\kappa_0)-\kappa_0}{1-4\kappa_0}.
	\end{align*}
	Moreover, choose $C_{a,\infty}$ such that
	\begin{align*}
		\max_{1 \leq i \leq m}{  C_{m}^{i}\left\| {\partial_x^{i} a} \right\|_{L^\infty_x}}=	C_{a,\infty}\le
		\frac{1-c(\kappa_0)-\kappa_0-(1-4\kappa_0)\bar\rho}{6(m-1)\bar\rho}.
	\end{align*}
	 Then there exists a positive constant $C$ such that for any
	$
	0<2\delta\le \gamma_0\le \gamma<1,
	t\le T^{\#},
	$
	one has
	\begin{align}
		\begin{split}
			&\frac{d}{dt}\left( \langle t\rangle^{\frac52-2\kappa}\left\| e^{\Psi}  G\right\|_{H^{5,0} }^2\right)
			+\kappa_0 \sum^{5}_{j=0}
			\left\| \frac{1}{\sqrt{\rho }}\langle t\rangle^{\frac54-\kappa}e^{\Psi} \partial_y \partial_x^j G\right\|_{L^2 }^2
			\\
			&\leq
			C\theta_2 \left\| \langle t\rangle^{\frac54-\kappa} e^{\Psi}  G\right\|_{H^{5,0} }^2
			+\kappa_0 \sum_{n=0}^{m}\big\| e^\Psi  \partial_y g_{n} \big\|_{L^2 }^2
			\\
			&\quad
			+C\langle t\rangle^{\frac52-2\kappa}
			\left( \left\| a\right\|_{H^{6}_{x} }^2
			+\left\| \frac{\partial_x a}{a+1}\right\|_{H^{6}_{x} }^2\right)
			\left( \left\| e^{\Psi}  w^s\right\|_{L^{2}_y }^2
			+\left\| e^{\gamma_0\Psi}  w^s\right\|_{L^{2}_y }^2
			+\varepsilon^2|f(t)|^2\right),
		\end{split}
		\label{G0}
	\end{align}
where
\begin{align}
\begin{split}
\theta_2(t)&= \left\langle t\right\rangle ^{\frac{3}{4}}\left\| e^{\gamma \Psi} \partial_y w\right\|_{H^{2,0} }
+\left\| e^{\gamma\Psi}\partial_x\varphi  \right\|_{L^{\infty} }^2 
+\langle t\rangle^{-\frac{4m(1-\kappa)-31+24\kappa}{2(2m-11)}}
\Big(
\langle t\rangle^{-\frac12}
\Big)^{\frac{2(m-4)}{2m-11}}
\\
&\quad+\langle t\rangle^{-\frac{4m(1-\kappa)-31+24\kappa}{2(2m-11)}}
\Big(
\langle t\rangle^{\frac94-\kappa}
\Big(
\|e^{\gamma\Psi}\partial_y^2 w\|_{L^2}
+\|e^{\gamma_0\Psi}\partial_y^2 w^s\|_{L_y^2}
\Big)
+\varepsilon \langle t\rangle^{2-\kappa}|f(t)|
\Big)^{\frac{2(m-4)}{2m-11}}
\\
&
\quad+
\left\|  \frac{\partial_x a}{a+1}\right\|_{H^{6}_{x} }^2 
\left( \left\langle t\right\rangle^2 \left\| e^{\gamma_0\Psi}(w+w^s+\varepsilon f(t)\chi'(y)) \right\|_{L^{\infty} }^2+
\left\| e^{\gamma\Psi}\partial_x \varphi \right\|_{L^{\infty} }
+
\left\langle t\right\rangle^{\frac{1}{2}}\left\| e^{\gamma\Psi}w \right\|_{L^{\infty} }^2\right) 
%\\
%&\quad
%+
%\varepsilon^2
%\left( 
%\left\| e^{\gamma_0\Psi} \partial_y w^s\right\|_{L^{2} }^2 
%+ |f(t)|^2  \left\| e^{\Psi}\chi''(y)\right\|_{L^2_y }^2 \right)
.
\end{split}
\label{theta2}
\end{align}
\end{proposition}
\begin{proof}
	Notice that
	\begin{align*}
		&v \partial_y G + v \partial_y \left( u^s + \varepsilon f(t)\chi(y) \right) 
		-\frac{1}{2} \langle t \rangle^{-1} v \partial_y (y\varphi) \\
		&=-\partial_x\varphi \left( w+w^s + \varepsilon f(t)\chi'(y) \right)  - \frac{\partial_x a}{a+1} \int_y^\infty \left[ u (t, x, y') + u^s (t, y')- \varepsilon f(t)(1 - \chi(y')) \right] dy'\left( w+w^s + \varepsilon f(t)\chi'(y) \right),
	\end{align*}
	then we can rewrite the equation $\eqref{G uprandtl}_1$ as follows
	\begin{align}
		\begin{split}
			& \partial_t G - (a+1)\partial_y^2 G + \langle t \rangle^{-1} G + (u + u^s + \varepsilon f(t)\chi(y)) \partial_x G -\partial_x\varphi \left( w+w^s + \varepsilon f(t)\chi'(y) \right)  \\ 
			&\quad 
			- \frac{\partial_x a}{a+1} \int_y^\infty \left[ u (t, x, y') + u^s (t, y')- \varepsilon f(t)(1 - \chi(y')) \right] dy'\left( w+w^s + \varepsilon f(t)\chi'(y) \right)\\
			&\quad+ \frac{y}{\langle t \rangle} \int_y^\infty 
			\left( \partial_{y'} (u + u^s + \varepsilon f(t)\chi(y')) \partial_x \varphi \right) dy' 
			\\ 
			&\quad - a \partial_y^2 (u + u^s) -a \varepsilon f(t)\chi''(y) - \frac{y}{2\langle t \rangle}(F_1 + F_2) = 0.
		\end{split}
		\label{G1}
	\end{align}	
	Applying the operator $\partial^5_x $ in \eqref{G1}, and 
	multiplying the resulting equation  by  $e^{2\Psi}   \partial_x^5 G$, we have
	\begin{align}
		\begin{split}
			& \int_{\mathbb{R}_+^2} e^{2\Psi}   \partial_x^5 G \left( \partial_t-  \partial_y^2+\langle t\rangle^{-1} \right)\partial_x^5 G 
			\\
			&
			=
			\int_{\mathbb{R}_+^2} e^{2\Psi}   \partial_x^5 G  [\partial_x^5, (a+1) ]\partial_y^2 G 
			\\
			&\quad
			- \int_{\mathbb{R}_+^2} e^{2\Psi}   \partial_x^5 G \partial_x^5\left\lbrace (u + u^s + \varepsilon f(t)\chi(y)) \partial_x G\right\rbrace 
			\\
			&\quad +\int_{\mathbb{R}_+^2} e^{2\Psi}   \partial_x^5 G \partial_x^5\left\lbrace\partial_x\varphi \left( w+w^s + \varepsilon f(t)\chi'(y) \right) \right\rbrace  \\ 
			&\quad 
			+\int_{\mathbb{R}_+^2} e^{2\Psi}   \partial_x^5 G \partial_x^5\left\lbrace \frac{\partial_x a}{a+1} \int_y^\infty \left[ u (t, x, y') + u^s (t, y')- \varepsilon f(t)(1 - \chi(y')) \right] dy'\left( w+w^s + \varepsilon f(t)\chi'(y) \right)\right\rbrace \\
			&\quad- \int_{\mathbb{R}_+^2} e^{2\Psi}   \partial_x^5 G  \frac{y}{\langle t \rangle} \int_y^\infty \partial_x^5
			\left( \partial_{y'} (u + u^s + \varepsilon f(t)\chi(y')) \partial_x \varphi \right) dy' 
			\\ 
			&\quad +\int_{\mathbb{R}_+^2} e^{2\Psi}   \partial_x^5 G \partial_x^5 (a \partial_y^2   u^s) 
			\\ 
			&\quad
			+\int_{\mathbb{R}_+^2} e^{2\Psi}   \partial_x^5 G (\partial_x^5 a) \varepsilon f(t)\chi''(y) 
			\\ 
			&\quad
			+ \frac{y}{2\langle t \rangle}\int_{\mathbb{R}_+^2} e^{2\Psi}   \partial_x^5 G  \partial_x^5(F_1 + F_2)\\
			&=\sum_{i=1}^{8}J_{i}.
		\end{split}
		\label{G2}
	\end{align}	
	For the term on the left-hand side of the equality \eqref{G2}, by using Lemma \eqref{heat pro}, we can directly obtain 
	\begin{align}
		\begin{split}
			\frac{1}{2}\frac{d}{dt}\left\| e^{\Psi} \partial_x^5 G\right\|_{L^2 }^2
			+\frac{1}{2}
			\left\| \frac{1}{\sqrt{\rho }}e^{\Psi} \partial_y\partial_x^5 G\right\|_{L^2 }^2
			+\frac{1}{\left\langle t\right\rangle }
			\left\| e^{\Psi} \partial_x^5 G\right\|_{L^2 }^2
			\lesssim   \sum_{i=1}^{8}J_{i}.
		\end{split}
		\label{G3}
	\end{align}

	Now, we deal with the right-hand side of \eqref{G3} term by term as follows.
	
	\textbf{Dealing with $J_1$  term :} 
	Since $m\ge m_0>\frac{17-24\kappa}{2(1-2\kappa)}>5$, we have $C_5^i\le C_m^i$ for $1\le i\le 5$. Hence, by the definition of $C_{a,\infty}$ and the bound $1\le \rho\le \bar\rho$, we infer that 
	\begin{align}
		\begin{split}
			|J_1|
			&\le
			\sum_{i=1}^5 C_5^i \|\partial_x^i a\|_{L_x^\infty}
			\left\|e^\Psi \partial_y\partial_x^5 G\right\|_{L^2}
			\left\|e^\Psi \partial_y\partial_x^{5-i} G\right\|_{L^2} \\
			&\le
			\bar\rho\, C_{a,\infty}\, \left\|\frac1{\sqrt\rho}e^\Psi \partial_y\partial_x^5 G(t)\right\|_{L^2} \sum_{j=0}^4 \left\|\frac1{\sqrt\rho}e^\Psi \partial_y\partial_x^j G(t)\right\|_{L^2} \\
			&\le
			\frac52 \bar\rho\, C_{a,\infty}\sum_{j=0}^5 \left\|\frac1{\sqrt\rho}e^\Psi \partial_y\partial_x^j G(t)\right\|_{L^2}^2 \\
			&\le
			6(m-1)\bar\rho\, C_{a,\infty}\sum_{j=0}^5
			\left\|\frac1{\sqrt\rho}e^\Psi \partial_y\partial_x^j G\right\|_{L^2}^2.
		\end{split}
		\label{G40}
	\end{align}

	\textbf{Dealing with $J_2$  term :}
	Integrating by parts in the $x$-variable, and using \eqref{m2} in Lemma \ref{moser},
	we have
	\begin{align}
		\begin{split}
			J_2&=- \int_{\mathbb{R}_+^2} e^{2\Psi}  \partial_x^5 G \partial_x^5\left\lbrace (u + u^s + \varepsilon f(t)\chi(y)) \partial_x G\right\rbrace   \\
			&=\int_{\mathbb{R}_+^2} e^{2\Psi}  (\partial_x^5 G)^2 \partial_xu  
			-\int_{\mathbb{R}_+^2} e^{2\Psi}  \partial_x^5 G\left[ \partial_x^{j}, 
			(u + u^s + \varepsilon f(t)\chi(y))\right] \partial_x G \\
			& \lesssim  \left(  \left\| \partial_x u\right\|_{L^\infty }
			\left\| e^{\Psi}G\right\|_{H^{{5},0} }
			+\left\| e^{\delta \Psi}\partial_x G\right\|_{L^\infty }\left\|  e^{(1-\delta)\Psi}u\right\|_{H^{{5},0} }
			\right)  \left\| e^{\Psi}G\right\|_{H^{{5},0} }\\
			& \lesssim  \left(  \left\| \partial_x u\right\|_{L^\infty }
			+\left\| e^{\delta \Psi}\partial_x G\right\|_{L^\infty }
			\right) \left\| e^{\Psi}G\right\|_{H^{{5},0} }^2
			\\
			&\lesssim \left\langle t\right\rangle ^{\frac{3}{4}}\left\| e^{\gamma \Psi} \partial_y w\right\|_{H^{2,0} } 
			\left\| e^{\Psi}G\right\|_{H^{{5},0} }^2 ,
		\end{split}
		\label{G4}
	\end{align}
	where,  in the last step, we have used the following fact  by means of the definition of $G$ and $\varphi(t, x, y)=-\int_{y}^\infty u(t, x, y')dy'$
	\begin{align*}
		\left\| e^{\delta \Psi}\partial_x G\right\|_{L^\infty } 
		&\leq
		\left\| e^{\delta \Psi}\partial_x u\right\|_{L^\infty }
		+\left\| \frac{y}{2\left\langle t\right\rangle }e^{\delta \Psi} \int^\infty_y e^{-2\delta\Psi}\cdot(e^{2\delta\Psi}\partial_x u)dy' \right\|_{L^\infty }\\
		&\lesssim 
		\left\| e^{\delta \Psi}\partial_x u\right\|_{L^\infty }
		+\left\| \frac{y}{2\left\langle t\right\rangle }e^{\delta \Psi}\int^\infty_y e^{-2\delta\Psi(y')}dy'\right\|_{L^\infty_y }\sup\limits_{ y'\geq y}\left| e^{2\delta\Psi}\partial_x u\right|
		\\
		&\lesssim \left\| e^{2\delta \Psi}\partial_x u\right\|_{L^\infty }\leq C_\gamma \left\langle t\right\rangle ^{\frac{3}{4}}\left\| e^{\gamma \Psi}\partial_y w\right\|_{H^{2,0} }
		.
	\end{align*}
	
	\textbf{Dealing with $J_3$  and $J_5$ terms :} To estimate $J_3$, we divide $J_3$ into the following two parts.

	\begin{align*}
		|J_3|&=\left| \int_{\mathbb{R}_+^2} e^{2\Psi} \partial_x^5 G \partial_x^6 \varphi  \left(w+ w^s + \varepsilon f(t)\chi'(y) \right)\right| 
		+\left| \int_{\mathbb{R}_+^2} e^{2\Psi} \partial_x^5 G [\partial_x^5, w]\partial_x \varphi\right|  \\
		&\leq     \left\| e^{(1-\delta)\Psi}\partial_x^{6}\varphi\right\|_{L^{2} } \left\| e^{2\delta\Psi}(w+w^s+ \varepsilon f(t)\chi'(y))\right\|_{L^{\infty} } \left\| e^{\Psi}G\right\|_{H^{{5},0} }\\
		&\quad  +   \left\| e^{(1+\delta)\Psi}\left[ \partial_x^{5}, w\right]  \partial_x\varphi\right\|_{L^{2} } \left\| e^{\Psi}G\right\|_{H^{{5},0} }.
		%+
		%\left\| e^{(1+\delta)\Psi} \partial_x\varphi(w+ w^s + \varepsilon f(t)\chi'(y))\right\|_{H^{{5}-1} }
		%\left\| e^{\Psi}G\right\|_{H^{{5},0} }
	\end{align*}
	To handle the first term, define
	\begin{align}
		X(t,x,y)
		:= \Bigg\{
		\langle t\rangle^{\frac94-\kappa}
		\Big(
		\|e^{\gamma\Psi}\partial_y^2 w\|_{L^2}
		+\|e^{\gamma_0\Psi}\partial_y^2 w^s\|_{L_y^2}
		\Big)
		+\varepsilon \langle t\rangle^{2-\kappa}|f(t)|
		\Bigg\}.
		\label{XT}
	\end{align}
	Using the Gagliardo--Nirenberg inequality in Lemma \ref{GN} and Lemmas \ref{uG}--\ref{infty}, we have
	\begin{align*}
		\big\| e^{(1-\delta)\Psi} \partial_x^{6} \varphi \big\|_{L^2}
		&\lesssim
		\big\| e^{(1-\delta)\Psi} \varphi \big\|_{H^{5,0}}^{1-\frac{1}{m-5}}
		\big\| e^{(1-\delta)\Psi} \partial_x^m \varphi \big\|_{L^2}^{\frac{1}{m-5}}
		\\
		&\lesssim
		\langle t \rangle^{\frac12\left(1-\frac{1}{m-5}\right)}
		\big\| e^\Psi G \big\|_{H^{5,0}}^{\,1-\frac{1}{m-5}}
		\langle t \rangle^{\frac{1}{m-5}}
		\big\| e^\Psi \partial_x^m w \big\|_{L^2}^{\frac{1}{m-5}}
		\\
		&\lesssim
		\langle t \rangle^{\frac12+\frac{1}{m-5}}
		\left(
		1+\langle t\rangle^{\frac52-\kappa}
		\big\| e^{\gamma_0\Psi}\partial_y\big(w+w^s+\varepsilon f(t)\chi'(y)\big)\big\|_{L^\infty}
		\right)^{\frac1{m-5}}
		\\
		&\qquad\qquad\qquad
		\times
		\big\| e^\Psi G \big\|_{H^{5,0}}^{\,1-\frac{1}{m-5}}
		\big\| e^\Psi \partial_y g_m \big\|_{L^2}^{\frac{1}{m-5}}
		\\
		&\lesssim
		\langle t \rangle^{\frac12+\frac{3}{2(m-5)}}
		\Big( \langle t\rangle^{-\frac12}+ X(t,x,y) \Big)^{\frac1{m-5}}
		\big\| e^\Psi G \big\|_{H^{5,0}}^{\,1-\frac{1}{m-5}}
		\big\| e^\Psi \partial_y g_m \big\|_{L^2}^{\frac{1}{m-5}}.
	\end{align*}
	Moreover,
	\begin{align*}
		&\left\| e^{2\delta\Psi}(w+w^s+\varepsilon f(t)\chi'(y))\right\|_{L^{\infty}}
		\\
		&\le
		\langle t \rangle^{-\frac{3}{2} +\kappa}
		\Bigg\{
		\langle t\rangle^{\frac94-\kappa}
		\Big(
		\|e^{\gamma\Psi}\partial_y^2 w\|_{L^2}
		+\|e^{\gamma_0\Psi}\partial_y^2 w^s\|_{L_y^2}
		\Big)
		+\langle t\rangle^{-\frac14}\varepsilon \langle t\rangle^{2-\kappa}|f(t)|
		\Bigg\}
		\\
		&\leq
		\langle t \rangle^{-\frac{3}{2} +\kappa}X(t,x,y)
\leq
		\langle t \rangle^{-\frac{3}{2} +\kappa}
		\Big( \langle t\rangle^{-\frac12}+ X(t,x,y) \Big).
	\end{align*}
	Thus, we arrive at
	\begin{align}
		\begin{split}
			&\big\| e^{(1-\delta)\Psi}\partial_x^{6}\varphi\big\|_{L^{2}}
			\big\| e^{2\delta\Psi}(w+w^s+\varepsilon f(t)\chi'(y))\big\|_{L^{\infty}}
			\big\| e^{\Psi}G\big\|_{H^{5,0}}
			\\
			&\lesssim
			\langle t\rangle^{-1+\kappa+\frac{3}{2(m-5)}}
			\Big( \langle t\rangle^{-\frac12}+ X(t,x,y) \Big)^{1+\frac{1}{m-5}}
			\big\| e^\Psi G \big\|_{H^{5,0}}^{\,2-\frac{1}{m-5}}
			\big\| e^\Psi \partial_y g_m \big\|_{L^2}^{\frac{1}{m-5}}
			\\
			&\le
			\frac14\kappa_0 \langle t\rangle^{-\frac52+2\kappa}
			\big\| e^\Psi \partial_y g_m \big\|_{L^2}^{2}
			+ C_{\kappa,\kappa_0}\,
			\langle t\rangle^{-\frac{4m(1-\kappa)-31+24\kappa}{2(2m-11)}}
			\Big( \langle t\rangle^{-\frac12}+ X(t,x,y) \Big)^{\frac{2(m-4)}{2m-11}}
			\big\| e^\Psi G \big\|_{H^{5,0}}^{2}.
		\end{split}
		\label{G5}
	\end{align}
	And for the remaining one term, it follows from \eqref{m2} in Lemma \ref{moser} and Lemma \ref{uG} that
	\begin{align}
		&  \left\| e^{(1+\delta)\Psi}\left[ \partial_x^{5}, w\right]  \partial_x\varphi\right\|_{L^{2} } \left\| e^{\Psi}G\right\|_{H^{{5},0} }\nonumber\\
		&\lesssim   \left( \left\| e^{2\delta\Psi}w  \right\|_{L^{\infty} }
		\left\| e^{(1-\delta)\Psi}\varphi\right\|_{H^{{5},0} }
		+\left\| e^{2\delta\Psi}\partial_x\varphi  \right\|_{L^{\infty} }
		\left\| e^{(1-\delta)\Psi}w\right\|_{H^{{5},0} }\right)  \cdot  \left\| e^{\Psi}G\right\|_{H^{{5},0} }
		\nonumber\\
		&\lesssim \frac{1}{12 C_2} \sum^{5}_{j=0} \left\| \frac{1}{\sqrt{\rho }}e^{\Psi} \partial_y \partial_x^j G\right\|_{L^2 }^2+\left(  \left\langle t\right\rangle^{\frac{1}{2}}  \left\| e^{\gamma\Psi}w  \right\|_{L^{\infty} }
+
		\left\| e^{\gamma\Psi}\partial_x\varphi  \right\|_{L^{\infty} }^2\right) 
		\left\| e^{\Psi} G\right\|_{H^{{5},0} }^2.
		\label{G6}
	\end{align}
	By direct computations, we find that
	\begin{align}
		|J_5| &=
		\left| \int_{\mathbb{R}_+^2} e^{2\Psi}  \partial_x^5 G  \frac{y}{\langle t \rangle} \int_y^\infty \partial_x^5
		\left(  (w + w^s + \varepsilon f(t)\chi'(y')) \partial_x \varphi \right) dy'\right| 
		\nonumber\\
		&\leq 
		\left\| \left( e^{-\frac{\delta}{2}\Psi} \frac{y}{\langle t \rangle}\right) e^{(1+\delta)\Psi} \left( \int_y^\infty e^{-2(1+\delta)\Psi} dy' \right)^{\frac{1}{2}} \right\|_{L_y^\infty   }  \left\| e^{-\frac{\delta}{2}\Psi} \right\|_{L_y^2 } 
		\nonumber\\
		&\quad\cdot\left\| e^{(1+\delta)\Psi}\partial_x^5\left(  (w + w^s + \varepsilon f(t)\chi'(y)) \partial_x \varphi \right) \right\|_{L^{2} }  \left\| e^\Psi G \right\|_{H^{5,0}}
		\nonumber\\
		&\lesssim 
		\left\| e^{(1+\delta)\Psi}\partial_x^5\left(  (w + w^s + \varepsilon f(t)\chi'(y)) \partial_x \varphi \right) \right\|_{L^{2} }  \left\| e^\Psi G \right\|_{H^{5,0}},
		\label{G7}
	\end{align}
	which implies that $J_5$ has the same upper bound as $J_3$. 
	Collecting the above three inequalities \eqref{G5}-\eqref{G7}, we find
	\begin{align}
		\begin{split}
			&|J_3 |+ |J_5 |
			\\
			&\leq \frac{1}{2}\kappa_0\left\langle t\right\rangle^{-\frac{5}{2}+2k}\big\| e^\Psi  \partial_y g_{m} \big\|_{L^2 }^2
			+\frac{1}{6C_2 } \sum^{5}_{j=0} \left\| \frac{1}{\sqrt{\rho }}e^{\Psi} \partial_y \partial_x^j G\right\|_{L^2 }^2
			\\
			&
			\quad
			+C\left(\langle t\rangle^{-\frac{4m(1-\kappa)-29+24\kappa}{2(2m-11)}}
			\Big(
			\langle t\rangle^{-\frac12}+
			X(t,x, y)
			\Big)^{\frac{2(m-4)}{2m-11}}
			+ \left\langle t\right\rangle^{\frac{1}{2}} \left\| e^{\gamma\Psi}w  \right\|_{L^{\infty} }
			+
			\left\| e^{\gamma\Psi}\partial_x\varphi  \right\|_{L^{\infty} }^2\right) \left\| e^{\Psi}G\right\|_{H^{{5},0} }^2
			.
		\end{split}
		\label{G8}
	\end{align}

	\textbf{Dealing with $J_4$  term :}
	Since the order of the terms in $J_4$ does not exceed $5$, it is very easy to infer
	\begin{align}
		\begin{split}
			|J_4|
			& \lesssim \int_{\mathbb{R}_+^2} e^{2\Psi}  \partial_x^5 G \partial_x^5\left( \frac{\partial_x a}{a+1}\right) \\
			&\quad\cdot\int_y^\infty \left[ u (t, x, y') + u^s (t, y')- \varepsilon f(t)(1 - \chi(y')) \right] dy'\left( w+w^s + \varepsilon f(t)\chi'(y) \right)
			\\
			&\quad
			+
			\int_{\mathbb{R}_+^2} e^{2\Psi}  \partial_x^5 G  \frac{\partial_x a}{a+1} \int_y^\infty \left[ u (t, x, y') + u^s (t, y')- \varepsilon f(t)(1 - \chi(y')) \right] dy' \partial_x^5 w
			\\
			&\quad
			+
			\sum_{i=1}^{5}\left| \int_{\mathbb{R}_+^2} e^{2\Psi}  \partial_x^5 G \partial_x^{5-i}\left(  \frac{\partial_x a}{a+1} \left( w+w^s + \varepsilon f(t)\chi'(y) \right)\right) \partial_x^i \varphi\right| 
			\\
			&
			\lesssim 
			\left\langle t\right\rangle  \left\|  \frac{\partial_x a}{a+1}\right\|_{H^{6}_{x} } 
			\left\| e^{2\delta\Psi}(w+w^s+\varepsilon f(t)\chi'(y)) \right\|_{L^{\infty} }
			\left\| e^{(1-2\delta)\Psi}(w+w^s+\varepsilon f(t)\chi'(y)) \right\|_{L^{2} }
			\left\| e^{\Psi}G\right\|_{H^{{5},0} }
			\\
			&\quad+
			\left\langle t\right\rangle  \left\| \frac{\partial_xa}{a+1}\right\|_{L^\infty_x } 
			\left\| e^{2\delta\Psi}(w+w^s+\varepsilon f(t)\chi''(y)) \right\|_{L^{\infty} } 
			\left\| e^{\Psi}\partial_y G\right\|_{H^{{5},0} }
			\left\| e^{\Psi}G\right\|_{H^{{5},0} }
			\\
			&\quad
			+
			\left\|  \frac{\partial_x a}{a+1}\right\|_{H^{6}_{x} } 
			\left( \left\| e^{2\delta\Psi}\partial_x \varphi \right\|_{L^{\infty} } \left\| e^{(1-\delta)\Psi}\partial_x^4 w \right\|_{L^{2} }+\left\| e^{2\delta\Psi}w \right\|_{L^{\infty} } \left\| e^{(1-\delta)\Psi} \partial_x^5 \varphi \right\|_{L^{2} }\right)
			\\
			&\qquad \cdot
			\left\| e^{\Psi}G\right\|_{H^{{5},0} }
			\\
			&
			\leq \frac{1}{12C_2 } \sum^{5}_{j=0} \left\| \frac{1}{\sqrt{\rho }}e^{\Psi} \partial_y \partial_x^j G\right\|_{L^2 }^2
			+
			C
			\left\|  \frac{\partial_x a}{a+1}\right\|_{H^{6}_{x} }^2\left( 
			\left\| e^{\gamma_0\Psi}  w^s\right\|_{L^{2}_y }^2 
			+ \varepsilon^2|f(t)|^2   \right)
			\\
			&\qquad +
			\left( \left\langle t\right\rangle^2 \left\| e^{\gamma_0\Psi}(w+w^s+\varepsilon f(t)\chi'(y)) \right\|_{L^{\infty} }^2
			+
			\left\| e^{\gamma_0\Psi}\partial_x \varphi \right\|_{L^{\infty} }
			+
			\left\langle t\right\rangle^{\frac{1}{2}}\left\| e^{\gamma_0\Psi}w \right\|_{L^{\infty} }^2\right) \left\| e^{\Psi}G\right\|_{H^{{5},0} }^2
			.
		\end{split}
		\label{G9}
	\end{align}

	\textbf{Dealing with $J_6$ and $J_7$ terms :} By the similar arguments, we have
	\begin{align}
		\begin{split}
			&J_6+J_7
			\\
			&= \int_{\mathbb{R}_+^2} e^{2\Psi}  \partial_x^5 G \partial_x^5 a \partial_y w^s
			+\int_{\mathbb{R}_+^2} e^{2\Psi}  \partial_x^5 G (\partial_x^5 a) \varepsilon f(t)\chi''(y) 
			\\
			&
			\leq
			\frac{1}{12C_2 } \sum^{5}_{j=0} \left\| \frac{1}{\sqrt{\rho }}e^{\Psi} \partial_y \partial_x^j G\right\|_{L^2 }^2
			+
			C
			\left\| a\right\|_{H^{6}_{x} }^2\left( 
			\left\| e^{\Psi}  w^s\right\|_{L^{2}_y }^2 
			+ \varepsilon^2|f(t)|^2   \right).
		\end{split}
		\label{G10}
	\end{align}

	\textbf{Dealing with $J_8$  term :}
	As for $J_8$, we have
	\begin{align*}
		J_8& \leq  
		\left\| e^{\Psi}   \partial_x^5F_1 \right\|_{L^{{2}} }
		\left\| \frac{y}{2\langle t \rangle}e^{\Psi}G\right\|_{H^{{5},0} }
		+ \left\| \frac{y}{2\langle t \rangle}e^{\Psi}   \partial_x^5F_2 \right\|_{L^{{2}} }
		\left\| e^{\Psi}G\right\|_{H^{{5},0} }  \\
		&\leq J_8^1+J_8^2,
	\end{align*}
	where
	\begin{align*}
		\begin{split}
			F_1 &{=} a\partial_y u^s + a\varepsilon f(t)\chi'(y), \\
			F_2 &{=} -\frac{\partial_x a}{a+1} \int_y^\infty \left[ u + u^s - \varepsilon f(t)(1 - \chi(y')) \right] \left[ u + u^s + \varepsilon f(t)\chi(y') \right] dy' \\
			&\quad + \frac{\partial_x a}{a+1} \int_y^\infty \left[ u + u^s - \varepsilon f(t)(1 - \chi(y')) \right] dy' \left[ u + u^s + \varepsilon f(t)\chi(y') \right].
		\end{split}
	\end{align*}
	Since
	\begin{align*}
		\left\| e^{\Psi}   \partial_x^5F_1 \right\|_{L^{{2}} }
		\leq
		\left\| e^{\Psi}   \partial_x^5(a w^s + a\varepsilon f(t)\chi'(y)) \right\|_{L^{{2}} },
	\end{align*}
	by the exactly same procedure with $J_6+J_7$, $J_8^1$ has the same bound
	\begin{align*}
		\begin{split}
			J_8^1
			\leq
			\frac{1}{12C_2 } \sum^{5}_{j=0}
			\left\| \frac{1}{\sqrt{\rho }}\langle t\rangle^{\frac54-\kappa}e^{\Psi} \partial_y \partial_x^j G\right\|_{L^2 }^2
			+
			C
			\left\| a\right\|_{H^{6}_{x} }^2\left( 
			\left\| e^{\Psi}  w^s\right\|_{L^{2}_y }^2 
			+ \varepsilon^2|f(t)|^2   \right).
		\end{split}
	\end{align*}
	On the other hand, according to the definition of $F_2$, we know
	\begin{align*}
		\begin{split}
			\partial_y F_2&{=} \frac{\partial_x a}{a+1} \left[ u + u^s - \varepsilon f(t)(1 - \chi(y')) \right] \left[ u + u^s + \varepsilon f(t)\chi(y') \right]  \\
			&\quad - \frac{\partial_x a}{a+1}  \left[ u + u^s - \varepsilon f(t)(1 - \chi(y')) \right]  \left[ u + u^s + \varepsilon f(t)\chi(y') \right]
			\\
			&\quad - \frac{\partial_x a}{a+1} \int_y^\infty \left[ u + u^s - \varepsilon f(t)(1 - \chi(y')) \right] dy' \left[ w + w^s + \varepsilon f(t)\chi'(y') \right]\\
			&=- \frac{\partial_x a}{a+1} \int_y^\infty \left[ u + u^s - \varepsilon f(t)(1 - \chi(y')) \right] dy' \left[ w + w^s + \varepsilon f(t)\chi'(y') \right].
		\end{split}
	\end{align*}
	In a similar way with $J_4$,  using Lemma \ref{uG}, we have
	\begin{align*}
		J_8^2 &\leq \frac{1}{12C_2 } \sum^{5}_{j=0} \left\| \frac{1}{\sqrt{\rho }}e^{\Psi} \partial_y \partial_x^j G\right\|_{L^2 }^2
		+
		C
		\left\|  \frac{\partial_x a}{a+1}\right\|_{H^{6}_{x} }^2\left( 
		\left\| e^{\gamma_0\Psi}  w^s\right\|_{L^{2}_y }^2 
		+ \varepsilon^2|f(t)|^2   \right)
		\\
		&\qquad +
		\left( \left\langle t\right\rangle^2 \left\| e^{\gamma_0\Psi}(w+w^s+\varepsilon f(t)\chi'(y)) \right\|_{L^{\infty} }^2
		+
		\left\| e^{\gamma_0\Psi}\partial_x \varphi \right\|_{L^{\infty} }
		+
		\left\langle t\right\rangle^{\frac{1}{2}}\left\| e^{\gamma_0\Psi}w \right\|_{L^{\infty} }^2\right) \left\| e^{\Psi}G\right\|_{H^{{5},0} }^2
		,
	\end{align*}
	which gives immediately
	\begin{align}
		\begin{split}
			J_8
			&
			\leq \frac{1}{6C_2 } \sum^{5}_{j=0} \left\| \frac{1}{\sqrt{\rho }}e^{\Psi} \partial_y \partial_x^j G\right\|_{L^2 }^2
			+
			C
			\left( 	\left\|  a\right\|_{H^{6}_{x} }^2+\left\|  \frac{\partial_x a}{a+1}\right\|_{H^{6}_{x} }^2\right) \left( \left\| e^{\Psi}  w^s\right\|_{L^{2}_y }^2 +
			\left\| e^{\gamma_0\Psi}  w^s\right\|_{L^{2}_y }^2 
			+ \varepsilon^2|f(t)|^2   \right)
			\\
			&\qquad +
			\left( \left\langle t\right\rangle^2 \left\| e^{\gamma_0\Psi}(w+w^s+\varepsilon f(t)\chi'(y)) \right\|_{L^{\infty} }^2
			+
			\left\| e^{\gamma_0\Psi}\partial_x \varphi \right\|_{L^{\infty} }
			+
			\left\langle t\right\rangle^{\frac{1}{2}}\left\| e^{\gamma_0\Psi}w \right\|_{L^{\infty} }^2\right) \left\| e^{\Psi}G\right\|_{H^{{5},0} }^2 .
		\end{split}
		\label{G11}
	\end{align}
	
	Substituting the above estimates \eqref{G40}, \eqref{G4}, and \eqref{G8}--\eqref{G11} for $J_1-J_8$ into \eqref{G3}, we obtain
	\begin{align*}
		\frac{d}{dt}\left\| e^{\Psi} \partial_x^5 G\right\|_{L^2 }^2
		&+
		\left\|\frac{1}{\sqrt{\rho }} e^{\Psi} \partial_y\partial_x^5 G\right\|_{L^2 }^2
		+\frac{2}{\langle t\rangle}
		\left\| e^{\Psi} \partial_x^5 G\right\|_{L^2 }^2 \\
		&\le
		C\theta_2 \left\| e^{\Psi}  G\right\|_{H^{5,0} }^2
		+\kappa_0\langle t\rangle^{-\frac52+2\kappa}\big\| e^\Psi  \partial_y g_{m} \big\|_{L^2 }^2
		\\
		&\quad
		+C
		\left( \left\|  a\right\|_{H^{6}_{x} }^2+\left\|  \frac{\partial_x a}{a+1}\right\|_{H^{6}_{x} }^2\right)
		\left( \left\| e^{\Psi}  w^s\right\|_{L^{2}_y }^2
		+\left\| e^{\gamma_0\Psi}  w^s\right\|_{L^{2}_y }^2
		+ \varepsilon^2|f(t)|^2 \right) \\
		&\quad
		+\frac{1}{C_2} \sum^{5}_{j=0} \left\| \frac{1}{\sqrt{\rho }}e^{\Psi} \partial_y \partial_x^j G\right\|_{L^2 }^2
		+6(m-1)\bar\rho\, C_{a,\infty}
		\sum_{j=0}^5
		\left\|\frac1{\sqrt\rho}e^\Psi \partial_y\partial_x^j G\right\|_{L^2}^2.
	\end{align*}
	
	Applying $\partial_x^j$ with $0\le j\le 4$ to \eqref{G}, one obtains the same estimate. Summing over $0\le j\le 5$ and multiplying by $\langle t\rangle^{\frac52-2\kappa}$, we infer that
	\begin{align*}
		&\frac{d}{dt}\left( \langle t\rangle^{\frac52-2\kappa}\left\| e^{\Psi}  G\right\|_{H^{5,0} }^2\right) \\
		&\quad+
		\Bigl(
		1-\frac6{C_2}-\bar\rho(1-4\kappa)_+
		-6(m-1)\bar\rho\, C_{a,\infty}
		\Bigr)
		\sum^{5}_{j=0}
		\left\| \frac{1}{\sqrt{\rho }}\langle t\rangle^{\frac54-\kappa}e^{\Psi} \partial_y \partial_x^j G\right\|_{L^2 }^2 \\
		&\le
		C\theta_2 \left\| \langle t\rangle^{\frac54-\kappa} e^{\Psi}  G\right\|_{H^{5,0} }^2
		+\kappa_0 \sum_{n=0}^{m}\big\| e^\Psi  \partial_y g_{n} \big\|_{L^2 }^2
		\\
		&\quad
		+C\langle t\rangle^{\frac52-2\kappa}
		\left( \left\|  a\right\|_{H^{6}_{x} }^2+\left\|  \frac{\partial_x a}{a+1}\right\|_{H^{6}_{x} }^2\right)
		\left( \left\| e^{\Psi}  w^s\right\|_{L^{2}_y }^2
		+\left\| e^{\gamma_0\Psi}  w^s\right\|_{L^{2}_y }^2
		+ \varepsilon^2|f(t)|^2 \right).
	\end{align*}
	
	Now choose
	$
	C_2:=\frac{6}{c(\kappa_0)}.
	$
	Under the assumptions
	$$
	1\le \rho(t,x)\le \bar\rho<
	\frac{1-c(\kappa_0)-\kappa_0}{1-4\kappa_0},
	$$
	and
	$$
	C_{a,\infty}\le
	\frac{1-c(\kappa_0)-\kappa_0-(1-4\kappa_0)\bar\rho}{6(m-1)\bar\rho},
	$$
	we have
	\begin{align*}
		\frac6{C_2}+\bar\rho(1-4\kappa)_+
		+6(m-1)\bar\rho\, C_{a,\infty}
		&\le
		c(\kappa_0)+(1-4\kappa_0)\bar\rho
		+6(m-1)\bar\rho\, C_{a,\infty} \\
		&\le
		c(\kappa_0)+(1-4\kappa_0)\bar\rho
		+1-c(\kappa_0)-\kappa_0-(1-4\kappa_0)\bar\rho \\
		&=
		1-\kappa_0.
	\end{align*}
	Therefore,
	\begin{align*}
		1-\frac6{C_2}-\bar\rho(1-4\kappa)_+
		-6(m-1)\bar\rho\, C_{a,\infty}
		\ge \kappa_0.
	\end{align*}
	Consequently,
	\begin{align*}
		&\frac{d}{dt}\left( \langle t\rangle^{\frac52-2\kappa}\left\| e^{\Psi}  G\right\|_{H^{5,0} }^2\right)
		+\kappa_0
		\sum^{5}_{j=0}
		\left\| \frac{1}{\sqrt{\rho }}\langle t\rangle^{\frac54-\kappa}e^{\Psi} \partial_y \partial_x^j G\right\|_{L^2 }^2 \\
		&\le
		C\theta_2 \left\| \langle t\rangle^{\frac54-\kappa} e^{\Psi}  G\right\|_{H^{5,0} }^2
		+\kappa_0\sum_{n=0}^{m}\big\| e^\Psi  \partial_y g_{n} \big\|_{L^2 }^2
		\\
		&\quad
		+C\langle t\rangle^{\frac52-2\kappa}
		\left( \left\|  a\right\|_{H^{6}_{x} }^2+\left\|  \frac{\partial_x a}{a+1}\right\|_{H^{6}_{x} }^2\right)
		\left( \left\| e^{\Psi}  w^s\right\|_{L^{2}_y }^2
		+\left\| e^{\gamma_0\Psi}  w^s\right\|_{L^{2}_y }^2
		+ \varepsilon^2|f(t)|^2 \right),
	\end{align*}
	which completes the proof of Proposition \ref{pG}.
\end{proof}

\section{Energy estimate of $w$} \label{s5}
In this section, we are going to give the weighted  energy estimate of $w$
and to obtain its decay-in-time estimate. To be specific, we have the following proposition.
\begin{proposition}\label{pw}
There exists a positive constant $C$ such that, for any
$\kappa \in (0, \frac{1}{2})$, $0<2\delta \leq \gamma_0 \leq \gamma <1$,and for all $t \leq T^{\#}$, one has
\begin{align}
	\begin{split}
	& \frac{d}{dt} \left\| \left\langle t\right\rangle^{\frac{7}{4}-\kappa}e^{\gamma\Psi} w\right\|_{H^{4,0} }^2
+  \sum_{j=0}^{4}\left\| \frac{1}{\sqrt{\rho}}\left\langle t\right\rangle^{\frac{7}{4}-\kappa} e^{\gamma \Psi} \partial_y \partial_x^j w \right\|_{L^2 }^2  
\\
&
\leq C\theta_3\left\| \left\langle t\right\rangle^{\frac{7}{4}-\kappa}e^{\gamma\Psi} w\right\|_{H^{4,0} }^2
+C\left\langle t\right\rangle^2 \left\|  w + w^s + \varepsilon f(t)\chi'(y)  \right\|_{L^\infty }^2
\left( 1+\left\| \frac{\partial_x a}{a+1}\right\|_{H^{4,0}_x }^2\right) 
\left\| \left\langle t\right\rangle^{\frac{5}{4}-\kappa}e^{ \Psi}G
\right\|_{H^{5,0} }^2
\\
& \quad+C\left\langle t\right\rangle^{\frac{7}{2}-2\kappa}\left\| a\right\|_{H^{5,0}_x }^2\left(\langle t\rangle \|e^{\gamma\Psi}\partial_y^2 w\|_{H^{3,0}}^2+ \|e^{\gamma \Psi} \partial_y w^s \|_{L^2_y }^2+\varepsilon^2|f(t)|^2\right) 
+\left( {\frac72-2\kappa}\right)\left\| \langle t\rangle^{\frac54-\kappa}e^{\Psi} \partial_y  G\right\|_{H^{5,0} }^2,
	\end{split}
	\label{w0}
\end{align}
where
\begin{align}
	\begin{split}
		\theta_3&= \left\langle t\right\rangle^{\frac{9}{4}-2\kappa} \left\| e^{\gamma\Psi}\partial_y^2 w\right\|_{H^{2,0}}^2
		+\langle t \rangle^{\frac{1}{2}} \| e^{\gamma\Psi}\partial_x w \|_{L^\infty}+\left\| \partial_x \varphi \right\|_{L^\infty }^2
		\\
		&\quad+\left\langle t\right\rangle \left\| \frac{\partial_x a}{a+1}\right\|_{H^{4,0}_x }^2\left\|  u  + u^s - \varepsilon f(t)(1 - \chi(y)) \right\|_{L^\infty }^2
		\\
		&
		\quad
		+\left\| \frac{\partial_x a}{a+1}\right\|_{H^{5,0}_x }
	\left\langle t\right\rangle\left\|e^{\gamma_0 \Psi}\partial_y \left(  w + w^s +\varepsilon f(t) \chi(y) \right) \right\|_{L^\infty }
	 .
	\end{split}
	\label{theta3}
\end{align}
\end{proposition}
\begin{proof}
Since $a=a(t,x)$ is independent of $y$, we rewrite $\eqref{vprandtl}_1$ as
\begin{align}
	\begin{split}
	&\partial_t w + (u + u^s + \varepsilon f(t)\chi(y)) \partial_x w + v \partial_y \left( w + w^s + \varepsilon f(t)\chi'(y) \right) \\
\\
&\quad
-\partial_y\bigl((a+1)\partial_y w + a\partial_y w^s\bigr)
-a\varepsilon f(t)\chi'''(y) = F_3.
	\end{split}
	\label{vprandtl-flux}
\end{align}
	
	Taking the $H^{4,0}$ inner product of \eqref{vprandtl-flux} with $e^{2\gamma\Psi}w$, we obtain
\begin{align}
	\begin{split}
		0
		&=
		\left\langle e^{\gamma\Psi}\left(\partial_t w-\partial_y\bigl((a+1)\partial_y w+a\partial_y w^s\bigr)\right),e^{\gamma\Psi}w\right\rangle_{H^{4,0}}
		\\
		&\quad
		+
		\left\langle e^{\gamma\Psi}(u + u^s + \varepsilon f(t)\chi(y)) \partial_x w ,e^{\gamma\Psi}w\right\rangle_{H^{4,0}}
		\\
		&\quad
		+
		\left\langle e^{\gamma\Psi}v \partial_y \left( w + w^s + \varepsilon f(t)\chi'(y) \right) ,e^{\gamma\Psi}w\right\rangle_{H^{4,0}}
		\\
		&\quad
		-
		\left\langle e^{\gamma\Psi}(a\varepsilon f(t)\chi'''(y)+F_3),e^{\gamma\Psi}w\right\rangle_{H^{4,0}}
		\\
		&=: K_1+K_2+K_3+ K_4.
	\end{split}
	\label{w1}
\end{align}
Then, we estimate $ K_1$, $K_2$, $K_3$ and $ K_4$ term by term.

For the first term $ K_1$, we have
\begin{align*}
 K_1
	&=
	\sum_{j\leq 4}\int_{\mathbb R_+^2} e^{2\gamma\Psi}\partial_t\partial_x^j w\,\partial_x^j w
	\\
	&\quad
	-
	\sum_{j\leq 4}\int_{\mathbb R_+^2}
	e^{2\gamma\Psi}
	\partial_y\Bigl(
	(a+1)\partial_y\partial_x^j w
	+[\partial_x^j,a+1]\partial_y w
	+\partial_x^j a\,\partial_y w^s
	\Bigr)\partial_x^j w .
\end{align*}
Since
$$
\Bigl(
(a+1)\partial_y\partial_x^j w
+[\partial_x^j,a+1]\partial_y w
+\partial_x^j a\,\partial_y w^s
\Bigr)\Big|_{y=0}
=
\partial_x^j\Bigl((a+1)\partial_y w+a\partial_y w^s\Bigr)\Big|_{y=0}
=0,
$$
the boundary term vanishes after integrating by parts in the $y$-variable. Hence, 
\begin{align*}
 K_1
	&\ge
	\frac12\frac{d}{dt}\|e^{\gamma\Psi}w\|_{H^{4,0}}^2
	+\left(1-\frac{\gamma}{2}\right)
	\sum_{j=0}^4
	\left\|
	\frac1{\sqrt\rho}e^{\gamma\Psi}\partial_y\partial_x^j w
	\right\|_{L^2}^2
	\\
	&\quad
	-\sum_{j\le4}
	\int_{\mathbb R_+^2}
	e^{2\gamma\Psi}
	\Bigl(
	[\partial_x^j,a+1]\partial_y w
	+\partial_x^j a\,\partial_y w^s
	\Bigr)
	\Bigl(
	2\gamma\partial_y\Psi\,\partial_x^j w+\partial_y\partial_x^j w
	\Bigr),
\end{align*}
by Corollary \ref{heat cor}. Therefore,
\begin{align}
	\begin{split}
		&\frac12\frac{d}{dt}\|e^{\gamma\Psi}w\|_{H^{4,0}}^2
		+\left(1-\frac{\gamma}{2}\right)
		\sum_{j=0}^4
		\left\|
		\frac1{\sqrt\rho}e^{\gamma\Psi}\partial_y\partial_x^j w
		\right\|_{L^2}^2
		\\
		&\le
		 K_1
		+C_3\gamma
		\sum_{j=0}^{4}\left\|\frac{1}{\sqrt{\rho}} e^{\gamma \Psi} \partial_y \partial_x^j w \right\|_{L^2}^2
		\\
		&\quad
		+C\langle t\rangle \|a\|_{H_x^{5,0}}^2\|e^{\gamma\Psi}\partial_y^2 w\|_{H^{3,0}}^2
		+C_\gamma \|a\|_{H_x^{5,0}}^2\|e^{\gamma\Psi}\partial_y w^s\|_{L_y^2}^2 .
	\end{split}
	\label{w2}
\end{align}

For the second term $K_2$, along the same line as $J_4$ in \eqref{G4}, we can deduce that 
\begin{align}
	\begin{split}
		|K_2|&=\left| \left\langle e^{{\gamma}\Psi}(u + u^s + \varepsilon f(t)\chi(y)) \partial_x w ,e^{{\gamma}\Psi}w\right\rangle_{H^{4,0} }\right| 
		\\
		&\leq \left( \| \partial_x u \|_{L^\infty } \| e^{\gamma \Psi} w \|_{H^{4,0} } + \| \partial_x w \|_{L^\infty } \| e^{\gamma \Psi} u \|_{H^{4,0} }\right)  \| e^{\gamma \Psi} w \|_{H^{4,0} } \\
		&\lesssim  \left(\| \partial_x u \|_{L^\infty }+ \langle t \rangle^{\frac{1}{2}} \| \partial_x w \|_{L^\infty}  \right) \| e^{\gamma \Psi} w \|_{H^{4,0} }^2
		\\
		&\lesssim  \langle t \rangle^{\frac{1}{2}} \|e^{\gamma\Psi} \partial_x w \|_{L^\infty}  \| e^{\gamma \Psi} w \|_{H^{4,0} }^2.
	\end{split}
	\label{w3}
\end{align}
Next, to establish the estimate of $K_3$, by the integration by parts, we separate $K_3$ into two parts
\begin{align*}
	K_3&=\left\langle e^{{\gamma}\Psi}v \partial_y \left( w + w^s + \varepsilon f(t)\chi'(y) \right) ,e^{{\gamma}\Psi}w\right\rangle_{H^{4,0} }
	\\
	&=\left\langle   \left( w + w^s + \varepsilon f(t)\chi'(y) \right) , \partial_y \left( e^{{\gamma}\Psi}\partial_x\varphi ~e^{{\gamma}\Psi}w\right) \right\rangle_{H^{4,0} }\\
	&\quad+\left\langle   \left( w + w^s + \varepsilon f(t)\chi'(y) \right) , \partial_y \left( e^{{\gamma}\Psi}\left(  \frac{\partial_x a}{a+1} \int_y^\infty \left[ u  + u^s - \varepsilon f(t)(1 - \chi(y')) \right] dy'\right) e^{{\gamma}\Psi}w\right) \right\rangle_{H^{4,0} }
	\\
	&
	=K_3^1+K_3^2.
\end{align*}
For $K_3^1$, using Lemma  \ref{uG}--\ref{moser}, we have
\begin{align*}
	|K_3^1|
	&\lesssim \left\| e^{\gamma \Psi} \left( w + w^s + \varepsilon f(t)\chi'(y) \right) \partial_x \varphi \right\|_{H^{4,0} } \left\| e^{\gamma \Psi} \partial_y w \right\|_{H^{4,0} } \\
	&\quad+ \left\| e^{\gamma \Psi} \left( w + w^s + \varepsilon f(t)\chi'(y) \right) \partial_x u \right\|_{H^{4,0} } \left\| e^{\gamma \Psi} w \right\|_{H^{4,0} } \\
	&\leq  C_3\gamma   \sum_{j=0}^{4}\left\|\frac{1}{\sqrt{\rho}} e^{\gamma \Psi} \partial_y \partial_x^j w \right\|_{L^2}^2 
	+C\left\langle t\right\rangle \left\|  w + w^s + \varepsilon f(t)\chi'(y)  \right\|_{L^\infty }^2\left\| e^{ \Psi}G \right\|_{H^{5,0} }^2
	\\
	&\quad+ C \left(\left\| \partial_x \varphi \right\|_{L^\infty }^2+ \left\| \partial_x u \right\|_{L^\infty } \right) \left\| e^{\gamma \Psi} w \right\|_{H^{4,0} }^2 .
\end{align*}
A similar derivation of the above inequality yields that
\begin{align*}
	|K_3^2| &\leq  C_3\gamma   \sum_{j=0}^{4}\left\|\frac{1}{\sqrt{\rho}} e^{\gamma \Psi} \partial_y \partial_x^j w \right\|_{L^2}^2 
	+C\left\langle t\right\rangle 
	\left\|  w + w^s + \varepsilon f(t)\chi'(y)  \right\|_{L^\infty }^2
	\left\| \frac{\partial_x a}{a+1}\right\|_{H^{4,0}_x }^2
	\left\| e^{ \Psi}G \right\|_{H^{4,0} }^2
	\\
	&\quad+ C \left(\left\langle t\right\rangle \left\|  u  + u^s - \varepsilon f(t)(1 - \chi(y)) \right\|_{L^\infty }^2+ \left\| u  + u^s - \varepsilon f(t)(1 - \chi(y)) \right\|_{L^\infty } \right) 
	\\
	&\qquad\cdot
	\left\| \frac{\partial_x a}{a+1}\right\|_{H^{4,0}_x }^2
	\left\| e^{\gamma \Psi} w \right\|_{H^{4,0} }^2. 
\end{align*}
Combining the above two estimates, 
it holds
\begin{align}
	\begin{split}
		|K_3|
		&\leq  2C_3\gamma   \sum_{j=0}^{4}\left\|\frac{1}{\sqrt{\rho}} e^{\gamma \Psi} \partial_y \partial_x^j w \right\|_{L^2}^2  
		\\
		&\quad
		+C\left\langle t\right\rangle \left\|  w + w^s + \varepsilon f(t)\chi'(y)  \right\|_{L^\infty }^2
		\left( 1+\left\| \frac{\partial_x a}{a+1}\right\|_{H^{4,0}_x }^2\right) 
		\left\| e^{ \Psi}G
		\right\|_{H^{5,0} }^2
		\\
		&\quad+ C \left(\left\| \partial_x \varphi \right\|_{L^\infty }^2+ \left\| \partial_x u \right\|_{L^\infty } \right) \left\| e^{\gamma \Psi} w \right\|_{H^{4,0} }^2 
		\\
		&\quad+ C \left(\left\langle t\right\rangle \left\|  u  + u^s - \varepsilon f(t)(1 - \chi(y)) \right\|_{L^\infty }^2+ \left\| u  + u^s - \varepsilon f(t)(1 - \chi(y)) \right\|_{L^\infty } \right) 
		\\
		&\qquad\cdot
		\left\| \frac{\partial_x a}{a+1}\right\|_{H^{4,0}_x }^2
		\left\| e^{\gamma \Psi} w \right\|_{H^{4,0} }^2. 
	\end{split}
	\label{w4}
\end{align}Finally, for last term $K_4$, we arrive at 
\begin{align}
	\begin{split}
		|K_4|&=\left| \left\langle e^{{\gamma}\Psi}( a\varepsilon f(t)\chi'''(y)+F_3) ,e^{{\gamma}\Psi}w\right\rangle\right| _{H^{4,0} }\\
		& \leq  C_3\gamma       \sum_{j=0}^{4}\left\|\frac{1}{\sqrt{\rho}} e^{\gamma \Psi} \partial_y \partial_x^j w \right\|_{L^2}^2  
		+	C\left\| a\right\|_{H^{5,0}_x }^2\varepsilon^2|f(t)|^2\|e^{\gamma \Psi} \chi''(y)\|_{L^2_y }^2
		\\ 
		&\quad+ C  \left\| \frac{\partial_x a}{a+1}\right\|_{H^{5,0}_x }
	 \left\langle t\right\rangle\left\|e^{\gamma \Psi}\partial_y \left(  w + w^s +\varepsilon f(t) \chi(y) \right) \right\|_{L^\infty }  \left\| e^{\gamma \Psi} w \right\|_{H^{4,0} }^2
,
	\end{split}
	\label{w5}
\end{align}
where we have used the following fact
\begin{align*}
	\left\| e^{\gamma \Psi} F_3 \right\|_{H^{4,0} }	
	&=\left\| e^{\gamma \Psi} \frac{\partial_x a}{a+1} \left[ u + u^s - \varepsilon f(t)(1 - \chi(y)) \right] (w + w^s + \varepsilon f(t)\chi'(y))\right\|_{H^{4,0} }
	\\
	&\leq C\left\| \frac{\partial_x a}{a+1}\right\|_{H^{5,0}_x }
	\left\|e^{\gamma_0 \Psi} \left[ u + u^s - \varepsilon f(t)(1 - \chi(y)) \right]\right\|_{L^\infty } \left\| e^{\gamma \Psi} w \right\|_{H^{4,0} }
	\\
	&\quad
	+C\left\| \frac{\partial_x a}{a+1}\right\|_{H^{5,0}_x }
	\left\|e^{\gamma_0 \Psi} \left[ w + w^s +\varepsilon f(t) \chi(y) \right]\right\|_{L^\infty } \left\| e^{\gamma \Psi} u \right\|_{H^{4,0} }
	\\
	&\leq C  \left\| \frac{\partial_x a}{a+1}\right\|_{H^{5,0}_x }
	\left\langle t\right\rangle\left\|e^{\gamma_0 \Psi}\partial_y \left(  w + w^s +\varepsilon f(t) \chi(y) \right) \right\|_{L^\infty }  \left\| e^{\gamma \Psi} w \right\|_{H^{4,0} }.
\end{align*}
Consequently, summing up all the estimates \eqref{w2}-\eqref{w5} for $K_1-K_4$ together, we can obtain
\begin{align}
	\begin{split}
		& \frac{1}{2}\frac{d}{dt} \left\| e^{\gamma\Psi} w\right\|_{H^{4,0} }^2
		+ \left( 1-\frac{1}{2}\gamma-4C_3\gamma \right)      \sum_{j=0}^{4}\left\|\frac{1}{\sqrt{\rho}} e^{\gamma \Psi} \partial_y \partial_x^j w \right\|_{L^2}^2  
		\\
		&
		\leq C\theta_3\left\| e^{\gamma\Psi} w\right\|_{H^{4,0} }^2
		+C\left\langle t\right\rangle \left\|  w + w^s + \varepsilon f(t)\chi'(y)  \right\|_{L^\infty }^2
		\left( 1+\left\| \frac{\partial_x a}{a+1}\right\|_{H^{4,0}_x }^2\right) 
		\left\| e^{ \Psi}G
		\right\|_{H^{5,0} }^2
		\\
		& \quad+C\left\| a\right\|_{H^{5,0}_x }^2\left( \langle t\rangle \|e^{\gamma\Psi}\partial_y^2 w\|_{H^{3,0}}^2+\|e^{\gamma \Psi} \partial_y w^s \|_{L^2_y }^2+\varepsilon^2|f(t)|^2\right) 
		,
	\end{split}
	\label{w6}
\end{align}
where $\theta_3$ is defined in \eqref{theta3}.
Choosing $\left( 1-\frac{1}{2}\gamma-4C_3\gamma \right) =\frac{1}{2}$ in \eqref{w6}, and
multiplying the  equation \eqref{w6} by $\left\langle t\right\rangle^{\frac{7}{2}-2\kappa}$,  we have
\begin{align*}
	& \frac{d}{dt} \left\| \left\langle t\right\rangle^{\frac{7}{4}-\kappa}e^{\gamma\Psi} w\right\|_{H^{4,0} }^2
	+  \sum_{j=0}^{4}\left\| \frac{1}{\sqrt{\rho}}\left\langle t\right\rangle^{\frac{7}{4}-\kappa} e^{\gamma \Psi} \partial_y \partial_x^j w \right\|_{L^2 }^2  
	\\
	&
	\leq C\theta_3\left\| \left\langle t\right\rangle^{\frac{7}{4}-\kappa}e^{\gamma\Psi} w\right\|_{H^{4,0} }^2
	+C\left\langle t\right\rangle^2 \left\|  w + w^s + \varepsilon f(t)\chi'(y)  \right\|_{L^\infty }^2
	\left( 1+\left\| \frac{\partial_x a}{a+1}\right\|_{H^{4,0}_x }^2\right) 
	\left\| \left\langle t\right\rangle^{\frac{5}{4}-\kappa}e^{ \Psi}G
	\right\|_{H^{5,0} }^2
	\\
	& \quad+C\left\langle t\right\rangle^{\frac{7}{2}-2\kappa}\left\| a\right\|_{H^{5,0}_x }^2\left(\langle t\rangle \|e^{\gamma\Psi}\partial_y^2 w\|_{H^{3,0}}^2+ \|e^{\gamma \Psi} \partial_y w^s \|_{L^2_y }^2+\varepsilon^2|f(t)|^2\right) 
	+\left( {\frac72-2\kappa}\right)\left\| \langle t\rangle^{\frac54-\kappa}e^{\Psi} \partial_y  G\right\|_{H^{5,0} }^2 .
\end{align*}
This completes the proof of Proposition \ref{pw}.

\end{proof}

\section{Energy estimate of $\partial_y w$}\label{s6}

	In this section, we establish the weighted energy estimates for the normal
derivative of the vorticity. Since the boundary condition of $w$ is not
homogeneous for $\partial_y w$ itself, we estimate instead the
boundary-compatible quantity
$
\partial_y w+\frac{a}{a+1}\partial_y w^s .
$
Indeed, since $\chi''(0)=0$, the boundary condition in \eqref{vprandtl} gives
\begin{align}
	\left(
	\partial_y w+\frac{a}{a+1}\partial_y w^s
	\right)\Big|_{y=0}
	=0 .
\end{align}
Consequently, for every $0\leq j\leq 3$,
\begin{align}
	\partial_x^j
	\left(
	\partial_y w+\frac{a}{a+1}\partial_y w^s
	\right)\Big|_{y=0}
	=0 .
\end{align}

\begin{proposition}\label{pyw}
	There exists a positive constant $C$ such that, for any
	$\kappa \in (0, \frac{1}{2})$, $0<2\delta \leq \gamma_0 \leq \gamma <1$,and for all $t \leq T^{\#}$, one has
\begin{align}
	\begin{split}
		&\frac{d}{dt}\left\| \left\langle t\right\rangle^{\frac{9}{4}-\kappa}e^{\gamma\Psi}
		\left( \partial_y w+\frac{a}{a+1}\partial_y w^s\right) \right\|_{H^{3,0} }^2
		+\sum_{j=0}^{3}\left\| \frac{1}{\sqrt{\rho}} e^{\gamma\Psi}
		\partial_y \partial_x^j
		\left( \partial_y w+\frac{a}{a+1}\partial_y w^s\right)\right\|_{L^2 }^2
		\\
		&\leq C\theta_4
		\left\| \left\langle t\right\rangle^{\frac{9}{4}-\kappa}e^{\gamma\Psi}
		\left( \partial_y w+\frac{a}{a+1}\partial_y w^s\right) \right\|_{H^{3,0} }^2
		+C\|a\|_{H^{4}_x}^2
		\left\| \left\langle t\right\rangle^{\frac{9}{4}-\kappa}e^{\gamma\Psi}\partial_y^2 w^s\right\| _{L^2_y}^2
		\\
		&\quad
		+ C\|a\|_{H^{4}_x}^2 \left\langle t\right\rangle^{\frac{9}{2}-2\kappa}
		\left\| e^{\gamma\Psi} \partial_y
		\left( \partial_y w+\frac{a}{a+1}\partial_y w^s\right)\right\|_{H^{3,0} }^2
		\\
		&\quad+
		C \left\langle t\right\rangle^3
		\left( 1+\left\| \frac{\partial_x a}{a+1}\right\|_{H^{3}_x}^2\right)
		\|e^{\gamma_0 \Psi}\partial_y(w + w^s + \varepsilon f(t)\chi'(y))\|_{L^\infty }^2
		\|\left\langle t\right\rangle^{\frac{7}{4}-\kappa}e^{{\gamma}\Psi} w\|_{H^{4,0} }^2
		\\
		&\quad+
		C \left\langle t\right\rangle^{-1+2\kappa}
		\left( 1+\left\| \frac{\partial_x a}{a+1}\right\|_{H^{3}_x}^2\right)
		\left( 1+\| \left\langle t\right\rangle^{\frac{7}{4}-\kappa}e^{\gamma\Psi}w \|_{H^{4,0} }^2\right)
		\left\| \left\langle t\right\rangle^{\frac{7}{4}-\kappa}e^{{\gamma}\Psi} \partial_y w\right\| _{H^{4,0} }^2
		\\
		&\quad
		+C \varepsilon^2 \left\langle t\right\rangle^{\frac{9}{2}-2\kappa}
		\left( 1+\left\|
		a
		\right\|_{H^{4}_x}^2\right)^2 \left\|
		a
		\right\|_{H^{4}_x}^2
		\left(|f(t)|^2+\left\|e^{ \gamma \Psi} 
		\partial_y^2 w^s
		\right\|_{L^2_y}^2\right) \left( 1+\left\|e^{ \gamma \Psi} 
		w^s
		\right\|_{L^2_y}^2+|f(t)|^2+|f'(t)|^2\right),
	\end{split}
	\label{yw0}
\end{align}
where
\begin{align}
	\theta_4
	&=
	\left\langle t \right\rangle^2
	\|w + w^s + \varepsilon f(t)\chi'(y)\|_{L^\infty }^2
	+\left\langle t \right\rangle^2
	\|e^{{\gamma}\Psi}\partial_x \partial_y w\|_{L^\infty }^2
	+\left\langle t \right\rangle^2
	\|e^{{\gamma}\Psi}\partial_x w\|_{L^\infty }^2 .
	\label{theta4}
\end{align}
\end{proposition}
	\begin{proof}
			For simplicity, in the proof we set
		\begin{align*}
			q
			:=
			\partial_y w+\frac{a}{a+1}\partial_y w^s
			.
		\end{align*}
		Thus the first equation of \eqref{vprandtl} can be rewritten as
		\begin{align}
			\begin{split}
				&\partial_t w
				+
				(u+u^s+\varepsilon f(t)\chi(y))\partial_x w
				+
				v\partial_y(w+w^s+\varepsilon f(t)\chi'(y))
				\\
				&\quad
				-
				(a+1)\partial_y q
				-
				a\varepsilon f(t)\chi'''(y)
				=
				F_3 .
			\end{split}
			\label{newpyw6.1}
		\end{align}
		Applying $\partial_y$ to \eqref{newpyw6.1}, 
		we obtain
		\begin{align}
			\begin{split}
				&
				\partial_t
				q
				-
				(a+1)\partial_y^2 q
				\\
				&=
				-\partial_y
				\left\{
				(u+u^s+\varepsilon f(t)\chi(y))\partial_xw
				\right\}
				-
				\partial_y
				\left\{
				v\partial_y(w+w^s+\varepsilon f(t)\chi'(y))
				\right\}
				\\
				&\quad
				+
				\partial_yF_3
				+
				a\varepsilon f(t)\chi^{(4)}(y)
				+
				\partial_t
				\left(
				\frac{a}{a+1}\partial_yw^s
				\right).
			\end{split}
			\label{pyw-eq1}
		\end{align}
		Taking the $H^{3,0}$ inner product of \eqref{pyw-eq1} with
		$
		 q
		$
		weighted by $e^{2\gamma\Psi}$, we get
				\begin{align}
					\begin{split}
							&-\left\langle e^{\gamma\Psi}\left( \partial_t -(a+1)\partial_y^2\right)  q ,e^{\gamma\Psi} q\right\rangle_{H^{3,0} }
						\\
						&=
						\left\langle e^{\gamma\Psi}\partial_y\left\lbrace (u + u^s + \varepsilon f(t)\chi(y)) \partial_x w\right\rbrace  ,e^{\gamma\Psi} q\right\rangle_{H^{3,0} }
						\\
						&\quad
						+
						\left\langle e^{\gamma\Psi}\partial_y \left\lbrace v \partial_y \left( w + w^s + \varepsilon f(t)\chi'(y) \right)\right\rbrace  ,e^{\gamma\Psi}  q\right\rangle_{H^{3,0} }
						\\
						&\quad
						-\left\langle e^{\gamma\Psi}\partial_y\left( a\varepsilon f(t)\chi^{(3)}(y)+ F_3+\partial_t
						\left(
						\frac{a}{a+1}w^s
						\right)\right)  ,e^{\gamma\Psi} q\right\rangle_{H^{3,0} } =\sum_{i=1}^{3} L_i.
					\end{split}
						\label{newpyw6.2}
		\end{align}
 By Corollary \ref{heat cor} and a direct integration by parts in the normal variable, we have
\begin{align}
	\begin{split}
		&\frac12 \frac{d}{dt}\|e^{\gamma\Psi}q\|_{H^{3,0}}^2
		+\left(1-\frac{\gamma}{2}\right)\sum_{j=0}^{3}
		\left\|
		\frac1{\sqrt{\rho}}e^{\gamma\Psi}\partial_y \partial_x^j q
		\right\|_{L^2}^2
		\\
		&\leq
		\sum_{j=1}^{3}\int_{\mathbb R_+^2}
		e^{2\gamma\Psi}[\partial_x^j,a+1]\partial_y^2 q\,\partial_x^j q\,dxdy
		+\sum_{i=1}^{3}|L_i|
		.
	\end{split}
	\label{yw1-new}
\end{align}

	We first estimate the  integral term in \eqref{yw1-new}.
By integration by parts and Young's inequality,
\begin{align}
	\begin{split}
		&\sum_{j=1}^{3}\int_{\mathbb R_+^2}
		e^{2\gamma\Psi}[\partial_x^j,a+1]\partial_y^2 q\,\partial_y\partial_x^j q\,dxdy
		\\
		&\le
		C\sqrt{\bar\rho}\|a\|_{H^{4}_x}
		\sum_{j=0}^{3}
		\left\|
		\frac1{\sqrt{\rho}}e^{\gamma\Psi}\partial_y\partial_x^j q
		\right\|_{L^2}
		\|e^{\gamma\Psi}\partial_y q\|_{H^{3,0}}
		\\
		&\leq
		\eta\sum_{j=0}^{3}
		\left\|
		\frac1{\sqrt{\rho}}e^{\gamma\Psi}\partial_y\partial_x^j q
		\right\|_{L^2}^2
		+
		C_{\eta,\bar\rho}\|a\|_{H^{4}_x}^2
		\|e^{\gamma\Psi}\partial_y q\|_{H^{3,0}}^2 .
	\end{split}
	\label{yw-comm}
\end{align}

Next, we deal with the three terms $L_1$--$L_3$. In what follows, we shall frequently use some standard lemmas stated in Section \ref{s2}; some routine details will therefore be omitted for brevity.

For the first term, it
can be treated like $K_2$.
\begin{align}
	\begin{split}
		\left|  L_1  \right| &\leq 2\eta \sum_{j=0}^{3}\left\| \frac{1}{\sqrt{\rho}} e^{\gamma\Psi} \partial_y \partial_x^j q\right\|_{L^2 }^2+C_{\eta}\|a\|_{H^{4}_x}^2	\|e^{\gamma\Psi}\partial_y^2 w^s\|_{L^2_y}^2
		+\left\langle t\right\rangle ^{-1}\|e^{{\gamma}\Psi}\partial_y w\|_{H^{4,0} }^2
		\\
		&\quad + C_\eta \left(\langle t \rangle^2\|w + w^s + \varepsilon f(t)\chi'(y)\|_{L^\infty }^2 +\langle t \rangle^2\|e^{{\gamma}\Psi}\partial_x \partial_y w\|_{L^\infty }^2 +\langle t \rangle^2\|e^{{\gamma}\Psi}\partial_x  w\|_{L^\infty }^2\right)\|e^{{\gamma}\Psi}q\|_{H^{3,0} }^2,
	\end{split}
	\label{yw5}
\end{align}
where we have used the following facts
\begin{align*}
	&\left| \left\langle e^{{\gamma}\Psi} (u + u^s + \varepsilon f(t)\chi(y)) \partial_x\partial_y w ,e^{{\gamma}\Psi} q\right\rangle_{H^{3,0} }\right| 
	\\
	&\lesssim \left(\|\partial_x u\|_{L^\infty } + \langle t \rangle\|\partial_x \partial_y w\|_{L^\infty }\right)\|e^{\gamma\Psi}\partial_y w\|_{H^{3,0} }\|e^{\gamma\Psi}q\|_{H^{3,0} }
	\\
	&\leq \eta \sum_{j=0}^{3}\left\| \frac{1}{\sqrt{\rho}} e^{\gamma\Psi} \partial_y^2 \partial_x^j w\right\|_{L^2 }^2+C_{\eta, \gamma, \bar{\rho}}\langle t \rangle^2\|e^{{\gamma}\Psi}\partial_x \partial_y w\|_{L^\infty }^2\|e^{\gamma\Psi}q\|_{H^{3,0} }^2
	\\
	&\leq \eta \sum_{j=0}^{3}\left\| \frac{1}{\sqrt{\rho}} e^{\gamma\Psi} \partial_y \partial_x^j q\right\|_{L^2 }^2+C_{\eta,\bar\rho}\|a\|_{H^{4}_x}^2
	\|e^{\gamma\Psi}\partial_y^2 w^s\|_{L^2_y}^2+C_{\eta, \gamma, \bar{\rho}}\langle t \rangle^2\|e^{{\gamma}\Psi}\partial_x \partial_y w\|_{L^\infty }^2\|e^{\gamma\Psi}q\|_{H^{3,0} }^2
	,
\end{align*}
and
\begin{align*}
	\begin{split}
		&\left|	\left\langle e^{{\gamma}\Psi} (w + w^s + \varepsilon f(t)\chi'(y)) \partial_x  w ,e^{{\gamma}\Psi}q\right\rangle_{H^{3,0} }\right| \\
		&\lesssim \left(\|w + w^s + \varepsilon f(t)\chi'(y)\|_{L^\infty }\|e^{{\gamma}\Psi}\partial_x w\|_{H^{3,0} } + \|\partial_x w\|_{L^\infty }\|e^{{\gamma}\Psi} w\|_{H^{3,0} }\right)\|e^{{\gamma}\Psi}q\|_{H^{3,0} } \\
		& \leq \eta \sum_{j=0}^{3}\left\| \frac{1}{\sqrt{\rho}} e^{\gamma\Psi} \partial_y \partial_x^j q\right\|_{L^2 }^2+\left\langle t\right\rangle ^{-1}\|e^{{\gamma}\Psi}\partial_y w\|_{H^{4,0} }^2
		\\
		&\quad
		+  C_{\eta, \gamma, \bar{\rho}}\left(\langle t \rangle^2\|w + w^s + \varepsilon f(t)\chi'(y)\|_{L^\infty }^2 + \langle t \rangle^2\|e^{{\gamma}\Psi}\partial_x w\|_{L^\infty }^2\right)\|e^{{\gamma}\Psi} q\|_{H^{3,0} }^2.
	\end{split}
\end{align*}

For the second term, it
can be dealt like $K_3$. Then, by Young's inequality, we have
\begin{align}
	\begin{split}
		\left| L_2\right| 
		&\leq 
		\eta \sum_{j=0}^{3}\left\| \frac{1}{\sqrt{\rho}} e^{\gamma\Psi} \partial_y \partial_x^j q\right\|_{L^2 }^2
		\\
		&
		\quad+
		C_\eta \left\langle t\right\rangle^2 \left( 1+\left\| \frac{\partial_x a}{a+1}\right\|_{H^{3}_x}^2\right) 
		\|\partial_y(w + w^s + \varepsilon f(t)\chi'(y))\|_{L^\infty }^2\|e^{{\gamma}\Psi} w\|_{H^{4,0} }^2
		\\
		&\quad
		+
		C_\eta \left\langle t\right\rangle^{\frac32}\left( 1+\left\| \frac{\partial_x a}{a+1}\right\|_{H^{3}_x}^2\right) 
		\|e^{\gamma\Psi}w \|_{H^{4,0} }^2\|e^{{\gamma}\Psi} \partial_y w\|_{H^{4,0} }^2 
		,
	\end{split}
	\label{yw6}
\end{align}
where we have used the  definition of $v$  
\begin{align*}
	v \stackrel{\text{def}}{=} -\partial_x \varphi - \frac{\partial_x a}{a+1} \int_y^\infty \left[ u + u^s - \varepsilon f(t)(1 - \chi(y')) \right] dy',
\end{align*}
and the following facts:
\begin{align*}
	&\left| \left\langle e^{{\gamma}\Psi} \partial_x\varphi \partial_y \left( w + w^s + \varepsilon f(t)\chi'(y) \right)  ,e^{{\gamma}\Psi} \partial_y q\right\rangle_{H^{3,0} }\right| 
	\\
	&\lesssim
	\|e^{{\gamma}\Psi}\partial_x\varphi \partial_y \left( w + w^s + \varepsilon f(t)\chi'(y) \right)\|_{H^{3,0} } \|e^{{\gamma}\Psi}\partial_y q\|_{H^{3,0} } 
	\\
	&
	\lesssim
	\left(\left\langle t\right\rangle \|\partial_y(w + w^s + \varepsilon f(t)\chi'(y))\|_{L^\infty }\|e^{{\gamma}\Psi}\partial_x w\|_{H^{3,0} } + \|\partial_x \varphi\|_{L^\infty }\|e^{{\gamma}\Psi} \partial_y w\|_{H^{3,0} }\right)
	\|e^{{\gamma}\Psi}\partial_y q\|_{H^{3,0} } ,
\end{align*}
and
\begin{align*}
	&\left| \left\langle e^{{\gamma}\Psi} \left(  \frac{\partial_x a}{a+1} \int_y^\infty \left[ u  + u^s - \varepsilon f(t)(1 - \chi(y')) \right] dy'\right)  \partial_y \left( w + w^s + \varepsilon f(t)\chi'(y) \right)  ,e^{{\gamma}\Psi} \partial_y q\right\rangle_{H^{3,0} }\right|
	\\
	&\lesssim
	\left\|e^{{\gamma}\Psi}\left(  \frac{\partial_x a}{a+1} \int_y^\infty \left[ u  + u^s - \varepsilon f(t)(1 - \chi(y')) \right] dy'\right)  \partial_y \left( w + w^s + \varepsilon f(t)\chi'(y) \right)\right\|_{H^{3,0} } \|e^{{\gamma}\Psi}\partial_y q\|_{H^{3,0} } 
	\\
	&
	\lesssim
	\left\langle t\right\rangle\left\| \frac{\partial_x a}{a+1}\right\|_{H^{3}_x}
	\|e^{\gamma_0\Psi}\partial_y(w + w^s + \varepsilon f(t)\chi'(y))\|_{L^\infty }\|e^{{\gamma}\Psi}  w\|_{H^{3,0} }\|e^{{\gamma}\Psi}\partial_y q\|_{H^{3,0} } 
	\\
	&\quad+
\left\| \frac{\partial_x a}{a+1}\right\|_{H^{3}_x}
	\|e^{\gamma_0\Psi}\varphi\|_{L^\infty }\|e^{{\gamma}\Psi} \partial_y w\|_{H^{3,0} }\|e^{{\gamma}\Psi}\partial_y q\|_{H^{3,0} }  
	.
\end{align*}

For last term, along the same line as $K_4$, we can obtain
\begin{align}
	\begin{split}
		\left| L_3\right| &\leq \eta \sum_{j=0}^{3}\left\| \frac{1}{\sqrt{\rho}} e^{\gamma\Psi} \partial_y \partial_x^j q\right\|_{L^2 }^2
		+
		C_\eta \varepsilon^2
 \left\|
		a
		\right\|_{H^{4}_x}^2
		|f(t)|^2
		\\
		&
		\quad+
		C_\eta \left\langle t\right\rangle^2 \left\| \frac{\partial_x a}{a+1}\right\|_{H^{3,0}_x }^2
		\left\|e^{\gamma_0 \Psi} \partial_y\left( w + w^s +\varepsilon f(t) \chi(y) \right)\right\|_{L^\infty }^2 \left\| e^{\gamma \Psi}  w \right\|_{H^{4,0} }^2
		\\
		&
		\quad
		+C_\eta \varepsilon^2
		\left( 1+\left\|
		a
		\right\|_{H^{4}_x}^2\right)^2 \left\|
		a
		\right\|_{H^{4}_x}^2
		\left(|f(t)|^2+\left\|e^{ \gamma \Psi} 
		\partial_y^2 w^s
		\right\|_{L^2_y}^2\right) \left( 1+\left\|e^{ \gamma \Psi} 
		w^s
		\right\|_{L^2_y}^2+|f(t)|^2+|f'(t)|^2\right) ,
	\end{split}
	\label{yw7}	
\end{align}
where we have used the following facts
\begin{align*}
\left\|e^{ \gamma \Psi} \partial_t
\left(
\frac{a}{a+1}w^s
\right)\right\|_{H^{3,0}} 
&\lesssim \left\|e^{ \gamma \Psi} \partial_t \rho
~w^s
\right\|_{H^{3}_x}
+\left\|e^{ \gamma \Psi} 
\frac{a}{a+1} \partial_t w^s
\right\|_{H^{3,0}}
	\\
& \leq C \varepsilon \|\partial_x\rho \|_{H^3_x} \left\|e^{ \gamma \Psi} 
 w^s
\right\|_{L^2_y}
|f(t)|
+
C_{\chi}\varepsilon\left\|
\frac{a}{a+1} 
\right\|_{H^{3}_x}
\left\|e^{ \gamma \Psi} 
 \partial_y^2 w^s
\right\|_{L^2_y}(|f(t)|+|f'(t)|)
\\
& \lesssim \varepsilon
\left( 1+\left\|
a
\right\|_{H^{4}_x}^2\right) \left\|
a
\right\|_{H^{4}_x}
 \left\|e^{ \gamma \Psi} 
 w^s
\right\|_{L^2_y}|f(t)|
+
 \left\|
a
\right\|_{H^{4}_x}
\left\|e^{ \gamma \Psi} 
\partial_y^2 w^s
\right\|_{L^2_y} (|f(t)|+|f'(t)|)
,
\end{align*}
and
\begin{align*}
	\left\| e^{ \gamma \Psi}  F_3 \right\|_{H^{3,0} }	
	&=\left\| e^{ \gamma \Psi} \frac{\partial_x a}{a+1}  \left( u + u^s - \varepsilon f(t)(1 - \chi(y)) \right) (w + w^s + \varepsilon f(t)\chi'(y)) \right\|_{H^{3,0} }
	\\
	&\leq
	C\left\langle t\right\rangle \left\| \frac{\partial_x a}{a+1}\right\|_{H^{3,0}_x }
	\left\|e^{ \gamma_0 \Psi} \partial_y\left( w + w^s +\varepsilon f(t) \chi(y) \right)\right\|_{L^\infty } \left\| e^{\gamma \Psi}  w \right\|_{H^{4,0} }
	.
\end{align*}

Substituting the above estimates on \eqref{yw-comm}--\eqref{yw7} into \eqref{yw1-new}, and multiplying the resulting inequality by $\langle t\rangle^{\frac92-2\kappa}$, we obtain
\begin{align*}
	\begin{split}
		&\frac{d}{dt}\left\| \langle t\rangle^{\frac94-\kappa}e^{\gamma\Psi} q\right\|_{H^{3,0} }^2
		+\left( 2-\gamma-10\eta \right)\langle t\rangle^{\frac92-2\kappa}
		\sum_{j=0}^{3}\left\| \frac{1}{\sqrt{\rho}} e^{\gamma\Psi} \partial_y \partial_x^j q\right\|_{L^2 }^2
		\\
		& \leq C \theta_4 \left\|\langle t\rangle^{\frac94-\kappa}e^{{\gamma}\Psi}q\right\|_{H^{3,0} }^2
		+ C_{\eta}\|a\|_{H^{4}_x}^2
		\|e^{\gamma\Psi}\partial_y q\|_{H^{3,0}}^2
		+C_{\eta}\|a\|_{H^{4}_x}^2	\|e^{\gamma\Psi}\partial_y^2 w^s\|_{L^2_y}^2
			\\
		&
		\quad+
		C_\eta \left\langle t\right\rangle^3 \left( 1+\left\| \frac{\partial_x a}{a+1}\right\|_{H^{3}_x}^2\right) 
		\|e^{\gamma_0 \Psi}\partial_y(w + w^s + \varepsilon f(t)\chi'(y))\|_{L^\infty }^2\|\left\langle t\right\rangle^{\frac{7}{4}-\kappa}e^{{\gamma}\Psi} w\|_{H^{4,0} }^2
		\\
		&\quad
		+
		C_\eta \left\langle t\right\rangle^{-1+2\kappa}\left( 1+\left\| \frac{\partial_x a}{a+1}\right\|_{H^{3}_x}^2\right) 
		\left( 1+\| \left\langle t\right\rangle^{\frac{7}{4}-\kappa}e^{\gamma\Psi}w \|_{H^{4,0} }^2\right) \left\| \left\langle t\right\rangle^{\frac{7}{4}-\kappa}e^{{\gamma}\Psi} \partial_y w\right\| _{H^{4,0} }^2
		\\
		&
		\quad
		+C_\eta \varepsilon^2 \left\langle t\right\rangle^{\frac{9}{2}-2\kappa}
	\left( 1+\left\|
	a
	\right\|_{H^{4}_x}^2\right)^2 \left\|
	a
	\right\|_{H^{4}_x}^2
	\left(|f(t)|^2+\left\|e^{ \gamma \Psi} 
	\partial_y^2 w^s
	\right\|_{L^2_y}^2\right) \left(1+ \left\|e^{ \gamma \Psi} 
	w^s
	\right\|_{L^2_y}^2+|f(t)|^2+|f'(t)|^2\right)
		\\
		&\quad
		+\left( \frac92-2\kappa\right) \left\|\langle t\rangle^{\frac74-\kappa}e^{{\gamma}\Psi} \partial_y w \right\|_{H^{3,0} }^2
,
	\end{split}
\end{align*}
where  $\theta_4$ is defined in \eqref{theta4}. Taking
$\eta=\frac{1-\gamma}{10}>0$, we have
$2-\gamma-10\eta=1$. Hence the above inequality yields the asserted estimate, and the proof of the proposition is complete.

\end{proof}

\section{Proof of Theorem 1.1} \label{s7}
\paragraph{Proof of Theorem 1.1.}
With the a priori estimates established in Sections \ref{s3}--\ref{s6}, we now complete the proof of Theorem \ref{t1} by a standard continuation argument. As in Section \ref{s3}, let $\mathcal T^*$ be the lifespan of a sufficiently smooth solution $(\rho,u,v)$ to system \eqref{uprandtl}. We introduce the bootstrap time
\begin{equation}
	\begin{aligned}
		T^\star \overset{\mathrm{def}}{=}
		\sup \Bigl\{\, t\in [0,\mathcal T^*) :\;
		\|\langle \tau\rangle^{\frac94-\kappa}e^{\gamma\Psi} \partial_y w
		\|_{L_\tau^\infty(0,t;H^{3,0})}^2
+
		\int_0^t
		\langle \tau\rangle^{\frac92-2\kappa}
		\left\|
		e^{\gamma\Psi}\partial_y^2
	 w(\tau)
		\right\|_{H^{3,0}}^2\,d\tau
		<2c_1
		\Bigr\},
	\end{aligned}
	\label{8.1}
\end{equation}
where $c_1>0$ is a small constant to be chosen later.

We first show that $T^\star\le T^{\#}$, where $T^{\#}$ is defined in \eqref{T}. Indeed, for any $t\le T^\star$, Sobolev embedding together with \eqref{8.1} yields
\begin{align*}
\|e^{\gamma_0\Psi}w(t)\|_{L^\infty}
&\le C\langle t\rangle^{\frac14}\|e^{\gamma\Psi}\partial_y w(t)\|_{H^{1,0}}
\le C c_1^{\frac12}\langle t\rangle^{-2+\kappa}.
\end{align*}
Choosing $c_1$ sufficiently small so that $Cc_1^{1/2}\le c_0/2$, and using \eqref{CM}, we obtain
\begin{align*}
&w(t,x,y)+w^s(t,y)+\varepsilon f(t)\chi'(y)
\\
&\qquad\ge \bigl(w^s(t,y)+\varepsilon f(t)\chi'(y)\bigr)-|w(t,x,y)|
\ge \frac{c_0}{2}\langle t\rangle^{-2+\kappa}e^{-\gamma_0\Psi(t,y)},
\end{align*}
for all $t\le T^\star$. Hence $T^\star\le T^{\#}$.

We next estimate $\int_0^T\theta_i(t)\,dt$ for $i=1,2,3,4$ and any $T\le T^\star$. By the support properties of $\chi$, namely $\mathrm{Supp}\,\chi'\subset [1,2]$ and $\mathrm{Supp}\,(1-\chi)\subset [0,2]$, we have for all $0\le t\le T$,
\begin{align}
\|\chi\|_{L_y^\infty}+\|e^{\Psi(t)}\chi'\|_{H^{0,2}}+\|e^{\Psi(t)}(1-\chi)\|_{L_y^2}+\|e^{\Psi(t)}(1-\chi)\|_{L_y^\infty}\le C.
\label{8.2}
\end{align}
Moreover, by \eqref{Ca},
\begin{align}
 \sup_{t \in [0, +\infty)} \left\|  \left\langle t\right\rangle^{\frac12+\kappa_0}\frac{\partial_x a}{a+1}(t) \right\|_{H^m} \lesssim C_a.
 \label{8.3}
\end{align}
Using \eqref{CM} together with \eqref{u0}, we further get
\begin{equation}\label{initial}
\sum_{j=0}^m \left\| e^{\frac{y^2}{8}} g_j(0)\right\|_{L^2}
+\left\| e^{\frac{y^2}{8}} G_0\right\|_{H^{5,0}}
+\left\| e^{\frac{\gamma y^2}{8}} w_0\right\|_{H^{4,0}}
+\left\| e^{\frac{\gamma y^2}{8}} \partial_y w_0\right\|_{H^{3,0}}
\le C_\gamma \varepsilon_0.
\end{equation}

Recalling the definition of $\theta_1$ in \eqref{theta1}, and using the conditions \eqref{Cf}, \eqref{CM1}, and \eqref{8.1}-\eqref{8.3}, we obtain for any $T\le T^\star$,
\begin{align}
\begin{split}
\int_0^T \theta_1(t)\,dt
&\lesssim
(1+C_a^2)(c_1^2+C_M^2+C_f^2)+c_1^2+C_a^2
\\
&\quad +(1+C_a^2)\int_0^T
\Bigl( \langle t\rangle^{\frac34}\| e^{\gamma\Psi} \partial_y w \|_{H^{2,0}}
+\langle t\rangle^{\frac52}\| e^{\gamma\Psi} \partial_y w \|_{H^{2,0}}^2
+\left\langle t\right\rangle \left\| e^{\gamma_0\Psi} \partial_y w^s \right\|_{L^\infty_y }
\\
&\qquad\qquad\qquad
 +\left\langle t\right\rangle^{3}\left\| e^{\gamma_0\Psi}\partial_y w^s \right\|_{L^\infty_y }^2
+\langle t\rangle^{-4+2\kappa}
+\langle t\rangle^{-\frac32+\kappa}\Bigr)\,dt
\\
&\quad +C_a\int_0^T \langle t\rangle^{-\frac12}
\Bigl(\varepsilon \langle t\rangle^{-\frac12}|f(t)|\Bigr)\,dt
+C_a^2\int_0^T \langle t\rangle^{-1}
\Bigl(\varepsilon^2 |f(t)|^2\Bigr)\,dt
\\
&\quad +(1+C_a^2)(c_1^2+C_M^2+C_f^2)
\sup_{0\le t\le T}\Bigl\{\langle t\rangle
\bigl(\langle t\rangle^{\frac34}\| e^{\gamma\Psi} \partial_y w \|_{H^{2,0}}
+\langle t\rangle^{\frac52}\| e^{\gamma\Psi} \partial_y w \|_{H^{2,0}}^2
\\
&\qquad\qquad\qquad
\left\langle t\right\rangle \left\| e^{\gamma_0\Psi} \partial_y w^s \right\|_{L^\infty_y }
+\left\langle t\right\rangle^{3}\left\| e^{\gamma_0\Psi}\partial_y w^s \right\|_{L^\infty_y }^2
+\langle t\rangle^{-4+2\kappa}
+\langle t\rangle^{-\frac32+\kappa}\bigr)\Bigr\}
\\
&\quad +(c_1^2+C_M^2+C_f^2)
\Bigl(C_a\sup_{0\le t\le T}\bigl(\langle t\rangle^{-\frac12-\kappa_0}\langle t\rangle\varepsilon \langle t\rangle^{-\frac12}|f(t)|\bigr)
+C_a^2\sup_{0\le t\le T}\bigl(\langle t\rangle^{-1-2\kappa_0}\langle t\rangle \varepsilon^2 |f(t)|^2\bigr)\Bigr)
\\
&\lesssim (1+C_a^2)(c_1^2+C_M^2+C_f^2)(1+c_1^2+C_M^2+C_f^2)\lesssim (1+C_a^2)(1+c_1^2+C_M^2+C_f^2)^2,
\end{split}
\label{8.4}
\end{align}
where we have used the facts
\begin{align*}
 &\int_{0}^{T}\langle t \rangle^{-1}\langle t \rangle^{{5}-2\kappa} \|e^{\gamma_0 \Psi} \partial_y(w + w^s + \varepsilon f(t)\chi'(y))\|_{L^\infty }^2 dt \\
&\lesssim \int_{0}^{T}\langle t \rangle^{\frac{9}{2}-2\kappa}\Bigl( \|e^{\gamma_0 \Psi} \partial_y^2w \|_{H^{3,0}}^2
+\|e^{\gamma_0 \Psi} \partial_y^2 w^s \|_{L^2_y }^2
+\langle t\rangle^{-\frac12} \varepsilon^2 |f(t)|^2\Bigr) \, dt
\\
& \lesssim \left( c_1^2+C_M^2+C_f^2\right) ,
\end{align*}
and
\begin{align*}
	\sup_{0\le t\le T}
	\Bigl(
	\left\langle t\right\rangle^{4}
	\left\| e^{\gamma_0\Psi}\partial_y w^s \right\|_{L^\infty_y }^2
	\Bigr)
	\leq 
C	\sup_{0\le t\le T}
	\Bigl(
	\left\langle t\right\rangle^{-2\kappa}
	\left\| e^{(\gamma_0+\gamma-2)\Psi} \right\|_{L^\infty_y }^2
	\Bigr) \leq C.
\end{align*}
In addition, we have also
\begin{align}
\begin{split}
&\int_0^T \Bigl( \bigl( 1+\langle t\rangle^{\frac32-\kappa} \bigr) \| a\|_{H^{m+1}_x }^2\| e^{\Psi}\partial_y w^s\|_{L^2 }^2
+\varepsilon^2 \left\langle t\right\rangle|f(t)|^2 \Bigr)\,dt
\\
&\lesssim C_a^2 	C_M^2 \int_0^T  \bigl( \left\langle t\right\rangle^{-1-2\kappa_0} +\langle t\rangle^{\frac12-\kappa-2\kappa_0} \bigr) \left\langle t\right\rangle^{-4+2\kappa}  \| e^{(-1+\gamma)\Psi}\|_{L^2 }^2
 \,dt
+\varepsilon^2   C_f^2
\\
&\lesssim \varepsilon_0^2(C_M^2+C_f^2).
\end{split}
\label{8.5}
\end{align}
Substituting \eqref{8.4} and \eqref{8.5} into \eqref{g0}, and then applying Gronwall's inequality and \eqref{initial}, we obtain
\begin{align}
\begin{split}
&\sum_{j=0}^m\| e^{\Psi} g_j\|_{L_T^\infty (L^2)}^2
+\sum_{j=0}^m\int_{0}^{T}
\| e^{\Psi} \partial_y g_j\|_{L^2}^2\,dt
\\
&
\leq
e^{C \int_0^T \theta_1 \,dt} \left\lbrace  \sum_{j=0}^m \left\| e^{\frac{y^2}{8}} g_j(0)\right\|_{L^2}^2
+ \int_0^T\left( \left( 1+\left\langle t\right\rangle^{\frac32-\kappa} \right) \left\| a\right\|_{H^{m+1}_x }^2\left\| e^{\Psi}\partial_y w^s\right\|_{L^2_y }^2 \right\rbrace 
+
\varepsilon^2 \left\langle t\right\rangle |f(t)|^2 \,dt \right)  
\\
&\leq Ce^{C(1+C_a^2)(1+c_1^2+C_M^2+C_f^2)^2} \left( \varepsilon_0^2+\varepsilon_0^2(C_M^2+C_f^2)\right)
\leq Ce^{C(1+C_a^2)(1+c_1^2+C_M^2+C_f^2)^2}  \varepsilon_0^2 .
\end{split}
\label{8.6}
\end{align}

For $\int_0^T \theta_2 \, dt$, recalling the definition of $X(t)$ in \eqref{XT}, and using the conditions \eqref{Cf}, \eqref{CM1}, and \eqref{8.1}-\eqref{8.3} again, we infer that
\begin{align}
\begin{split}
\int_0^T \theta_2(t)\,dt
&\lesssim 1+c_1^2
+\int_0^T
\langle t\rangle^{-\frac{4m(1-\kappa)-31+24\kappa}{2(2m-11)}}
\Bigl(\langle t\rangle^{-\frac12}\Bigr)^{\frac{2(m-4)}{2m-11}}\,dt
\\
&\quad +\int_0^T
\langle t\rangle^{-\frac{4m(1-\kappa)-31+24\kappa}{2(2m-11)}}
\Bigl(
\langle t\rangle^{\frac94-\kappa}
\bigl(
\|e^{\gamma\Psi}\partial_y^2 w\|_{L^2}
+\|e^{\gamma_0\Psi}\partial_y^2 w^s\|_{L_y^2}
\bigr)
+\varepsilon \langle t\rangle^{2-\kappa}|f(t)|
\Bigr)^{\frac{2(m-4)}{2m-11}}\,dt
\\
&\quad +(1+C_a^2)(1+c_1^2+C_M^2+C_f^2)
\\
&\lesssim (1+C_a^2)(1+c_1^2+C_M^2+C_f^2),
\end{split}
\label{8.7}
\end{align}
where we have used
\begin{align*}
\int_0^T
\langle t\rangle^{-\frac{4m(1-\kappa)-31+24\kappa}{2(2m-11)}}
\Bigl(\langle t\rangle^{-\frac12}\Bigr)^{\frac{2(m-4)}{2m-11}}\,dt
&=
\int_0^T
\langle t\rangle^{-\frac{6m-4m\kappa-39+24\kappa}{2(2m-11)}}\,dt
\lesssim 1,
\end{align*}
and
\begin{align*}
\int_0^T
\langle t\rangle^{-\frac{4m(1-\kappa)-31+24\kappa}{2(2m-11)}}
X(t)^{\frac{2(m-4)}{2m-11}}\,dt
&\lesssim
\left(
\int_0^T
\langle t\rangle^{-\frac{4m(1-\kappa)-31+24\kappa}{2(m-7)}}\,dt
\right)^{\frac{m-7}{2m-11}}
\\
&\quad\times
\left(
\int_0^T X(t)^2\,dt
\right)^{\frac{m-4}{2m-11}}
\lesssim (c_1^2+C_M^2+C_f^2).
\end{align*}
These estimates are valid provided that
\begin{align*}
m\ge m_0>\frac{17-24\kappa}{2(1-2\kappa)}
=5+\frac{7-4\kappa}{2(1-2\kappa)}.
\end{align*}
Moreover, note that condition \eqref{CM1}, together with the elementary Gaussian tail estimate, implies
\begin{align*}
	\|e^\Psi w^s(t)\|_{L_y^2}^2
	&=
	\left\|
	e^\Psi
	\int_y^\infty \partial_z^2u^s(t,z)\,dz
	\right\|_{L_y^2}^2
	\\
	&\lesssim
	C_M^2\langle t\rangle^{-4+2\kappa}
	\left\|
	e^\Psi
	\int_y^\infty e^{-(2-\gamma)\Psi(t,z)}\,dz
	\right\|_{L_y^2}^2
	\\
	&\le
	C_\gamma C_M^2\langle t\rangle^{-4+2\kappa}
	\left\|
	\langle t\rangle^{\frac12}
	e^{-(1-\gamma)\Psi}
	\right\|_{L_y^2}^2
	\\
	&\le
	C_\gamma C_M^2
	\langle t\rangle^{-4+2\kappa}
	\langle t\rangle
	\langle t\rangle^{\frac12}
	\\
	&=
	C_\gamma C_M^2
	\langle t\rangle^{-\frac52+2\kappa},
\end{align*}
and hence
\begin{align}
	\begin{split}
		&\int_0^T \langle t\rangle^{\frac52-2\kappa}
		\Bigl( \| a\|_{H^{6}_{x} }^2
		+\left\| \frac{\partial_x a}{a+1}\right\|_{H^{6}_{x} }^2\Bigr)
		\Bigl( \| e^{\Psi}  w^s\|_{L^{2}_y }^2
		+\| e^{\gamma_0\Psi}  w^s\|_{L^{2}_y }^2
		+\varepsilon^2|f(t)|^2\Bigr)\,dt
		\\
		&\lesssim C_a^2 \int_0^T \langle t\rangle^{-1-\kappa_0}  \,dt+C_a^2 ( C_M^2+C_f^2)
		\\
		&\lesssim C_a^2 ( C_M^2+C_f^2).
	\end{split}
	\label{8.8}
\end{align}

Collecting \eqref{8.6}--\eqref{8.8}, using \eqref{initial}, and applying   Gronwall's inequality to \eqref{G0}, we obtain
\begin{align}
\begin{split}
\sup_{0\le t\le T}\| \langle t\rangle^{\frac54-\kappa} e^{\Psi}  G(t)\|_{H^{5,0}}^2
+\int_{0}^{T}
\| \langle t\rangle^{\frac54-\kappa}e^{\Psi} \partial_y  G\|_{H^{5,0} }^2\,dt
\le Ce^{(1+C_a^2)( 1+c_1^2+C_M^2+C_f^2)^2} \varepsilon_0^2.
\end{split}
\label{8.9}
\end{align}

By a similar process of argument, from \eqref{w0} and \eqref{yw0}, we deduce
\begin{align}
\begin{split}
\sup_{0\le t\le T}\| \langle t\rangle^{\frac74-\kappa} e^{\gamma\Psi}  w(t)\|_{H^{4,0}}^2
+\int_{0}^{T}
\| \langle t\rangle^{\frac74-\kappa}e^{\gamma\Psi} \partial_y^2 w\|_{H^{4,0} }^2\,dt
\le Ce^{(1+C_a^2)( 1+c_1^2+C_M^2+C_f^2)^2} \varepsilon_0^2,
\end{split}
\label{8.10}
\end{align}
and
\begin{align}
\begin{split}
&\sup_{0\le t\le T} \left\| \left\langle t\right\rangle^{\frac{9}{4}-\kappa}e^{\gamma\Psi}
\left( \partial_y w+\frac{a}{a+1}\partial_y w^s\right) \right\|_{H^{3,0} }^2
+\int_{0}^{T}
\langle t\rangle^{\frac92-2\kappa}
\left\|  e^{\gamma\Psi}
\partial_y 
\left( \partial_y w+\frac{a}{a+1}\partial_y w^s\right)\right\|_{H^{3,0} }^2
\,dt
\\
&
\le Ce^{(1+C_a^2)( 1+c_1^2+C_M^2+C_f^2)^2} \left( \varepsilon_0^2 +
\left\|  \frac{a}{a+1}(0) \right\|_{H^{3}_x }^2  \left\| e^{\gamma} \partial_y w^s(0) \right\|_{L^2_y }^2\right)
\\
&
\le Ce^{(1+C_a^2)( 1+c_1^2+C_M^2+C_f^2)^2} \left( \varepsilon_0^2 +
C_{\bar\rho}  C_a^2 M^2\right) 
\\
&
\le Ce^{(1+C_a^2)( 1+c_1^2+C_M^2+C_f^2)^2}  \varepsilon_0^2 .
\end{split}
\label{8.11}
\end{align}

Combining \eqref{8.10} with Lemma \ref{po}, we further obtain the following large-time decay estimate for $u$:
\begin{align}
	\begin{split}
		&\sup_{0\le t\le T}
		\left\| \langle t\rangle^{\frac54-\kappa} e^{\gamma\Psi}  u(t)\right\|_{H^{5,0}}^2
		+\int_{0}^{T}
		\left\| \langle t\rangle^{\frac74-\kappa}e^{\gamma\Psi} \partial_y u\right\|_{H^{5,0} }^2\,dt
		\\
		&\qquad\le
		Ce^{(1+C_a^2)( 1+c_1^2+C_M^2+C_f^2)^2} \varepsilon_0^2.
	\end{split}
	\label{8.12}
\end{align}

In addition, by \eqref{8.11} and \eqref{CM1}, we also have 
\begin{align*}
	&\sup_{0\le t\le T}
	\left\|
	\langle t\rangle^{\frac{9}{4}-\kappa}e^{\gamma\Psi}
	\partial_y w
	\right\|_{H^{3,0}}^2
	\\
	&\leq
	C e^{(1+C_a^2)(1+c_1^2+C_M^2+C_f^2)^2}
	\varepsilon_0^2
	+
	\sup_{0\le t\le T}
	\left\|
	\langle t\rangle^{\frac{9}{4}-\kappa}e^{\gamma\Psi}
	\frac{a}{a+1}\partial_y w^s
	\right\|_{H^{3,0}}^2
	\\
	&\leq
	C e^{(1+C_a^2)(1+c_1^2+C_M^2+C_f^2)^2}
	\varepsilon_0^2
	+
	\sup_{0\le t\le T}
	\Biggl(
	\left\|
	\langle t\rangle^{\frac{1}{2}+\kappa_0}
	\frac{a}{a+1}
	\right\|_{H_x^3}^2
	\left\|
	\langle t\rangle^{\frac{7}{4}-\kappa-\kappa_0}
	e^{\gamma\Psi}
	\partial_y w^s
	\right\|_{L_y^2}^2
	\Biggr)
		\\
	&\leq
	C e^{(1+C_a^2)(1+c_1^2+C_M^2+C_f^2)^2}
	\varepsilon_0^2
	+
	C_{\bar\rho}C_a^2 C_M^2 \sup_{0\le t\le T}
	\Biggl(
	\left\|
	\langle t\rangle^{-\frac{1}{4}-\kappa_0}
	e^{(-2+2\gamma)\Psi}
	\right\|_{L_y^2}^2
	\Biggr)
	\\
	&\leq
	C e^{(1+C_a^2)(1+c_1^2+C_M^2+C_f^2)^2}
	\varepsilon_0^2
	+
	C_{\bar\rho}\varepsilon_0^2 C_M^2\sup_{0\le t\le T}
	\Biggl(
	\langle t\rangle^{-2\kappa_0}
	\Biggr)\leq
	C e^{CC^*}
	\varepsilon_0^2,
\end{align*}
and
\begin{align*}
\int_{0}^{T}
\langle t\rangle^{\frac92-2\kappa}
\left\|  e^{\gamma\Psi}
\partial_y^2 w\right\|_{H^{3,0} }^2
\,dt  
&
\leq
C e^{(1+C_a^2)(1+c_1^2+C_M^2+C_f^2)^2}
\varepsilon_0^2
\\
&\quad
+
\int_{0}^{T}
\left\|
\langle t\rangle^{\frac{9}{4}-\kappa}e^{\gamma\Psi}
\frac{a}{a+1}\partial_y^2 w^s
\right\|_{H^{3,0}}^2 \,dt  
\leq 	C e^{CC^*}
\varepsilon_0^2,
\end{align*}
where $C^*>0$ is chosen so that
$$
(1+C_a^2)( 1+c_1^2+C_M^2+C_f^2)^2 \le C^*.
$$
In particular, if $\varepsilon_0$ is sufficiently small and satisfies
$$
Ce^{C C^*}\varepsilon_0^2 \le c_1,
$$
namely,
$$
\varepsilon_0 \le \sqrt{\frac{c_1}{C}}\,e^{-\frac12 C C^*},
$$
then, for every $T\le T^\star$,
\begin{align}
\|\langle t\rangle^{\frac94-\kappa}e^{\gamma\Psi} \partial_y w\|_{L_T^\infty(H^{3,0})}^2
+\int_0^T \langle t\rangle^{\frac92-2\kappa}\left\|  e^{\gamma\Psi}
\partial_y^2  w\right\|_{H^{3,0} }^2\,dt
\le c_1.
\label{8.13}
\end{align}

Hence the bootstrap bound in \eqref{8.1} is strictly improved. If $T^\star<\mathcal T^*$, the standard continuity argument allows us to extend the bootstrap inequality beyond $T^\star$, which contradicts the maximality of $T^\star$. Therefore $T^\star=+\infty$. Since $T^\star\le \mathcal T^*$ by definition, we also obtain $\mathcal T^*=+\infty$. Finally, combining \eqref{8.6}, \eqref{8.9}, \eqref{8.10}, \eqref{8.12}  and  \eqref{8.13}, we obtain the desired global weighted estimates stated in Theorem \ref{t1}. This completes the proof of Theorem \ref{t1}.

\begin{appendices}

\section{Equation for $g_m$} \label{ba}

In this appendix, we derive the evolution equation satisfied by $g_m$.
Recall that $g_m$ is defined by
\begin{align*}
	g_m
	=
	\partial_x^m w
	-
	\frac{\partial_y(w + w^s + \varepsilon f(t)\chi'(y))}{w + w^s + \varepsilon f(t)\chi'(y)}
	\partial_x^m u
	=
	\left(w + w^s + \varepsilon f(t)\chi'(y)\right)
	\partial_y\left(
	\frac{\partial_x^m u}{w + w^s + \varepsilon f(t)\chi'(y)}
	\right).
\end{align*}

Applying $\partial_x^m$ to both $\eqref{uprandtl}_1$ and $\eqref{vprandtl}_1$, we obtain
\begin{equation}
	\begin{cases}
		\partial_t \partial_x^m u
		+\left(u + u^s + \varepsilon f(t)\chi(y)\right)\partial_x^{m+1} u
		+\left( w + w^s + \varepsilon f(t)\chi'(y) \right)\partial_x^m v
		-(a+1)\partial_y^2 \partial_x^m u
		\\
		\qquad=
		\left[ u, \partial_x^m\right] \partial_x u
		+\left[ w, \partial_x^m\right] v
		+[a, \partial_x^m]\partial_y^2 (u + u^s)
		+\partial_x^m a\,\varepsilon f(t)\chi''(y),
		\\[1ex]
		\partial_t \partial_x^m w
		+\left(u + u^s + \varepsilon f(t)\chi(y)\right)\partial_x^{m+1}w
		+\partial_y \left( w + w^s + \varepsilon f(t)\chi'(y) \right) \partial_x^m v
		-(a+1)\partial_y^2 \partial_x^m w
		\\
		\qquad=
		\left[ u, \partial_x^m\right] \partial_x w
		+\left[ \partial_yw, \partial_x^m\right] v
		+[a, \partial_x^m]\partial_y^2 (w + w^s)
		+\partial_x^m a\,\varepsilon f(t)\chi'''(y)
		+\partial_x^m F_3.
	\end{cases}
	\label{A1}
\end{equation}

To eliminate the term $\partial_x^m v$, we subtract
$$
\frac{\partial_y(w + w^s + \varepsilon f(t)\chi'(y))}{w + w^s + \varepsilon f(t)\chi'(y)} \times
\eqref{A1}_1
$$
from $\eqref{A1}_2$.
Then we get
\begin{align}
	\begin{split}
		&\left(
		\partial_t \partial_x^m w
		-
		\frac{\partial_y(w + w^s + \varepsilon f(t)\chi'(y))}{w + w^s + \varepsilon f(t)\chi'(y)}
		\partial_t \partial_x^m u
		\right)
		\\
		&\quad
		+\left(u + u^s + \varepsilon f(t)\chi(y)\right)
		\left(
		\partial_x^{m+1}w
		-
		\frac{\partial_y(w + w^s + \varepsilon f(t)\chi'(y))}{w + w^s + \varepsilon f(t)\chi'(y)}
		\partial_x^{m+1} u
		\right)
		\\
		&\quad
		-(a+1)
		\left(
		\partial_y^2 \partial_x^m w
		-
		\frac{\partial_y(w + w^s + \varepsilon f(t)\chi'(y))}{w + w^s + \varepsilon f(t)\chi'(y)}
		\partial_y^2\partial_x^m u
		\right)
		\\
		&=
		\left[ u, \partial_x^m\right] \partial_x w
		+\left[ \partial_yw, \partial_x^m\right] v
		+[a, \partial_x^m]\partial_y^2 (w + w^s)
		\\
		&\quad
		-
		\frac{\partial_y(w + w^s + \varepsilon f(t)\chi'(y))}{w + w^s + \varepsilon f(t)\chi'(y)}
		\left(
		\left[ u, \partial_x^m\right] \partial_x u
		+\left[ w, \partial_x^m\right] v
		+[a, \partial_x^m]\partial_y^2 (u + u^s)
		\right)
		\\
		&\quad
		+\partial_x^m a\,\varepsilon f(t)
		\left(
		\chi'''(y)
		-
		\frac{\partial_y(w + w^s + \varepsilon f(t)\chi'(y))}{w + w^s + \varepsilon f(t)\chi'(y)}
		\chi''(y)
		\right)
		+\partial_x^m F_3 .
	\end{split}
	\label{A2}
\end{align}

Observing that
\begin{align*}
	&\left[ u, \partial_x^m\right] \partial_x w
	+\left[ \partial_yw, \partial_x^m\right] v
	+[a, \partial_x^m]\partial_y^2 (w + w^s)
	\\
	&\quad
	-
	\frac{\partial_y(w + w^s + \varepsilon f(t)\chi'(y))}{w + w^s + \varepsilon f(t)\chi'(y)}
	\left(
	\left[ u, \partial_x^m\right] \partial_x u
	+\left[ w, \partial_x^m\right] v
	+[a, \partial_x^m]\partial_y^2 (u + u^s)
	\right)
	\\
	&=
	\left(w + w^s + \varepsilon f(t)\chi'(y)\right)
	\partial_y\left(
	\frac{
		\left[ u, \partial_x^m\right] \partial_x u
		+\left[ w, \partial_x^m\right] v
		+[a, \partial_x^m]\partial_y^2 (u + u^s)
	}{
		w + w^s + \varepsilon f(t)\chi'(y)
	}
	\right),
\end{align*}
we can rewrite \eqref{A2} as
\begin{align}
	\begin{split}
		&\left(
		\partial_t \partial_x^m w
		-
		\frac{\partial_y(w + w^s + \varepsilon f(t)\chi'(y))}{w + w^s + \varepsilon f(t)\chi'(y)}
		\partial_t \partial_x^m u
		\right)
		\\
		&\quad
		+\left(u + u^s + \varepsilon f(t)\chi(y)\right)
		\left(
		\partial_x^{m+1}w
		-
		\frac{\partial_y(w + w^s + \varepsilon f(t)\chi'(y))}{w + w^s + \varepsilon f(t)\chi'(y)}
		\partial_x^{m+1} u
		\right)
		\\
		&\quad
		-(a+1)
		\left(
		\partial_y^2 \partial_x^m w
		-
		\frac{\partial_y(w + w^s + \varepsilon f(t)\chi'(y))}{w + w^s + \varepsilon f(t)\chi'(y)}
		\partial_y^2\partial_x^m u
		\right)
		\\
		&=
		\left(w + w^s + \varepsilon f(t)\chi'(y)\right)
		\partial_y\left(
		\frac{
			\left[ u, \partial_x^m\right] \partial_x u
			+\left[ w, \partial_x^m\right] v
			+[a, \partial_x^m]\partial_y^2 (u + u^s)
		}{
			w + w^s + \varepsilon f(t)\chi'(y)
		}
		\right)
		\\
		&\quad
		+\partial_x^m a\,\varepsilon f(t)
		\left(
		\chi'''(y)
		-
		\frac{\partial_y(w + w^s + \varepsilon f(t)\chi'(y))}{w + w^s + \varepsilon f(t)\chi'(y)}
		\chi''(y)
		\right)
		+\partial_x^m F_3 .
	\end{split}
	\label{A3}
\end{align}

We next rewrite the terms on the left-hand side of \eqref{A3}.
For the first term, we have
\begin{align*}
	&\partial_t \partial_x^m w
	-
	\frac{\partial_y(w + w^s + \varepsilon f(t)\chi'(y))}{w + w^s + \varepsilon f(t)\chi'(y)}
	\partial_t \partial_x^m u
	\\
	&=
	\partial_t g_m
	+\partial_x^m u\,
	\partial_y\left(
	\frac{\partial_t (w + w^s + \varepsilon f(t)\chi'(y))}{w + w^s + \varepsilon f(t)\chi'(y)}
	\right).
\end{align*}
For the second term, since $w^s$ and $\chi$ are independent of $x$, we have
\begin{align*}
	&\left(u + u^s + \varepsilon f(t)\chi(y)\right)
	\left(
	\partial_x^{m+1}w
	-
	\frac{\partial_y(w + w^s + \varepsilon f(t)\chi'(y))}{w + w^s + \varepsilon f(t)\chi'(y)}
	\partial_x^{m+1} u
	\right)
	\\
	&=
	\left(u + u^s + \varepsilon f(t)\chi(y)\right)\partial_x g_m
	\\
	&\quad
	+\left(u + u^s + \varepsilon f(t)\chi(y)\right)\partial_x^m u\,
	\partial_y\left(
	\frac{\partial_x w}{w + w^s + \varepsilon f(t)\chi'(y)}
	\right).
\end{align*}
For the third term, a direct computation gives
\begin{align*}
	&-(a+1)
	\left(
	\partial_y^2 \partial_x^m w
	-
	\frac{\partial_y(w + w^s + \varepsilon f(t)\chi'(y))}{w + w^s + \varepsilon f(t)\chi'(y)}
	\partial_y^2\partial_x^m u
	\right)
	\\
	&=
	-(a+1)\partial_y^2 g_m
	-2(a+1)g_m\,
	\partial_y\left(
	\frac{\partial_y(w + w^s + \varepsilon f(t)\chi'(y))}{w + w^s + \varepsilon f(t)\chi'(y)}
	\right)
	\\
	&\quad
	-(a+1)\partial_x^m u
	\left(
	\frac{\partial_y^3 (w + w^s + \varepsilon f(t)\chi'(y))}{w + w^s + \varepsilon f(t)\chi'(y)}
	-
	\frac{\partial_y^2 (w + w^s + \varepsilon f(t)\chi'(y))}{w + w^s + \varepsilon f(t)\chi'(y)}
	\frac{\partial_y (w + w^s + \varepsilon f(t)\chi'(y))}{w + w^s + \varepsilon f(t)\chi'(y)}
	\right).
\end{align*}

Combining the above identities, the left-hand side of \eqref{A3} becomes
\begin{align*}
	&\partial_t g_m
	+\left(u + u^s + \varepsilon f(t)\chi(y)\right)\partial_x g_m
	-(a+1)\partial_y^2 g_m
	\\
	&\quad
	-2(a+1)g_m\,
	\partial_y\left(
	\frac{\partial_y(w + w^s + \varepsilon f(t)\chi'(y))}{w + w^s + \varepsilon f(t)\chi'(y)}
	\right)
	\\
	&\quad
	+\partial_x^m u
	\left[
	\partial_y\left(
	\frac{
		\left(\partial_t+\left(u + u^s + \varepsilon f(t)\chi(y)\right)\partial_x-(a+1)\partial_y^2\right)
		(w + w^s + \varepsilon f(t)\chi'(y))
	}{
		w + w^s + \varepsilon f(t)\chi'(y)
	}
	\right)
	-\partial_x w
	\right].
\end{align*}

On the other hand, from $\eqref{wprandtl}_1$ and the vorticity form of the shear flow equation
\begin{align}\label{vshear flow}
		\partial_t w^s - \partial_y^2 w^s = \varepsilon \partial_y m,
\end{align}
we obtain
\begin{align*}
	&\left(\partial_t+\left(u + u^s + \varepsilon f(t)\chi(y)\right)\partial_x-(a+1)\partial_y^2\right)
	\left(w + w^s + \varepsilon f(t)\chi'(y)\right)
	\\
	&=
	-v\,\partial_y\left(w + w^s + \varepsilon f(t)\chi'(y)\right)+F_3.
\end{align*}
Therefore,
\begin{align*}
	&\partial_y\left(
	\frac{
		\left(\partial_t+\left(u + u^s + \varepsilon f(t)\chi(y)\right)\partial_x-(a+1)\partial_y^2\right)
		(w + w^s + \varepsilon f(t)\chi'(y))
	}{
		w + w^s + \varepsilon f(t)\chi'(y)
	}
	\right)
	-\partial_x w
	\\
	&=
	-\partial_x w
	-\partial_y\left(
	\frac{\partial_y(w + w^s + \varepsilon f(t)\chi'(y))}{w + w^s + \varepsilon f(t)\chi'(y)}\,v
	\right)
	+\partial_y\left(
	\frac{F_3}{w + w^s + \varepsilon f(t)\chi'(y)}
	\right).
\end{align*}
Substituting this identity into \eqref{A3}, we conclude that
\begin{align*}
	&\partial_t g_m
	+\left(u + u^s + \varepsilon f(t)\chi(y)\right)\partial_x g_m
	-(a+1)\partial_y^2 g_m
	\\
	&\quad
	-2(a+1)g_m
	\partial_y\left(
	\frac{\partial_y(w + w^s + \varepsilon f(t)\chi'(y))}{w + w^s + \varepsilon f(t)\chi'(y)}
	\right)
	\\
	&=
	\partial_x^m u
	\left(
	\partial_x w
	+\partial_y\left(
	\frac{\partial_y(w + w^s + \varepsilon f(t)\chi'(y))}{w + w^s + \varepsilon f(t)\chi'(y)}\,v
	\right)
	-\partial_y\left(
	\frac{F_3}{w + w^s + \varepsilon f(t)\chi'(y)}
	\right)
	\right)
	\\
	&\quad
	+\left(w + w^s + \varepsilon f(t)\chi'(y)\right)\partial_y\left(
	\frac{
		\left[ u, \partial_x^m\right] \partial_x u
		+\left[ w, \partial_x^m\right] v
		+[a, \partial_x^m]\partial_y^2 (u + u^s)
	}{
		w + w^s + \varepsilon f(t)\chi'(y)
	}
	\right)
	\\
	&\quad
	+\partial_x^m a\,\varepsilon f(t)
	\left(
	\chi'''(y)
	-
	\frac{\partial_y(w + w^s + \varepsilon f(t)\chi'(y))}{w + w^s + \varepsilon f(t)\chi'(y)}
	\chi''(y)
	\right)
	+\partial_x^m F_3.
\end{align*}
This is exactly the evolution equation satisfied by $g_m$.

It remains to determine the boundary condition satisfied by $g_m$.
Since
$$
g_m=\partial_x^m w- \frac{\partial_y(w + w^s + \varepsilon f(t)\chi'(y))}{w + w^s + \varepsilon f(t)\chi'(y)}\partial_{x}^{m}u,
$$
and $\partial_x^m u|_{y=0}=0$, we have
$$
g_m\big|_{y=0}=\partial_x^m w\big|_{y=0}.
$$
Moreover,
\begin{align*}
	\partial_y g_m
	&=
	\partial_x^m \partial_y w
	- \partial_y\left( \frac{\partial_y(w + w^s + \varepsilon f(t)\chi'(y))}{w + w^s + \varepsilon f(t)\chi'(y)}\right) \partial_{x}^{m}u
	\\
	&\quad
	-  \frac{\partial_y(w + w^s + \varepsilon f(t)\chi'(y))}{w + w^s + \varepsilon f(t)\chi'(y)} \partial_{x}^{m}w.
\end{align*}
Hence, using again $\partial_x^m u|_{y=0}=0$, we obtain
$$
\partial_x^m\partial_y w\big|_{y=0}
=
\left(
\partial_y g_m
+
\frac{\partial_y(w + w^s + \varepsilon f(t)\chi'(y))}{w + w^s + \varepsilon f(t)\chi'(y)}\, g_m
\right)\bigg|_{y=0}.
$$

On the other hand, differentiating the boundary condition
$$
\bigl((a+1)\partial_y w + a\partial_y w^s + a\varepsilon f(t)\chi''(y)\bigr)\big|_{y=0}=0
$$
with respect to $x$ yields
\begin{align*}
	&\bigl((a+1)\partial_y\partial_x^m w\bigr)\big|_{y=0}
	+\bigl([\partial_x^m,a+1]\partial_y w\bigr)\big|_{y=0}
	+\bigl(\partial_x^m a\,\partial_y w^s\bigr)\big|_{y=0}
	\\
	&\qquad
	+\bigl(\partial_x^m a\,\varepsilon f(t)\chi''(y)\bigr)\big|_{y=0}
	=0.
\end{align*}
Since $\chi''(0)=0$, this reduces to
\begin{align*}
	\bigl((a+1)\partial_y\partial_x^m w\bigr)\big|_{y=0}
	+\bigl([\partial_x^m,a+1]\partial_y w\bigr)\big|_{y=0}
	+\bigl(\partial_x^m a\,\partial_y w^s\bigr)\big|_{y=0}
	=0.
\end{align*}
Therefore, $g_m$ satisfies the boundary condition
\begin{align*}
	&\Biggl[
	(a+1)\left(
	\partial_y g_m
	+
	\frac{\partial_y(w + w^s + \varepsilon f(t)\chi'(y))}{w + w^s + \varepsilon f(t)\chi'(y)}\, g_m
	\right)
	\\
	&\qquad
	+[\partial_x^m,a+1]\partial_y w
	+\partial_x^m a\,\partial_y w^s
	\Biggr]\bigg|_{y=0}
	=0.
\end{align*}
This completes the derivation.

\end{appendices}

\section*{Statement about conflicting interests}
The authors declare that there are no conflicts of interest.

\section*{Acknowledgements}
The authors would like to thank the referee for the valuable suggestions on improving this manuscript.

	\phantomsection

\end{document}